\documentclass[preprint,12pt,number]{elsarticle}

\usepackage[a4paper,margin=25mm]{geometry}

\usepackage{graphicx}%
\usepackage{microtype}
\usepackage{multirow}%
\usepackage{amsmath,amssymb,amsfonts}%
\usepackage{mathrsfs}%
\usepackage[title]{appendix}%
\usepackage{xcolor}%
\usepackage{textcomp}%
\usepackage{manyfoot}%
\usepackage{booktabs}%
\usepackage{algorithm}%
\usepackage{algorithmicx}%
\usepackage{algpseudocode}%
\usepackage{listings}%

\usepackage{stmaryrd}%
\usepackage{subcaption}%
\usepackage{placeins}%
\usepackage{array}%
\usepackage{longtable}%
\usepackage{tikz}%
\usepackage{tikz-cd}%
\usetikzlibrary{shapes,arrows.meta,positioning,calc,shadows,fit,backgrounds,decorations.pathmorphing}%

\DeclareGraphicsExtensions{.pdf,.png,.jpg,.eps}%

\algrenewcommand\algorithmicrequire{\textbf{Input:}}
\algrenewcommand\algorithmicensure{\textbf{Output:}}
\algrenewcommand\algorithmicreturn{\textbf{Return}}
\floatname{algorithm}{Algorithm}

\newcommand{\bv}{\mathbf{v}}

\newcommand{\bx}{\mathbf{x}}

\newtheorem{theorem}{Theorem}
\newtheorem{proposition}[theorem]{Proposition}

\newdefinition{assumption}{Assumption}
\newdefinition{remark}{Remark}
\newproof{proof}{Proof}
\AddToHook{env/proof/end}{\qed}

\biboptions{sort&compress}
\usepackage[colorlinks=true,citecolor=blue,linkcolor=blue,urlcolor=blue]{hyperref}
\makeatletter
\providecommand{\theHALG@line}{}
\renewcommand{\theHALG@line}{\thealgorithm.\arabic{ALG@line}}
\makeatother

\usepackage{cleveref}%
\crefname{assumption}{Assumption}{Assumptions}%
\Crefname{assumption}{Assumption}{Assumptions}%

\begin{document}

\begin{frontmatter}
  \title{Learning Coarse Spaces for Two-Grid Methods from Smoothed Errors  with Applications to High-Contrast Porous Media Flow}
  \author[thu,bimsa]{Chenglong Bao}
  \ead{clbao@tsinghua.edu.cn}
  \author[thu]{Chen Cui}
  \ead{chencui@tsinghua.edu.cn}
  \author[uestc]{Ruohua Li}
  \ead{2024110902017@std.uestc.edu.cn}
  \affiliation[thu]{organization={Yau Mathematical Sciences Center, Tsinghua University},
  city={Beijing}, postcode={100084}, country={China}}
  \affiliation[bimsa]{organization={Yanqi Lake Beijing Institute of Mathematical Sciences and Applications},
  city={Beijing}, postcode={101408}, country={China}}
  \affiliation[uestc]{organization={School of Mathematical Sciences, University of Electronic Science and Technology of China},
  city={Chengdu}, postcode={611731}, country={China}}

  \begin{abstract}
  Efficient two-grid solvers for high-contrast parametric elliptic problems require coarse spaces that capture error components poorly reduced by local smoothing. As these components vary with the problem parameters, repeatedly constructing parameter-adaptive coarse spaces can be computationally expensive. We propose a coarse-space learning method based on finite-step errors generated by additive Schwarz smoothing. For each training parameter, proper orthogonal decomposition (POD) of smoothed-error snapshots in the discrete energy inner product defines a reference coarse space. A continuous-coordinate neural generator is then trained using a basis-invariant, whitened subspace loss, enabling the prediction of parameter-adaptive coarse spaces on unseen grids. We establish discrete perturbation estimates relating coarse-space discrepancies to the associated nonsymmetric two-grid operators and derive conditional finite-step Krylov residual-polynomial bounds without assuming stationary energy contraction. Numerical experiments on high-contrast interface problems, moving-interface Darcy flow, and nonlinear Forchheimer flow demonstrate effective solver performance across varying parameters and grid resolutions. On the finest static test grid, with approximately 264 times as many degrees of freedom as the training grid, the learned spaces achieve FGMRES iteration counts close to those of parameter-specific POD references and comparable to GenEO spaces, while substantially reducing repeated coarse-space construction costs.
  \end{abstract}

\end{frontmatter}

\section{Introduction}
\label{sec:introduction}

High-contrast parametric elliptic problems arise in applications such as porous-media flow with heterogeneous coefficients, moving-interface problems, and flows with nonlinear constitutive laws. In many-query and time-dependent settings, the associated discrete systems must be solved repeatedly, making efficient iterative solvers particularly important. Such coefficient heterogeneity can give rise to error components that are poorly reduced by local relaxation methods. Multigrid and two-level domain decomposition methods (DDM) address this difficulty through coarse correction, whose effectiveness depends critically on the choice of coarse space \cite{hackbusch1985multi,trottenberg2000multigrid,toselli2005domain,xu1992iterative}. As the coefficient distribution varies with the problem parameters, the error components requiring coarse correction also change, making it difficult for a fixed coarse space to remain effective across the parameter range.

Several adaptive coarse-space techniques have been developed to address coefficient heterogeneity and multiscale features. GenEO (generalized eigenproblems in the overlap) employs local generalized eigenvalue problems to identify modes essential for robust convergence \cite{spillane2014abstract}. Related two-level restricted additive Schwarz methods use optimal local approximation spaces constructed by multiscale spectral generalized finite element methods \cite{strehlow2026fast}. Generalized multiscale finite element methods (GMsFEM) construct approximation spaces through local spectral decompositions to capture fine-scale heterogeneity \cite{efendiev2013generalized,fu2024efficient}. Although these approaches provide effective coefficient-adaptive coarse spaces, their construction involves solving eigenvalue problems or carrying out other fine-scale computations. When the coefficient field varies repeatedly, the associated reconstruction costs can become substantial. This motivates the development of methods that predict parameter-adaptive coarse spaces without repeatedly performing expensive spectral computations.

Neural iterative methods seek to learn components of classical iterative algorithms while retaining established discretizations \cite{mao2026accurate,moore2025graph,zhang2024blending,meng2026overview}. Existing approaches have considered learned multigrid smoothers \cite{huang2022learning,chen2022meta,cui2025neural}, restriction and prolongation operators and coarsening strategies \cite{katrutsa2020black,greenfeld2019learning,luz2020learning,caldana2024deep,taghibakhshi2021optimization}, coarse-grid operators \cite{huang2024reducing}, transmission conditions and coarse corrections for DDMs \cite{taghibakhshi2022learning,taghibakhshi2023mg}, neural preconditioners for Krylov methods \cite{dimola2026numerical,chen2025graph,trifonov2025efficient}, and hybrid iterative solvers \cite{cui2022fourier,zhang2024blending,cui2024hybrid,hu2025hybrid,sun2026learning}. These studies demonstrate the potential of learning solver components that can be reused across related problems.
More recently, several studies have explored learning parameter-adaptive coarse
spaces. Fanaskov et al. \cite{fanaskov2026deep} consider slowly decaying
eigenspaces of the stationary error-propagation operator, while Liu et al.
\cite{liu2026approximating} learn basis functions associated with established
multiscale constructions, such as GMsFEM. These studies demonstrate the feasibility of learning coefficient-dependent coarse spaces, although their training targets remain tied to particular spectral or multiscale constructions. The effectiveness of coarse correction is strongly influenced by how well the coarse space captures the errors remaining after smoothing, rather than on the particular construction used to obtain it. This suggests a more direct approach based on approximating smoothed errors.

Reduced-order modeling (ROM) offers a natural approach to constructing such approximations. Classical projection-based ROM methods, including proper orthogonal decomposition (POD) and reduced-basis methods, construct low-dimensional approximations of parametric solution families \cite{quarteroni2016reduced,hesthaven2016certified}, while nonlinear approaches seek more flexible representations of solution manifolds \cite{fresca2022pod,CSIAM-AM-4-13}. This perspective can also be applied to errors remaining after a prescribed number of smoothing steps. Since these errors depend on both the discrete PDE operator and the smoother, the resulting reduced spaces are directly related to the error components requiring coarse correction.

Motivated by this perspective, we develop a parameter-adaptive coarse-space learning method based on finite-step smoothed-error ensembles generated by the additive Schwarz method (ASM). Rather than approximating selected eigenvectors of an iteration operator, we construct reference coarse spaces through POD in the discrete energy inner product. This approach captures finite-step error dynamics without requiring symmetry or normality of the one-level additive Schwarz error-propagation operator. We formulate coarse-space prediction as a Grassmann-valued learning problem using a basis-invariant, whitened subspace loss that balances the retained POD directions. A continuous-coordinate neural representation further enables evaluation on grids different from those used for training. We further establish a discrete perturbation analysis for nonsymmetric two-grid preconditioners, relating coarse-space approximation errors to perturbations of the corresponding two-grid operators. An operator-weighted discrepancy characterizes the effects of pre- and postsmoothing through low-dimensional Gram matrices, while conditional finite-step Krylov polynomial estimates connect these perturbations to preconditioned iterative solvers without assuming stationary energy contraction.

We test the performance of the proposed method on three representative problems: static interface equations, moving-interface Darcy flow, and nonlinear Forchheimer flow. The numerical experiments assess the effectiveness of POD whitening, the benefits of parameter adaptation, and the robustness of the learned coarse spaces across grid resolutions and changing operators. The learned spaces achieve Krylov iteration counts close to parameter-specific POD references and competitive with local spectral coarse spaces, while requiring substantially fewer coarse basis functions. The results demonstrate the effectiveness of the proposed method for high-contrast problems, including time-dependent and nonlinear applications.

The remainder of the paper is organized as follows. \Cref{sec:preliminaries} introduces the model problems, finite volume discretization, additive Schwarz smoother, and two-grid method. \Cref{sec:methods} presents the construction and learning of parameter-adaptive coarse spaces. \Cref{sec:num-experiments} reports the numerical experiments. \Cref{sec:convergence} develops the discrete subspace and Krylov perturbation analysis for nonsymmetric two-grid preconditioners. Finally, \Cref{sec:conclusion} concludes the paper.

\paragraph{Notation}
Bold uppercase and lowercase symbols denote matrices and vectors,
respectively. Scalars, individual matrix or vector entries, and scalar-valued
functions are written in ordinary italic type. Calligraphic $\mathcal U$ and
$\mathcal V$ denote subspaces, whereas $\mathbf U$ and $\mathbf V$ denote
basis matrices representing them:
\[
  \mathcal U=\operatorname{range}(\mathbf U),\qquad
  \mathcal V=\operatorname{range}(\mathbf V).
\]
Different full-column-rank basis matrices may represent the same subspace.
In particular, $\mathcal U_{\mathrm{ref},h}^{r}(\mu)$ denotes an
$r$-dimensional POD reference space and
$\mathcal V_{\theta,h}^{n}(\mu)$ an $n$-dimensional predicted coarse space.
Here $h$ is the grid spacing, $N_h$ the fine-grid dimension, and $n$ the
prescribed coarse dimension. The superscript $T$ denotes transpose,
$\mathbf I$ the identity matrix of the required size, and
$\|\cdot\|_2$ and $\|\cdot\|_F$ the Euclidean (or spectral) and Frobenius
norms, respectively. For a symmetric positive definite (SPD) matrix $\mathbf A$, the energy norm is
$\|\mathbf x\|_{\mathbf A}=(\mathbf x^T\mathbf A\mathbf x)^{1/2}$.

\section{Preliminaries}
\label{sec:preliminaries}
This section introduces the model problem, the finite volume method (FVM),
the ASM smoother, and the two-grid (TG) method used in coarse-space construction,
learning, and analysis.

\subsection{Model problem}
\label{sec:model-problem}

\begin{figure}[!htbp]
  \centering
  \begin{subfigure}[b]{0.45\textwidth}
    \centering
    \begin{tikzpicture}[scale=0.8]
      \draw[thick] (0,0) rectangle (5,5);
      \draw[step=0.5cm, gray!20, thin] (0,0) grid (5,5);
      \draw[
        ultra thick,
        blue!70!black,
        fill=blue!10,
        fill opacity=0.5,
        smooth,
        samples=240,
        domain=90:450,
        variable=\t
      ]
      plot ({2.5 + (1.04 + 0.39*cos(5*(\t-90))) * cos(\t)},
      {2.55 + (1.04 + 0.39*cos(5*(\t-90))) * sin(\t)}) -- cycle;
      \node at (2.50, 2.80) {\large $\Omega_1$};
      \node at (4.30, 4.30) {\large $\Omega_2$};
      \node[blue!70!black] at (4.05, 3.30) {$\Gamma$};
      \foreach \x in {0.25, 0.75, ..., 4.75}
      \foreach \y in {0.25, 0.75, ..., 4.75}
      \fill[black!60] (\x,\y) circle (0.8pt);
      \draw[red, thick] (0.5, 0.5) rectangle (1.0, 1.0);
      \fill[red] (0.75, 0.75) circle (1.5pt);
    \end{tikzpicture}
    \caption{Computational domain and cell-centered finite-volume grid}
    \label{fig:domain}
  \end{subfigure}
  \hfill
  \begin{subfigure}[b]{0.45\textwidth}
    \centering
    \begin{tikzpicture}[scale=0.8]
      \draw[thick] (0,0) rectangle (5,5);
      \fill[gray!20] (1.3, 0) rectangle (1.9, 5);
      \fill[gray!20] (3.1, 0) rectangle (3.7, 5);
      \fill[gray!20] (0, 1.3) rectangle (5, 1.9);
      \fill[gray!20] (0, 3.1) rectangle (5, 3.7);
      \foreach \x in {1.3, 1.9, 3.1, 3.7}
      \draw[thin, dashed, gray!60] (\x, 0) -- (\x, 5);
      \foreach \y in {1.3, 1.9, 3.1, 3.7}
      \draw[thin, dashed, gray!60] (0, \y) -- (5, \y);
      \node at (0.65, 4.35) {\footnotesize $D_1$};
      \node at (2.50, 2.50) {\footnotesize $D_5$};
      \node at (4.35, 0.65) {\footnotesize $D_9$};
      \draw[<->, thick] (1.3, 4.2) -- (1.9, 4.2);
      \node[above, scale=0.7] at (1.6, 4.2) {$\delta$};
      \draw[|-|, thick, blue] (0.15, 0.35) -- (1.75, 0.35);
      \node[above, blue, scale=0.8] at (0.95, 0.35) {$w$};
      \draw[|-|, thick, red!70!black] (0.15, 0.75) -- (1.25, 0.75);
      \node[above, red!70!black, scale=0.8] at (0.70, 0.75) {$d$};
    \end{tikzpicture}
    \caption{Overlapping additive Schwarz decomposition}
    \label{fig:ddm}
  \end{subfigure}

  \caption{Computational domain and domain decomposition for the interface
    problem. Panel~(a) shows the interior and exterior regions $\Omega_1$ and
    $\Omega_2$, separated by the interface $\Gamma$, together with the
    cell-centered FVM grid. Panel~(b) shows the overlapping decomposition, with the
    shaded bands indicating the overlap regions. The decomposition is
    determined by the subdomain side length $w$, stride $d$, and overlap width
  $\delta$}
  \label{fig:preliminaries}
\end{figure}
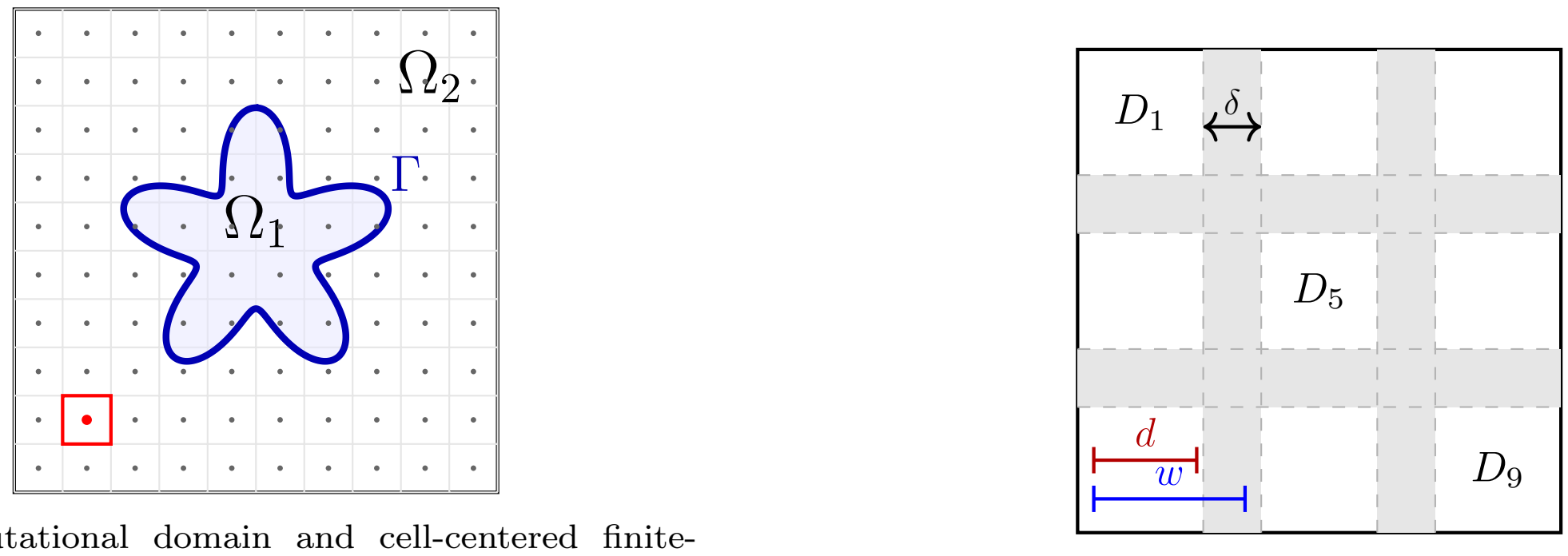

Consider the following two-dimensional elliptic pressure equation with a
discontinuous coefficient
\begin{equation}
  \begin{cases}
    -\nabla\cdot\bigl(\lambda(\boldsymbol{x})\nabla p\bigr)=f,
    & \text{in }\Omega,\\
    p=0, & \text{on }\partial\Omega,\\
    \llbracket p \rrbracket=0,\quad
    \llbracket \lambda\nabla p\cdot\mathbf n\rrbracket=0,
    & \text{on }\Gamma,
  \end{cases}
  \label{eq:jumping}
\end{equation}
where $\Omega\subset\mathbb R^2$ is a bounded computational domain and the
physical interface $\Gamma\subset\Omega$ partitions $\Omega$ into an interior
subdomain $\Omega_1$ and an exterior subdomain $\Omega_2$, as shown in
\Cref{fig:domain}. The notation $\llbracket\,\cdot\,\rrbracket$ denotes the jump
across $\Gamma$, and $\mathbf n$ is the unit normal on $\Gamma$. Both the
pressure $p$ and the normal flux $\lambda\nabla p\cdot\mathbf n$ are continuous
across the interface.
We consider the piecewise constant diffusion coefficient
\begin{equation}
  \lambda(\boldsymbol{x}) =
  \begin{cases}
    10^m, & \boldsymbol{x} \in \Omega_1,\\
    1,    & \boldsymbol{x} \in \Omega_2,
  \end{cases}
  \qquad m \in \mathbb{N}_0.
\end{equation}
For $m>0$, the coefficient jumps across the interface and the pressure
gradient is generally discontinuous there.

We represent the interface implicitly by a level set function
\cite{osher1988fronts}. Let $\phi:\Omega\to\mathbb R$ be continuous and satisfy
\[
  \Gamma=\{\boldsymbol{x}\in\Omega:\phi(\boldsymbol{x})=0\},\qquad
  \phi<0\ \text{in }\Omega_1,\qquad
  \phi>0\ \text{in }\Omega_2.
\]
We discretize \eqref{eq:jumping} using an interface-aware cell-centered FVM on
a uniform Cartesian grid. For a control-volume face intersected by the
physical interface, the intersection location is estimated by linear
interpolation of the level set function, and the transmissibility is computed
as a series resistance weighted by the actual distances on either side of the
interface.
Standard harmonic transmissibilities are used on all other interior faces,
and a half-cell-distance correction is used on Dirichlet boundaries. The
scheme retains the standard five-point stencil and produces the
parameterized symmetric positive definite system
\begin{equation}
  \mathbf{A}_h(\mu)\,\mathbf{p}_h(\mu)=\mathbf{f}_h,
  \qquad \mu\in\mathcal P,
  \label{eq:discrete_system}
\end{equation}
where $\mathbf{A}_h(\mu)\in\mathbb R^{N_h\times N_h}$ and
$N_h=N_{\mathrm{side}}^2$, with $N_{\mathrm{side}}$ denoting the number of
cells along each side of the domain. The complete five-point discretization,
boundary treatment, and a manufactured-solution convergence study are
provided in Appendix~\ref{sec:interface-aware-ccfvm}.
Here $\mu$ denotes the problem parameters that determine the discrete
operator. For example, $\mu_{\rm st}=(\lambda,\phi)$ for this interface
problem.

To assess the proposed method in many-query settings with changing parameters,
we also consider two applications that require sequences of linear diffusion
solves. The first is a moving-interface Darcy problem. Its interface
$\Gamma(t)$ satisfies the transport equation
\cite{sussman1994levelset}
\[
  \frac{\partial \phi}{\partial t} + \mathbf{v} \cdot \nabla \phi = 0,
\]
where $\mathbf v=-\lambda\nabla p$ is the Darcy velocity. At each time step,
\eqref{eq:discrete_system} is assembled and solved using the current level set
function, after which the resulting velocity advances the interface. The
discrete linear operator therefore changes in time.

The second application is nonlinear Forchheimer flow
\cite{whitaker1996forchheimer,girault2008darcyforchheimer}.
The pressure $p$ and flux $\mathbf q$ satisfy
\begin{equation}
  \begin{cases}
    -\nabla p=\dfrac{\mu_f}{\kappa(\boldsymbol{x},p)}\,\mathbf{q}
    +\beta\rho |\mathbf{q}|\,\mathbf{q}, & \text{in }\Omega,\\
    \nabla\cdot\mathbf{q}=f, & \text{in }\Omega,\\
    p=0, & \text{on }\partial\Omega,
  \end{cases}
  \label{eq:forchheimer-model}
\end{equation}
where $\mu_f$ is the fluid viscosity, $\rho$ is the fluid density, $\beta$ is
the Forchheimer coefficient, and $\kappa(\boldsymbol x,p)$ is the
permeability. Introducing the effective mobility
\[
  \mathbf{q}=-\lambda(\boldsymbol{x},p,|\nabla p|)\nabla p
\]
rewrites \eqref{eq:forchheimer-model} as
\begin{equation}
  -\nabla\cdot\Bigl(\lambda(\boldsymbol{x},p,|\nabla p|)\nabla p\Bigr)=f,
  \qquad
  \lambda=\frac{2}{\dfrac{\mu_f}{\kappa(\boldsymbol{x},p)}
    +\sqrt{\left(\dfrac{\mu_f}{\kappa(\boldsymbol{x},p)}\right)^2
  +4\beta\rho|\nabla p|}}.
  \label{eq:forchheimer-nonlinear-diffusion}
\end{equation}
Let $\kappa_0(\boldsymbol{x})=\kappa(\boldsymbol{x},p_0)>0$ denote the
prescribed spatial permeability field at the reference pressure $p_0=0$.
We consider both pressure-independent permeability
$\kappa(\boldsymbol{x},p)=\kappa_0(\boldsymbol{x})$ and pressure-dependent
permeability
$\kappa(\boldsymbol{x},p)=\kappa_0(\boldsymbol{x})
\exp\!\bigl(\alpha(p-p_0)\bigr)$, where $\alpha$ is the pressure-sensitivity
coefficient. Using
Picard iteration to linearize the nonlinear operator yields a sequence of
linearized systems of the form \eqref{eq:discrete_system}, with the operator
changing from one iteration to the next. At Picard iteration $r$, the
effective mobility is frozen at the current iterate
\[
  \lambda^{[r]}(\boldsymbol{x})
  =\lambda\bigl(\boldsymbol{x},p^{[r]},|\nabla p^{[r]}|\bigr),
\]
and the linearized system
\begin{equation}
  \mathbf{A}_h(\lambda^{[r]})\,\mathbf{p}_h^{[r+1]}=\mathbf{f}_h
  \label{eq:prelim-picard-system}
\end{equation}
is solved.

\subsection{Additive Schwarz method}
\label{sec:asm}

We use a multiplicity-weighted one-level overlapping ASM as the smoother. See \cite{smith1996domain} for the classical framework and \cite{cai1999restricted} for a related nonsymmetric restricted variant. After the fine grid is partitioned into the $K$ overlapping subdomains illustrated in \Cref{fig:ddm}, the additive Schwarz preconditioner and its error-propagation operator are
\begin{align}
  \mathbf{M}_{\mathrm{ASM}}^{-1}
  &=\mathbf{W}^{-1}\sum_{k=1}^{K}\mathbf{R}_k^T\mathbf{B}_k\mathbf{R}_k,
  \label{eq:asm-preconditioner}\\
  \mathbf{S}_h&=\mathbf{I}-\omega \mathbf{M}_{\mathrm{ASM}}^{-1}\mathbf{A}_h.
  \label{eq:error-propagation}
\end{align}
Here $\mathbf{R}_k$ is the restriction operator onto subdomain $k$, $\mathbf{R}_k^T$ is the corresponding zero-extension operator, and $\mathbf{B}_k$ denotes the local solver on that subdomain. The matrix $\mathbf{W}$ is the multiplicity matrix whose diagonal entries record the number of subdomains containing each global degree of freedom. We take $\omega=0.8$ in the numerical experiments. Details of the local solvers, restriction and prolongation operators, artificial subdomain boundaries, partition of unity, padding, and tensorized implementation are provided in Appendix~\ref{sec:asm-implementation}.

\subsection{Two-grid method}
\label{sec:two-grid}
The ASM rapidly
attenuates local high-frequency errors through parallel subdomain corrections (see Appendix~\ref{sec:asm-smoothing-visualization}), but it is less effective at reducing global low-frequency errors, particularly the slowly decaying errors induced
by high-contrast channels, interface geometry, and the connectivity of the
medium. We therefore introduce a low-dimensional coarse space
$\mathcal V_h(\mu)$ and use a Galerkin coarse
correction to exchange global information. The resulting iteration is given in
\Cref{alg:two-grid}.
\begin{algorithm}[!htbp]
  \caption{Two-grid iterative method with additive Schwarz smoothing.}
  \label{alg:two-grid}
  \begin{algorithmic}[1]
    \Require Initial iterate $\mathbf{p}_h^{(0)}$, operator $\mathbf{A}_h(\mu)$,
    smoother $\mathbf{M}_{\mathrm{ASM}}(\mu)^{-1}$, damping $\omega$, coarse matrix
    $\mathbf V_h(\mu)$, and smoothing steps $\nu_1,\nu_2$
    \Ensure Approximate solution $\mathbf p_h^{(\ell)}$ to
    $\mathbf{A}_h(\mu)\mathbf p_h=\mathbf f_h$
    \For{$\ell = 1, 2, \ldots$ until convergence}
    \State $\mathbf{w} \gets \mathbf{p}_h^{(\ell-1)}$
    \For{$i = 1, \ldots, \nu_1$} \Comment{Presmoothing}
    \State $\mathbf{w} \gets \mathbf{w} + \omega \mathbf{M}_{\mathrm{ASM}}(\mu)^{-1}\bigl(\mathbf{f}_h - \mathbf{A}_h(\mu)\mathbf{w}\bigr)$
    \EndFor
    \State $\mathbf{r} \gets \mathbf{f}_h - \mathbf{A}_h(\mu)\mathbf{w}$ \Comment{Compute residual}
    \State Solve $\bigl(\mathbf V_h(\mu)^{\top} \mathbf{A}_h(\mu) \mathbf V_h(\mu)\bigr)\mathbf{e}_c = \mathbf V_h(\mu)^{\top}\mathbf{r}$ \Comment{Coarse problem}
    \State $\mathbf{w} \gets \mathbf{w} + \mathbf V_h(\mu)\mathbf{e}_c$ \Comment{Coarse correction}
    \For{$i = 1, \ldots, \nu_2$} \Comment{Postsmoothing}
    \State $\mathbf{w} \gets \mathbf{w} + \omega \mathbf{M}_{\mathrm{ASM}}(\mu)^{-1}\bigl(\mathbf{f}_h - \mathbf{A}_h(\mu)\mathbf{w}\bigr)$
    \EndFor
    \State $\mathbf{p}_h^{(\ell)} \gets \mathbf{w}$
    \EndFor
  \end{algorithmic}
\end{algorithm}

Constructing the coarse space generally requires prior information about the
problem, such as its geometry and coefficient distribution. Classical
methods may combine prescribed geometry or near-nullspace information with
interface constraints and energy-minimizing extensions, or construct
coefficient-adaptive spectral spaces from local generalized eigenproblems and
snapshot compression. Representative examples include energy-minimizing
coarse spaces, GenEO, and GMsFEM
\cite{dohrmann2008family,spillane2014abstract,efendiev2013generalized}. These
methods are often numerically effective, but repeatedly rebuilding a coarse
space for a changing operator can be costly in a many-query setting. We
therefore ask a solver-oriented question: can a coarse space adapted to the
prescribed solver be constructed from the errors that actually remain after
ASM smoothing steps, without specifying those difficult error components in
advance?

\section{Learning parameter-adaptive coarse spaces}
\label{sec:methods}

To this end, we develop a parameter-adaptive coarse-space learning method based on errors remaining after additive Schwarz smoothing. For each training parameter, we first generate a finite-step smoothed-error ensemble by applying ASM to random initial errors. We then compress this ensemble through POD in the discrete energy inner product to obtain a parameter-specific reference coarse space. The resulting family of subspaces defines a Grassmann-valued learning problem, which we address using a basis-invariant, whitened subspace loss. Finally, a continuous-coordinate neural generator represents the learned basis functions, allowing the predicted coarse spaces to be evaluated on grids different from those used for training.

\subsection{POD reference coarse spaces constructed from smoothed-error snapshots}
\label{sec:smoothed-error-target}

Let $\mathbf{A}_h(\mu)$ denote the discrete operator and $\mathbf{S}_h(\mu)$ the homogeneous
ASM error-propagation operator. For fixed $\mu$ and grid spacing $h$,
 draw
\[
  \mathbf e_i^{(0)}
  \sim \mathcal N(\mathbf 0,\mathbf{I}_{N_h}),
  \qquad i=1,\ldots,q,
\]
and collect the initial errors in
\[
  \mathbf{E}^{(0)}
  =
  [\mathbf e_1^{(0)},\ldots,\mathbf e_q^{(0)}].
\]
After a prescribed number $\nu_{\rm snap}$ of ASM steps, we obtain the
finite-step error ensemble
\begin{equation}
  \mathbf{E}_\mu
  :=
  \mathbf{S}_h(\mu)^{\nu_{\rm snap}}\mathbf{E}^{(0)}
  =
  [\mathbf e_1(\mu),\ldots,\mathbf e_q(\mu)].
  \label{eq:smoothed-error-set}
\end{equation}
The construction uses the homogeneous error iteration and is therefore
independent of the right-hand side of an online linear system. Components that
are efficiently reduced by ASM contribute less to $\mathbf{E}_\mu$, whereas components
that remain significant over the selected smoothing horizon are retained.
We use this finite-step ensemble to define a POD reference space
for coarse correction.

The ensemble in \eqref{eq:smoothed-error-set} captures finite-step error
dynamics rather than directly approximating an asymptotic invariant subspace
of $\mathbf{S}_h(\mu)$. This distinction is relevant because the generally nonsymmetric ASM
error-propagation operator need not be normal in the $\mathbf{A}_h(\mu)$ inner product, so its
finite-step behavior need not be characterized by eigenvalues alone.
The snapshot construction also
avoids solving a parameter-adaptive global eigenproblem.

Since the subsequent Galerkin coarse correction is naturally measured in the
discrete energy norm, we compress $\mathbf{E}_\mu$ by POD in the
$\mathbf{A}_h(\mu)$ inner product. The complete construction is summarized in
\Cref{alg:reference-space}. The nominal coarse dimension is fixed at $n$,
and $r_{\rm w}(\mu)\le n$ denotes the reliable retained rank.
Here, $\boldsymbol{\Xi}_{r_{\rm w}}$ contains the first $r_{\rm w}$ eigenvectors of
$\mathbf{E}^T\mathbf{A}\mathbf{E}$, and
$\boldsymbol{\Sigma}_{r_{\rm w}}=\operatorname{diag}(\sigma_1,\ldots,\sigma_{r_{\rm w}})$.
The construction assumes $\mathbf{E}\ne0$, so the retained singular values are
positive and $\boldsymbol{\Sigma}_{r_{\rm w}}^{-1}$ is well defined.

\begin{algorithm}[!htbp]
  \caption{Construction of a POD reference coarse space}
  \label{alg:reference-space}
  \begin{algorithmic}[1]
    \Require
    Discrete operator $\mathbf{A}=\mathbf{A}_h(\mu)$,
    smoothed-error ensemble $\mathbf{E}=\mathbf{E}_\mu$ from \eqref{eq:smoothed-error-set},
    nominal coarse dimension $1\le n\le N_h$,
    relative singular-value threshold $0<\varepsilon_{\rm svd}\le1$

    \State Solve the  weighted Gram-matrix eigenproblem
    \[
      \mathbf{E}^T\mathbf{A}\mathbf{E}
      =
      \boldsymbol{\Xi}\boldsymbol{\Sigma}^2\boldsymbol{\Xi}^T,
      \qquad
      \sigma_1\ge\sigma_2\ge\cdots\ge0
    \]

    \State Set the reliable retained rank
    \[
      r_{\rm w}
      \gets
      \min\left\{
        n,\,
        \#\left\{
          k:
          \sigma_k\ge
          \varepsilon_{\rm svd}\sigma_1
        \right\}
      \right\}
    \]

    \State Construct the reliable $\mathbf{A}$-orthonormal POD modes
    \[
      \boldsymbol{\Phi}_{\mu,r_{\rm w}}
      \gets
      \mathbf{E}\boldsymbol{\Xi}_{r_{\rm w}}\boldsymbol{\Sigma}_{r_{\rm w}}^{-1}
    \]

    \State Define
    \[
      \mathcal U_{\mathrm{ref},h}^{r_{\rm w}}(\mu)
      :=
      \operatorname{range}(\boldsymbol{\Phi}_{\mu,r_{\rm w}})
    \]

    \State \Return
    $\boldsymbol{\Phi}_{\mu,r_{\rm w}}$ and
    $\mathcal U_{\mathrm{ref},h}^{r_{\rm w}}(\mu)$
  \end{algorithmic}
\end{algorithm}

The use of the energy metric is consistent with the Galerkin correction.
Indeed, for $\mathbf{A}=\mathbf{A}_h(\mu)$ and any full-column-rank matrix $\mathbf{V}$,
\begin{equation}
  \mathbf{P}_{\mathbf{V}}^{(\mathbf{A})}
  :=
  \mathbf{V}(\mathbf{V}^T\mathbf{A}\mathbf{V})^{-1}\mathbf{V}^T\mathbf{A}
  \label{eq:energy-projector}
\end{equation}
is the $\mathbf{A}$-orthogonal projector onto $\operatorname{range}(\mathbf{V})$, and the
post-correction error is
\[
  \mathbf e_h^+
  =
  (\mathbf{I}-\mathbf{P}_{\mathbf{V}}^{(\mathbf{A})})\mathbf e_h.
\]
Thus, minimizing the mean squared energy norm of the snapshot errors
remaining after coarse correction leads directly to an $\mathbf{A}$-weighted
projection problem. For a prescribed coarse dimension, energy-POD solves
this problem optimally. The following proposition states this optimality
and its worst-case counterpart on the finite-snapshot ellipsoid.
We write $\|\mathbf{M}\|_{\mathbf{A},F}^2:=\operatorname{tr}(\mathbf{M}^T\mathbf{A}\mathbf{M})$ for the squared
$\mathbf{A}$-weighted Frobenius norm and $d_n(\cdot,\mathbf{A})$ for the Kolmogorov
$n$-width in the $\mathbf{A}$-energy norm.

\begin{proposition}[Energy-POD optimality]
  \label[proposition]{prop:weighted-pod-optimality}
  Let $\sigma_1\ge\sigma_2\ge\cdots$ be the singular values of
  $\mathbf{A}^{1/2}\mathbf{E}$, extended by zeros as needed. For
  $1\le n\le\operatorname{rank}(\mathbf{E})$, the leading $\mathbf{A}$-weighted POD space
  $\mathcal U_n=\operatorname{range}(\boldsymbol{\Phi}_n)$ minimizes the empirical
  mean-square projection error
  \begin{equation}
    \min_{\substack{\mathbf{V}\in\mathbb R^{N_h\times n}\\\operatorname{rank}(\mathbf{V})=n}}
    \frac1q
    \left\|
    (\mathbf{I}-\mathbf{P}_{\mathbf{V}}^{(\mathbf{A})})\mathbf{E}
    \right\|_{\mathbf{A},F}^2
    =
    \frac1q
    \sum_{j>n}\sigma_j^2.
    \label{eq:empirical-pod-objective}
  \end{equation}
  Moreover, for the finite-snapshot ellipsoid
  \[
    \mathfrak E_{h,q}
    :=
    \{\mathbf{E}\mathbf c:\mathbf c\in\mathbb R^q,\ \|\mathbf c\|_2\le1\},
  \]
  the same space minimizes the worst-case $\mathbf{A}$-energy approximation error, with
  \begin{equation}
    d_n(\mathfrak E_{h,q},\mathbf{A})
    =
    \sigma_{n+1}.
    \label{eq:discrete-snapshot-ellipsoid-nwidth}
  \end{equation}
\end{proposition}

\begin{proof}
  Let $\boldsymbol{\Pi}_{\mathbf{W}}$ denote the Euclidean orthogonal projector onto
  $\operatorname{range}(\mathbf{A}^{1/2}\mathbf{V})$. Then
  $\mathbf{A}^{1/2}(\mathbf{I}-\mathbf{P}_{\mathbf{V}}^{(\mathbf{A})})=(\mathbf{I}-\boldsymbol{\Pi}_{\mathbf{W}})\mathbf{A}^{1/2}$, so the objective in
  \eqref{eq:empirical-pod-objective} equals
  $q^{-1}\|(\mathbf{I}-\boldsymbol{\Pi}_{\mathbf{W}})\mathbf{A}^{1/2}\mathbf{E}\|_F^2$, and the worst-case error over
  $\mathfrak E_{h,q}$ equals $\|(\mathbf{I}-\boldsymbol{\Pi}_{\mathbf{W}})\mathbf{A}^{1/2}\mathbf{E}\|_2$. By the
  Eckart--Young theorem \cite{eckart1936approximation} and its spectral-norm
  counterpart \cite{pinkus1985nwidths}, both are minimized when
  $\operatorname{range}(\mathbf{A}^{1/2}\mathbf{V})$ is spanned by the leading $n$ left singular
  vectors of $\mathbf{A}^{1/2}\mathbf{E}$, with minimal values $q^{-1}\sum_{j>n}\sigma_j^2$
  and $\sigma_{n+1}$, respectively. Since
  $\mathbf{A}^{1/2}\boldsymbol{\Phi}_n=\mathbf{A}^{1/2}\mathbf{E}\boldsymbol{\Xi}_n\boldsymbol{\Sigma}_n^{-1}$ consists of exactly these
  singular vectors, the minimizer is $\mathcal U_n$.
\end{proof}

For the primary setting considered below, $r_{\rm w}(\mu)=n$ for all
training parameters, so that
\[
  \mathcal U_{\mathrm{ref},h}^{n}(\mu)
  =
  \operatorname{range}(\boldsymbol{\Phi}_{\mu,n})
\]
defines a fixed-dimensional family of POD reference spaces. The case
$r_{\rm w}(\mu)<n$ is handled by the one-sided subspace loss introduced in
\Cref{rem:rank-deficient-targets}, while the predicted coarse dimension
remains fixed at $n$.

\begin{remark}[Relation to adaptive multilevel methods]
  \label{rem:smoothed-error-motivation}
  The use of errors exposed by relaxation to guide coarse-space construction
  is an established principle in adaptive multilevel methods. Adaptive
  smoothed aggregation uses errors remaining after application of the current
  iteration, including those generated from random initial errors, to improve
  the coarse representation \cite{brezina2005adaptive}. Related constructions
  use ensembles of smooth or low-energy test vectors, including local SVD
  compression \cite{chow2006aggregation} and bootstrap AMG
  \cite{brandt2011bootstrap}. A similar solver-oriented viewpoint appears in
  Schwarz methods, where slowly converging eigenmodes identify continuous and
  discontinuous coarse-space components \cite{cuvelier2023fundamental}, with
  further analysis of the continuous Robin--Schwarz spectrum in
  \cite{bennequintechniques}. Dominant singular vectors of local interface
  maps also provide Helmholtz coarse spaces \cite{gander2025hierarchical}.

Our construction follows this principle but uses the error information differently. Rather than identifying selected eigenmodes or adapting a coarse space directly for each online operator, we apply a prescribed number of ASM smoothing steps to random initial errors and compress the resulting global error ensemble through POD in the discrete energy inner product. This produces a parameter-specific reference coarse space without requiring spectral analysis of the ASM error-propagation operator. As the problem parameters vary, these reference spaces form a family of subspaces, which provides the starting point for the coarse-space learning problem considered next.
\end{remark}

\subsection{Basis-invariant subspace learning with POD whitening}
\label{sec:loss-spec}

A Galerkin coarse correction is invariant under a change of basis within the
coarse space. Indeed, for any invertible matrix \(\mathbf{R}\in\mathbb R^{n\times n}\)
\[
  \operatorname{range}(\mathbf{V})
  =
  \operatorname{range}(\mathbf{V}\mathbf{R}),
  \qquad
  \mathbf{P}_{\mathbf{V}\mathbf{R}}^{(\mathbf{A})}
  =
  \mathbf{P}_{\mathbf{V}}^{(\mathbf{A})}.
\]
Directly regressing individual POD basis vectors would therefore introduce an artificial dependence on the choice of basis. In particular, changes in sign, ordering, or scaling do not alter the coarse space but can affect a columnwise regression loss. The natural learning target is instead the
parameter-to-Grassmannian map
\begin{equation}
  F_{\theta,h}:
  \mathcal P
  \longrightarrow
  \operatorname{Gr}(n,N_h),
  \qquad
  \mu
  \longmapsto
  \mathcal V_{\theta,h}^{n}(\mu).
  \label{eq:param-to-coarse-space-map}
\end{equation}

Suppose that the network outputs a full-column-rank raw basis matrix
\[
  \mathbf{V}_\theta
  :=
  \mathbf{V}_{\theta,h}^{\rm raw}(\mu)
  \in\mathbb R^{N_h\times n}.
\]
Its column space is
$\mathcal V_{\theta,h}^{n}(\mu)=\operatorname{range}(\mathbf V_\theta)$.
The POD reference space is represented by $\boldsymbol\Phi_{\mu,n}$ when
the retained rank equals $n$.
During training, we use an $\mathbf{A}$-orthonormal representative to improve
numerical conditioning
\begin{equation}
  \mathbf{G}_\theta
  =
  \mathbf{V}_\theta^T\mathbf{A}\mathbf{V}_\theta
  =
  \mathbf{L}_\theta \mathbf{L}_\theta^T,
  \qquad
  \mathbf{Q}_\theta
  =
  \mathbf{V}_\theta \mathbf{L}_\theta^{-T}.
  \label{eq:training-energy-orthogonalization}
\end{equation}
Then
\[
  \mathbf{Q}_\theta^T\mathbf{A}\mathbf{Q}_\theta=\mathbf{I},
  \qquad
  \mathbf{P}_\theta^{(\mathbf{A})}
  :=
  \mathbf{P}_{\mathbf{V}_\theta}^{(\mathbf{A})}
  =
  \mathbf{Q}_\theta \mathbf{Q}_\theta^T\mathbf{A}.
\]

Training directly on all raw snapshots would give the normalized snapshot projection
loss. Let $r_{\mathrm{snap}}(\mu):=\operatorname{rank}(\mathbf{E}_\mu)$, then
\begin{equation}
  \mathcal L_{\rm raw}
  =
  \frac{
    \|(\mathbf{I}-\mathbf{P}_\theta^{(\mathbf{A})})\mathbf{E}_\mu\|_{\mathbf{A},F}^{2}
  }{
    \|\mathbf{E}_\mu\|_{\mathbf{A},F}^{2}
  }
  =
  \frac{
    \displaystyle
    \sum_{j=1}^{r_{\mathrm{snap}}(\mu)}
    \sigma_{\mu,j}^{2}
    \|(\mathbf{I}-\mathbf{P}_\theta^{(\mathbf{A})})
    \boldsymbol\phi_{\mu,j}\|_{\mathbf{A}}^{2}
  }{
    \displaystyle
    \sum_{j=1}^{r_{\mathrm{snap}}(\mu)}\sigma_{\mu,j}^{2}
  }.
  \label{eq:spectral-weighted-residual}
\end{equation}
When the snapshot singular values decay rapidly, the directions associated
with the largest singular values dominate both the loss and its gradient.
Such weighting is natural for reconstructing the raw snapshot matrix, but it
does not precisely match our objective of learning the retained
$n$-dimensional subspace $\mathcal U_{\rm ref,h}^{n}(\mu)$.

We therefore remove the singular-value weighting among the retained
directions after POD truncation. Equivalently, we define the whitened target
\[
  \widetilde{\mathbf{E}}_{\mu,n}
  =
  \mathbf{E}_\mu
  \boldsymbol{\Xi}_{\mu,n}
  \boldsymbol{\Sigma}_{\mu,n}^{-1}
  =
  \boldsymbol{\Phi}_{\mu,n},
\]
and use the basis-invariant whitened subspace loss
\begin{equation}
  \begin{aligned}
    \mathcal L_{\mathrm{white}}(\theta,\mu)
    &:=
    \frac1n
    \left\|
    (\mathbf{I}-\mathbf{P}_\theta^{(\mathbf{A})})
    \boldsymbol{\Phi}_{\mu,n}
    \right\|_{\mathbf{A},F}^{2}
    \\
    &=
    1-
    \frac1n
    \left\|
    \mathbf{Q}_\theta^T\mathbf{A}\boldsymbol{\Phi}_{\mu,n}
    \right\|_F^{2}.
  \end{aligned}
  \label{eq:white-projection-loss}
\end{equation}
If $\vartheta_1,\ldots,\vartheta_n$ are the principal angles between the
predicted and POD reference spaces in the $\mathbf{A}$-energy inner product, then
\begin{equation}
  \mathcal L_{\mathrm{white}}
  =
  \frac1n
  \sum_{j=1}^{n}\sin^2\vartheta_j.
  \label{eq:principal-angle-loss}
\end{equation}
The loss therefore depends only on the relative position of the two subspaces
and is invariant under changes of basis within either space.

Thus, POD singular values determine which directions are retained, whereas
whitening assigns equal weight to the retained directions during subspace
regression. When both spaces have dimension \(n\),
\[
  \mathcal L_{\mathrm{white}}=0
  \quad\Longleftrightarrow\quad
  \mathcal V_{\theta,h}^{n}(\mu)
  =
  \mathcal U_{\mathrm{ref},h}^{n}(\mu).
\]

\begin{remark}[Rank-deficient targets]
  \label{rem:rank-deficient-targets}
  When \(1\le r_{\rm w}(\mu)<n\), the predicted space remains \(n\)-dimensional and
  is trained only to contain the reliable reference directions. Specifically,
  \[
    \mathcal L_{\mathrm{white}}^{(r_{\rm w})}
    =
    1-\frac{1}{r_{\rm w}}
    \|\mathbf{Q}_\theta^T\mathbf{A}\boldsymbol{\Phi}_{\mu,r_{\rm w}}\|_F^2,
  \]
  so that
  \[
    \mathcal L_{\mathrm{white}}^{(r_{\rm w})}=0
    \iff
    \mathcal U_{\rm ref,h}^{r_{\rm w}}(\mu)
    \subseteq
    \mathcal V_{\theta,h}^{n}(\mu).
  \]
  The corresponding one-sided approximation estimate and the limitation of
  comparing unequal-rank projectors are given in
  \Cref{rem:discrete-unequal-ranks}.
\end{remark}

\subsection{Continuous-coordinate coarse-space generation and offline--online implementation}
\label{sec:network}

To decouple the predicted basis from the output dimension \(N_h\), we represent
each coarse basis function continuously,
\[
  v_{\theta,j}:\mathcal P\times\Omega\to\mathbb R,
  \qquad j=1,\ldots,n,
\]
and evaluate it at the cell centers of a target grid
\[
  \mathbf{V}_{\theta,h}(\mu)
  =
  \left[
    v_{\theta,j}
    (\mu,\boldsymbol x_i^{\,h})
  \right]_{\substack{i=1,\ldots,N_h\\j=1,\ldots,n}}.
\]

\begin{figure}[!htbp]
  \centering
  \resizebox{\textwidth}{!}{%
    \begin{tikzpicture}[
        x=1cm,y=1cm,
        >=Stealth,
        font=\sffamily\small,
        input/.style={draw=cyan!45!black, rounded corners=3pt, very thick, align=center,
        minimum height=1.05cm, minimum width=3.25cm, fill=cyan!7},
        branch/.style={draw=blue!55!black, rounded corners=3pt, very thick, align=center,
        minimum height=1.1cm, minimum width=3.25cm, fill=blue!7},
        trunk/.style={draw=teal!55!black, rounded corners=3pt, very thick, align=center,
        minimum height=1.1cm, minimum width=3.35cm, fill=teal!8},
        head/.style={draw=violet!55!black, rounded corners=3pt, very thick, align=center,
        minimum height=1.1cm, minimum width=3.35cm, fill=violet!8},
        tentative/.style={draw=orange!65!black, rounded corners=3pt, very thick, align=center,
        minimum height=1.2cm, minimum width=3.6cm, fill=orange!8},
        output/.style={draw=red!55!black, rounded corners=3pt, very thick, align=center,
        minimum height=1.2cm, minimum width=3.8cm, fill=red!7},
        arr/.style={->, very thick, draw=gray!65!black},
        cond/.style={->, very thick, dashed, draw=gray!55}
      ]

      \node[input] (param) at (0,1.15) {Parameter fields\\[-0.1em]$\log_{10}\lambda,\ \phi$};
      \node[input] (point) at (0,-1.15) {Point features\\[-0.1em]$\boldsymbol{x}_i$, Fourier, $\phi(\boldsymbol{x}_i)$};
      \node[branch] (branch) at (3.9,1.15) {Branch encoder\\[-0.1em]$\mathbf{z}_\mu$};
      \node[head] (head) at (7.8,1.15) {Coefficient head\\[-0.1em]$\mathbf C_\theta(\mu)$};
      \node[trunk] (trunk) at (7.8,-1.15) {Trunk network\\[-0.1em]$\mathbf T_{\theta,h}(\mu)$};
      \node[tentative] (vtilde) at (12.05,0) {Tentative basis\\[-0.1em]$\widetilde{\mathbf V}_{\theta,h}=\mathbf T_{\theta,h}\mathbf C_\theta^T$};
      \node[output] (vraw) at (12.05,-2.25) {Boundary-compatible raw basis\\[-0.1em]$\mathbf V^{\rm raw}_{\theta,h}=\mathbf{D}_{\partial,h}\widetilde{\mathbf V}_{\theta,h}$};

      \draw[arr] (param) -- (branch);
      \draw[arr] (branch) -- (head.west);
      \draw[cond] (branch.south east) -- node[above, sloped, font=\sffamily\scriptsize, text=gray!60!black] {conditioning} (trunk.north west);
      \draw[arr] (point) -- (trunk.west);
      \draw[arr] (head.east) -- (vtilde.north west);
      \draw[arr] (trunk.east) -- (vtilde.south west);
      \draw[arr] (vtilde.south) -- node[right, font=\sffamily\scriptsize] {boundary envelope} (vraw.north);
    \end{tikzpicture}
  }
  \caption{The continuous-coordinate network uses the parameter fields and
    target cell centers to generate a tentative basis
    $\widetilde{\mathbf V}_{\theta,h}$. The boundary envelope then produces the
    raw basis $\mathbf V_{\theta,h}^{\rm raw}$ used in the training loss and
  online Galerkin correction}
  \label{fig:deeponet-vraw}
\end{figure}
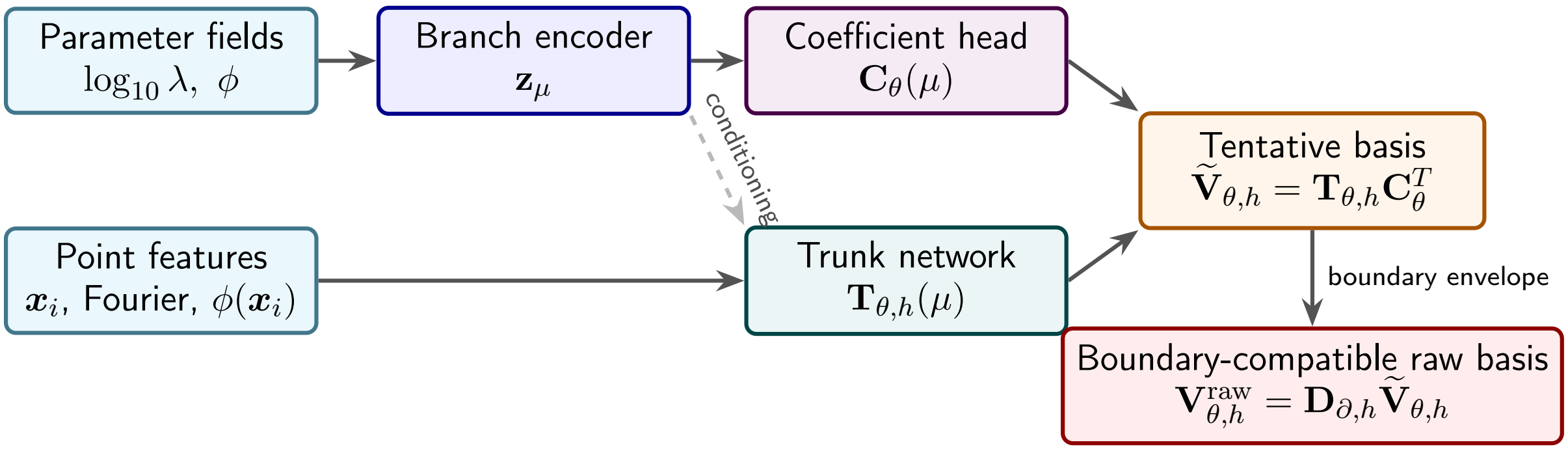

We implement the generator using the DeepONet-type branch--trunk architecture
shown in \Cref{fig:deeponet-vraw}. The branch encoder extracts a global
representation $\mathbf{z}_\mu$ from the parameter fields
$(\log_{10}\lambda,\phi)$. At any spatial location $\boldsymbol x$, the
coordinate trunk produces a dictionary conditioned on the parameter, while a
coefficient head maps $\mathbf{z}_\mu$ to the coefficients used to combine the
dictionary functions. Let
\[
  \mathbf T_{\theta,h}(\mu)
  \in
  \mathbb R^{N_h\times d_{\rm dict}},
  \qquad
  \mathbf C_\theta(\mu)
  \in
  \mathbb R^{n\times d_{\rm dict}}.
\]
The tentative basis is
\begin{equation}
  \widetilde{\mathbf V}_{\theta,h}(\mu)
  =
  \mathbf T_{\theta,h}(\mu)
  \mathbf C_\theta(\mu)^T.
  \label{eq:tentative-basis}
\end{equation}
This architecture separates parameter encoding from the target solver grid:
the parameter fields are first encoded in the network input representation,
while the coordinate trunk can be reevaluated at the cell centers of a
different target grid.

For problems with homogeneous Dirichlet conditions on the entire boundary,
the generated functions are multiplied by a continuous zero-trace
boundary envelope, sampled on the target grid
\begin{equation}
  \mathbf V_{\theta,h}^{\rm raw}(\mu)
  =
  \mathbf{D}_{\partial,h}
  \widetilde{\mathbf V}_{\theta,h}(\mu).
  \label{eq:boundary-envelope}
\end{equation}
This makes the continuous-coordinate representation compatible with the
boundary conditions. For the mixed-boundary moving-interface problem, no
boundary envelope is applied, so $\mathbf{D}_{\partial,h}=\mathbf{I}$. The network widths,
Fourier features, FiLM conditioning, and the explicit boundary-envelope formula
\eqref{eq:boundary-envelope-formula} are given in
Appendix~\ref{sec:network-implementation}.

\paragraph{Offline stage}
Given the training parameter set
$\mathcal P_{\rm train}
=\{\mu_j\}_{j=1}^{N_{\rm train}}$, let $r_j=r_{\rm w}(\mu_j)$. For each
parameter, we construct
\[
  \mu_j
  \longrightarrow
  \mathbf{E}_{\mu_j}
  \longrightarrow
  \boldsymbol{\Phi}_{\mu_j,r_j}
  \longrightarrow
  \mathcal U_{\rm ref,h}^{r_j}(\mu_j),
\]
and minimize the whitened subspace loss, using
\eqref{eq:white-projection-loss} when $r_j=n$ and its one-sided extension when
$r_j<n$. The offline training procedure is summarized in
\Cref{alg:offline-coarse-training}.
\begin{algorithm}[!t]
  \caption{Offline training of parameter-adaptive coarse spaces}
  \label{alg:offline-coarse-training}
  \begin{algorithmic}[1]
    \Require Training parameters
    $\mathcal P_{\rm train}=\{\mu_j\}_{j=1}^{N_{\rm train}}$,
    number $q$ of random initial errors, smoothing count $\nu_{\rm snap}$,
    coarse dimension $n$, number of training epochs $N_{\rm ep}$,
    and singular-value threshold $\varepsilon_{\rm svd}$
    \Ensure Trained network parameters $\theta^\star$

    \Statex Stage 1: Construct POD reference spaces
    \For{$j=1,\ldots,N_{\rm train}$}
    \State Assemble $\mathbf{A}_j=\mathbf{A}_h(\mu_j)$ and form $\mathbf{S}_j=\mathbf{S}_h(\mu_j)$
    \State Draw the random initial-error matrix
    $\mathbf{E}_j^{(0)}=[\mathbf e_{j,1}^{(0)},\ldots,\mathbf e_{j,q}^{(0)}]$
    \State $\mathbf{E}_j\gets \mathbf{S}_j^{\nu_{\rm snap}}\mathbf{E}_j^{(0)}$
    \State Compute the eigendecomposition
    $\mathbf{E}_j^T\mathbf{A}_j\mathbf{E}_j
    =\boldsymbol{\Xi}_j\boldsymbol{\Sigma}_j^2\boldsymbol{\Xi}_j^T$
    \State
    $r_j
    \gets
    \min\{n,\#\{k:\sigma_{j,k}
    \ge\varepsilon_{\rm svd}\sigma_{j,1}\}\}$
    \State
    $\boldsymbol{\Phi}_{j,r_j}
    \gets
    \mathbf{E}_j\boldsymbol{\Xi}_{j,r_j}\boldsymbol{\Sigma}_{j,r_j}^{-1}$
    \State Store
    $(\mu_j,\mathbf{A}_j,\boldsymbol{\Phi}_{j,r_j},r_j)$
    \EndFor

    \Statex Stage 2: Learn the parameter-to-coarse-subspace map
    \State Initialize the network parameters $\theta$
    \For{$t=1,\ldots,N_{\rm ep}$}
    \For{each minibatch $\mathcal B$}
    \ForAll{$j\in\mathcal B$}
    \State
    $\mathbf{V}_j^{\rm raw}
    \gets
    \mathbf{V}_{\theta,h}^{\rm raw}(\mu_j)$
    \State Compute the Cholesky factorization
    $(\mathbf{V}_j^{\rm raw})^T\mathbf{A}_j\mathbf{V}_j^{\rm raw}
    =
    \mathbf{L}_j\mathbf{L}_j^T$
    \State
    $\mathbf{Q}_j
    \gets
    \mathbf{V}_j^{\rm raw}\mathbf{L}_j^{-T}$
    \State
    $\ell_j(\theta)
    \gets
    1-\dfrac1{r_j}
    \|\mathbf{Q}_j^T\mathbf{A}_j\boldsymbol{\Phi}_{j,r_j}\|_F^2$
    \EndFor
    \State
    $\mathcal L_{\mathcal B}
    \gets
    |\mathcal B|^{-1}
    \sum_{j\in\mathcal B}\ell_j(\theta)$
    \State
    $\theta
    \gets
    \operatorname{Adam}
    (\theta,\nabla_\theta\mathcal L_{\mathcal B})$
    \EndFor
    \EndFor
    \State $\theta^\star\gets\theta$
    \State \Return $\theta^\star$
  \end{algorithmic}
\end{algorithm}

\paragraph{Online stage}
For a new parameter $\mu^\star$ and a new target grid with spacing $h^\star$,
the network is evaluated directly at the target cell centers:
\[
  \mathbf{V}^\star
  =
  \mathbf{V}_{\theta^\star,h^\star}^{\rm raw}(\mu^\star).
\]
For numerical stability, we form an orthonormal representative of this raw basis in the
$\mathbf{A}_{h^\star}(\mu^\star)$ inner product
\[
  (\mathbf{V}^\star)^T \mathbf{A}_{h^\star}(\mu^\star)\mathbf{V}^\star
  =
  \mathbf{L}^\star(\mathbf{L}^\star)^T,
  \qquad
  \mathbf{Q}^\star
  =
  \mathbf{V}^\star(\mathbf{L}^\star)^{-T}.
\]
Using this representative, the small Galerkin coarse matrix is
\begin{equation}
  \mathbf{A}_c^\star
  =
  (\mathbf{Q}^\star)^T
  \mathbf{A}_{h^\star}(\mu^\star)
  \mathbf{Q}^\star.
  \label{eq:online-coarse-matrix}
\end{equation}
For the current fine-grid residual $\mathbf r$, the coarse correction is
obtained from
\begin{equation}
  \mathbf{A}_c^\star\mathbf y
  =
  (\mathbf{Q}^\star)^T\mathbf r,
  \qquad
  \mathbf p_{h^\star}
  \leftarrow
  \mathbf p_{h^\star}+\mathbf{Q}^\star\mathbf y.
  \label{eq:online-raw-direct-galerkin}
\end{equation}
With exact energy orthonormalization, $\mathbf A_c^\star=\mathbf I$, and
this gives the same ideal Galerkin correction as the raw basis
$\mathbf V^\star$. The implemented regularized solve uses the Gram matrix
of the basis actually prepared, as specified in \Cref{sec:discrete-setting}.
Expensive random-error sampling and parameter-specific POD are therefore
confined to offline training.

\section{Numerical experiments}
\label{sec:num-experiments}
In this section, we evaluate the proposed coarse space learning method
on three representative problems: static elliptic interface problems,
moving-interface Darcy flow, and nonlinear Forchheimer flow.
We compare the learned coarse spaces with parameter-specific POD
reference spaces, parameter-independent coarse spaces, and GenEO
spaces. The experiments assess four aspects of the proposed method:
\begin{enumerate}
  \item[(i)] \emph{Solver performance relative to the POD reference}:
  the ability to reproduce the solver performance of the
  parameter-specific POD reference space at the same dimension.

  \item[(ii)] \emph{Parameter adaptation}:
  the benefits of parameter-dependent coarse spaces
  over fixed, parameter-independent alternatives.

  \item[(iii)] \emph{Space compactness}:
  the dimension of the learned coarse spaces relative to
  GenEO spaces with comparable convergence rates.

  \item[(iv)] \emph{Online construction efficiency}:
  the computational cost of constructing learned coarse spaces
  through network evaluation rather than repeated snapshot
  or spectral computations.
\end{enumerate}

FGMRES is right-preconditioned in the static interface and Forchheimer
experiments. For the static problems, final acceptance requires an explicitly
recomputed relative residual
$\|\mathbf f_h-\mathbf A_h\mathbf p_h\|_2/\|\mathbf f_h\|_2\le10^{-6}$
and an independent CPU sparse-matrix check. Forchheimer solves retain two
pressure components and use the independent linear and nonlinear checks
specified in Appendix~\ref{sec:finite-precision-verification}.

The moving-interface results in \Cref{tab:hs-metrics} use a separate
FGMRES(50) dataset with the solver's native stopping criterion. Its residual
is updated recursively during the iteration, and the stored success flags
are not an independent recomputation of the final relative residual. These
records therefore document native stopping-test success rather than the
independent residual acceptance used in the static and Forchheimer
experiments.

The online TG preconditioner uses
$(\nu_1,\nu_2)=(1,1)$, a damping factor of $0.8$, and
block-tridiagonal Cholesky local solves.
The local factorizations and coarse spaces remain fixed
within each linear solve.
FGMRES iteration counts represent actual Arnoldi iterations, while the
moving-interface FGMRES counts represent GCR updates. Both exclude the
initial residual evaluation.

\subsection{Comparison methods}
The static interface benchmark uses the following reference and baseline
methods. All methods in this comparison use the same discretization and a
$10\times10$ subdomain partition, with a relative Euclidean residual
tolerance of $10^{-6}$ and a 200-iteration cap. The moving-interface and
Forchheimer experiments use learned and GenEO spaces with the dimensions
and settings specified in their respective subsections.
\begin{itemize}
  \item ASM: the one-level baseline without coarse correction.
  \item POD reference: the parameter-specific POD
    reference space ($n=15$).
  \item Learned: the proposed parameter-adaptive
    learned space ($n=15$).
  \item Global-POD: a fixed parameter-independent
    error space ($n=30$).
  \item DeepONet trunk basis (TB)
    \cite{kopanivcakova2025deeponet}:
    a fixed parameter-independent solution-trained space ($n=15$).
  \item GenEO: a classical coefficient-adaptive
    local spectral coarse space.
\end{itemize}

The POD reference provides a benchmark for evaluating the quality
of the learned coarse spaces, while Global-POD and TB assess the
benefits of parameter adaptation. GenEO serves as a classical
adaptive baseline for comparing solver convergence, coarse-space
dimension, and computational efficiency.

\paragraph{Global-POD coarse space}
\label{sec:global-pod-30}

Global-POD-30 is constructed independently on each target grid using the same
50 offline parameter samples (seeds 91234--91283, disjoint from the held-out test seeds). For each parameter, 15 random initial errors are
propagated by 50 ASM steps, and the leading 15 POD directions are extracted in
the corresponding $\mathbf{A}_h(\mu_i)$ inner product. The resulting
$750$ vectors are Euclidean-normalized and concatenated, after which
a second POD retains 30 global modes. Thus, Global-POD-30 has twice the
dimension of the learned and POD reference spaces and tests whether a larger shared
space can compensate for the absence of parameter adaptation. Because it is
reconstructed independently on every target grid, this baseline specifically
tests the effect of using a parameter-independent space on each grid.

\paragraph{DeepONet trunk-basis space}
\label{sec:deeponet-trunk-basis}

Following \cite{kopanivcakova2025deeponet}, we construct a solution-trained
DeepONet baseline:
\begin{equation}
  \widehat u(\mu,\boldsymbol x)
  =
  \sum_{j=1}^{p} B_j(\mu)T_j(\boldsymbol x),
  \qquad p=128.
  \label{eq:tb-deeponet}
\end{equation}
The branch network takes the parameter fields
$[\log_{10}\lambda,\phi]$ as input, whereas the trunk depends only on the
spatial coordinate. The network is trained on 5000 parameters at
$N_{\mathrm{side}}=63$ by minimizing the relative discrete $\ell^2$ solution
error. 
Because the trunk is independent of $\mu$, any coarse space constructed solely
from its dictionary is parameter independent. We evaluate all 128 trunk
functions on the $N_{\mathrm{side}}=63$ selection grid and apply SVD to the
trunk matrix. Let $\mathbf{W}_{15}$ contain the leading 15 right singular vectors. On a
target grid,
\[
  \widetilde{\mathbf{V}}_h=\mathbf{T}_h\mathbf{W}_{15}.
\]
A thin QR factorization is used to remove numerically dependent columns with
$|R_{jj}|<\varepsilon_{\rm QR}=10^{-10}$. For each online parameter, the
resulting span is then orthogonalized in the current $\mathbf{A}_h(\mu)$ inner product
before forming the Galerkin coarse correction.

\paragraph{GenEO local spectral construction}
\label{sec:geneo-local-spectrum}
For the modified GenEO construction
\cite{spillane2014robust,spillane2014abstract}, we solve the
following generalized eigenproblem on each subdomain:
\begin{equation}
  \mathbf{A}_j^N \mathbf{p}=\xi\,\boldsymbol{\Pi}_j\mathbf{A}_j^N\boldsymbol{\Pi}_j\mathbf{p},
  \qquad
  \xi(\mathbf{p})=
  \frac{\mathbf{p}^T\mathbf{A}_j^N\mathbf{p}}{(\boldsymbol{\Pi}_j\mathbf{p})^T\mathbf{A}_j^N(\boldsymbol{\Pi}_j\mathbf{p})},
  \quad
  (\boldsymbol{\Pi}_j\mathbf{p})^T\mathbf{A}_j^N(\boldsymbol{\Pi}_j\mathbf{p})>0.
  \label{eq:modified-geneo-fvm}
\end{equation}
The matrix $\mathbf{A}_j^N$ is assembled from the FVM transmissibilities. It retains
contributions from interior faces and physical Dirichlet boundaries while
omitting fluxes across artificial subdomain boundaries. The diagonal matrix
$\boldsymbol{\Pi}_j$ contains the ASM partition-of-unity weights. We order the local
eigenmodes by increasing $\xi$ and retain all finite modes satisfying
$\xi\le\tau$, together with the constant null modes of Neumann subdomains.
Linearly dependent directions are removed after global extension. We set
$\tau=0.30$ on every grid throughout the main experiments.

\subsection{Static elliptic interface problem}
\label{sec:elliptic-interface-exp}

\subsubsection{Experimental setup}
\label{sec:static-training-testing}

We consider a family of elliptic interface problems with randomly
generated interfaces. Let $\boldsymbol c=(c_x,c_y)$ denote the
interface center. For $\boldsymbol x\ne\boldsymbol c$, define
$r=\|\boldsymbol x-\boldsymbol c\|_2$ and
$\vartheta=\operatorname{atan2}(y-c_y,x-c_x)$.
The interface is described by the level-set function
\begin{equation}
  \begin{aligned}
    c_x,c_y&\stackrel{\mathrm{iid}}{\sim}\mathcal U(0.4,0.6),
    &R_0&\sim\mathcal U(0.25,0.35),\\
    a_k,b_k&\stackrel{\mathrm{iid}}{\sim}\mathcal N(0,0.05^2),
    &k&=1,\ldots,4,\\
    R_\mu(\vartheta)&=R_0+\sum_{k=1}^{4}
    \bigl(a_k\cos(k\vartheta)+b_k\sin(k\vartheta)\bigr),
    &\phi_\mu(\boldsymbol x)&=r-R_\mu(\vartheta).
  \end{aligned}
  \label{eq:static-random-interface}
\end{equation}
The coefficient is set to $\lambda=10^4$ in
$\{\phi_\mu<0\}$ and $\lambda=1$ in
$\{\phi_\mu\ge0\}$, with a constant source term $f=1$.

\begin{figure}[!htbp]
 \centering
 \begin{subfigure}[t]{0.4\textwidth}
  \includegraphics[width=\linewidth]{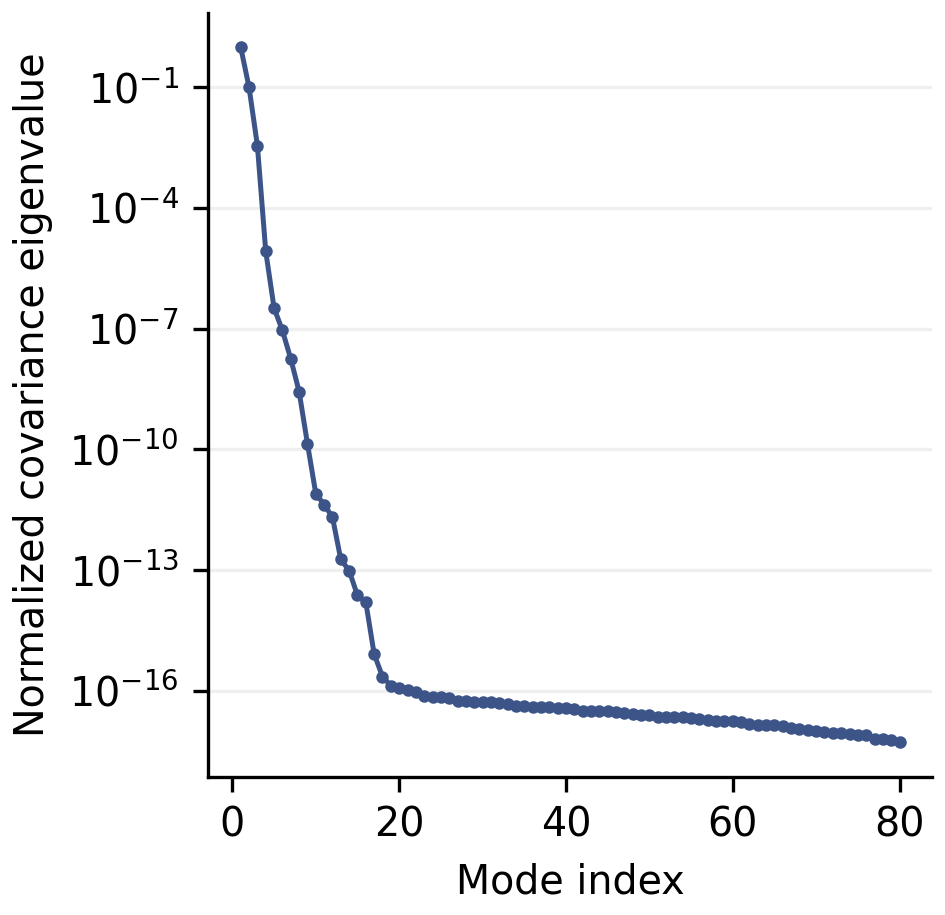}
  \caption{Finite-ensemble covariance spectrum}
  \label{fig:spectral}
 \end{subfigure}\,
 \begin{subfigure}[t]{0.4\textwidth}
  \includegraphics[width=\linewidth]{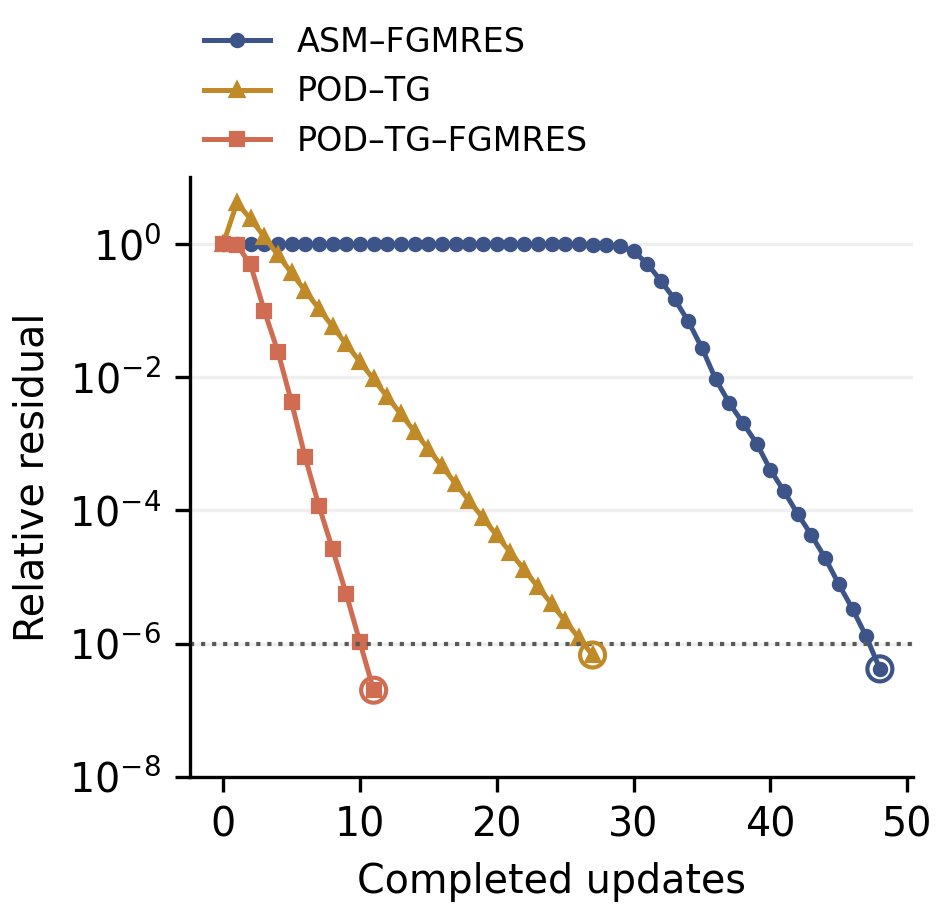}
  \caption{Relative residual histories}
  \label{fig:svd-example-iter}
 \end{subfigure}
 \caption{Smoothed-error covariance spectrum and relative residual
 histories for the static interface problem with a rank-15 POD
 coarse space. The covariance eigenvalues are normalized by
 the largest eigenvalue.}
 \label{fig:spectral-and-iterations}
\end{figure}

We first examine the choice of coarse-space dimension
through the smoothed-error ensemble shown in
\Cref{fig:spectral-and-iterations}.
This diagnostic uses $N_{\rm side}=63$,
$q=200$ independent standard Gaussian initial errors,
and $\nu_{\rm snap}=50$ smoothing steps.
The rapid decay of the energy-covariance eigenvalues
indicates that the dominant smoothed-error components
admit a low-dimensional representation.
Together with the residual histories, this observation
supports the choice of $n=15$ for the coarse space.
The eigenvalues characterize the finite smoothed-error
ensemble rather than the contraction spectrum of the
nonsymmetric ASM error-propagation operator.

The ablation study in
Appendix~\ref{sec:snapshot-strategy-ablation}
further investigates the effects of smoothing steps and
snapshot count. The results show that sufficient smoothing
substantially improves solver convergence, whereas further
increasing the smoothing steps or snapshot count yields
diminishing returns.
We therefore use $q=20$ and $\nu_{\rm snap}=50$
for the main experiments.

Training data are generated exclusively on the
$N_{\mathrm{side}}=63$ grid.
For each training parameter, independent random initial
errors are smoothed by ASM to construct the snapshot ensemble.
The leading $r_{\rm w}=n=15$ POD modes in the
$\mathbf{A}_h$-energy inner product define the
parameter-specific reference space used in the whitened
subspace loss, with the relative singular-value threshold
$\varepsilon_{\rm svd}=10^{-6}$, using the notation of
Algorithms~\ref{alg:reference-space} and~\ref{alg:offline-coarse-training}.
Both the learned and POD reference spaces have dimension $n=15$.

To evaluate cross-grid generalization, the trained model
is applied directly to
$N_{\mathrm{side}}\in\{63,127,255,511,1023\}$
without retraining.
Eight held-out parameter samples are used on each grid,
with all methods evaluated on identical systems.
This setting allows us to assess whether the learned
coarse spaces maintain their effectiveness across
different spatial resolutions.

\subsubsection{Effect of the whitened subspace loss}
\label{sec:static-learned-space-diagnostics}

\begin{figure}[!htbp]
  \centering
  \includegraphics[width=\textwidth]{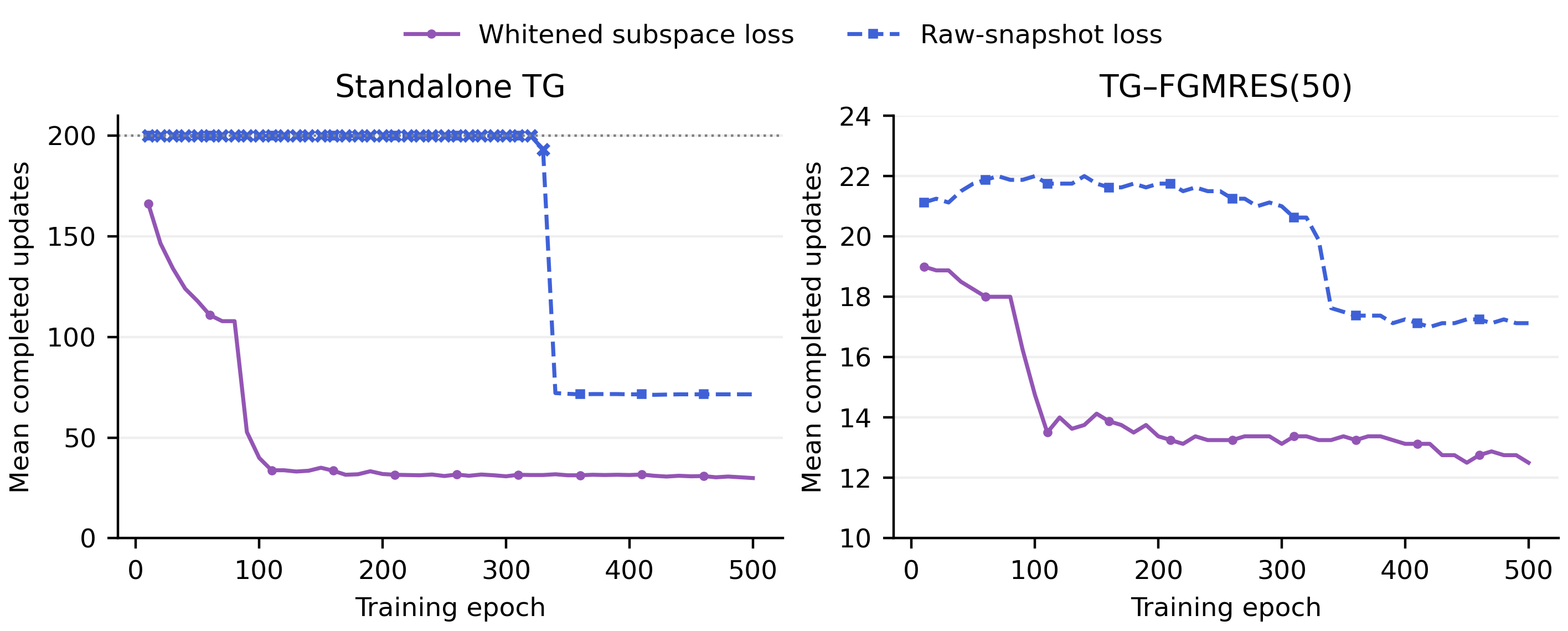}
  \caption{Solver performance over training epochs for the
  whitened subspace and raw-snapshot losses on eight
  independent validation fields.}
  \label{fig:static-whitening-validation-history}
\end{figure}

We investigate the effectiveness of POD whitening by
comparing the raw-snapshot loss
\eqref{eq:spectral-weighted-residual}
with the proposed whitened subspace loss
\eqref{eq:white-projection-loss}.
Both models use identical training data, network architecture,
initialization, optimizer, and training schedule,
thereby isolating the loss function as the primary experimental variable.
\Cref{fig:static-whitening-validation-history}
compares the convergence behavior of the resulting
coarse spaces throughout training on an independent
validation set.
The whitened loss produces coarse spaces with better
solver performance at the end of training, reducing
the iterations required by both standalone TG and
TG-preconditioned FGMRES.
This demonstrates that POD whitening improves the
solver effectiveness of the learned coarse spaces.

The results in \Cref{tab:static-whitening-ablation}
further demonstrate the advantage of the whitened
subspace loss across all tested grid resolutions.
Compared with the raw-snapshot loss, whitening
consistently reduces the iteration counts and their
variability, indicating faster and more stable convergence.
The improvement is particularly pronounced for standalone
TG on finer grids, where the raw-snapshot loss exhibits
greater sensitivity to the problem size.
\begin{table}[!htbp]
 \centering
 \caption{Effect of POD whitening on solver performance over eight paired
 fields per grid. Entries give mean $\pm$ population standard deviation
 over accepted solves. TG counts sweeps and FGMRES counts Arnoldi updates,
 excluding the initial check. Solves that do not meet the residual criterion
 within 200 updates are counted as failures.
 Unmarked entries have 8/8 accepted solves. $\dagger$: 7/8.}
 \label{tab:static-whitening-ablation}
 \small\setlength{\tabcolsep}{4pt}
 \begin{tabular}{@{}rcccc@{}}
 \toprule
 $N_{\rm side}$ & Raw TG & Whitened TG &
 Raw TG--FGMRES & Whitened TG--FGMRES \\
 \midrule
 $63$ & $67.50\pm14.44$ & $\boldsymbol{30.75\pm3.03}$ &
 $16.75\pm0.83$ & $\boldsymbol{13.38\pm0.48}$ \\
 $127$ & $98.38\pm22.18$ & $\boldsymbol{47.38\pm6.28}$ &
 $20.62\pm0.99$ & $\boldsymbol{16.62\pm0.70}$ \\
 $255$ & $120.12\pm12.95$ & $\boldsymbol{51.00\pm3.74}$ &
 $22.88\pm0.78$ & $\boldsymbol{17.50\pm0.50}$ \\
 $511$ & $125.00\pm17.19$ & $\boldsymbol{60.12\pm3.85}$ &
 $24.12\pm0.78$ & $\boldsymbol{18.88\pm0.78}$ \\
 $1023$ & $148.57\pm26.78^{\dagger}$ & $\boldsymbol{70.38\pm6.78}$ &
 $26.62\pm0.86$ & $\boldsymbol{20.88\pm0.60}$ \\
 \bottomrule
 \end{tabular}
\end{table}

\subsubsection{Cross-grid performance and online cost}
\label{sec:static-cross-grid-comparison}

We next compare the proposed method with the reference
and baseline coarse spaces across different grid resolutions.
All methods are evaluated on identical linear systems,
right-hand sides, and initial guesses.
The results in \Cref{tab:static-cross-grid-iterations}
demonstrate the effectiveness of the learned coarse spaces
in terms of solver convergence, parameter adaptation,
and coarse-space compactness.

\begin{table}[!htbp]
 \centering
 \caption{Right-preconditioned FGMRES(50) performance on eight paired
 fields per grid. Entries are actual Arnoldi-update mean $\pm$
 population standard deviation over accepted solves.
 The GenEO dimension is the mean constructed coarse-space dimension.
 Solves that do not meet the residual criterion within 200 updates
 are counted as failures.
 Unmarked entries have 8/8 accepted solves.
 $\dagger$: 7/8. $\ddagger$: 1/8.}
 \label{tab:static-cross-grid-iterations}
 \footnotesize\setlength{\tabcolsep}{3pt}
 \resizebox{\textwidth}{!}{
 \begin{tabular}{@{}rccccccc@{}}
 \toprule
 $N_{\rm side}$ &
 \multicolumn{5}{c}{FGMRES updates} &
 \multicolumn{2}{c}{GenEO} \\
 \cmidrule(lr){2-6}\cmidrule(lr){7-8}
 & ASM & POD-ref-15 & Learned-15 &
 Global-POD-30 & TB-SVD-15 & Updates & Dimension \\
 \midrule
 $63$ & $47.88\pm0.60$ & $10.50\pm0.50$ &
 $\boldsymbol{13.38\pm0.48}$ & $20.00\pm0.50$ &
 $19.75\pm0.66$ & $17.88\pm2.03$ & $52.1$ \\
 $127$ & $93.50\pm9.38$ & $13.62\pm0.48$ &
 $\boldsymbol{16.62\pm0.70}$ & $23.88\pm0.93$ &
 $23.75\pm0.97$ & $19.88\pm0.60$ & $86.9$ \\
 $255$ & $96.88\pm2.62$ & $14.62\pm0.48$ &
 $\boldsymbol{17.50\pm0.50}$ & $27.00\pm0.71$ &
 $26.62\pm0.70$ & $24.25\pm0.66$ & $112.8$ \\
 $511$ & $125.14\pm21.18^{\dagger}$ & $16.38\pm0.48$ &
 $\boldsymbol{18.88\pm0.78}$ & $30.25\pm0.83$ &
 $29.25\pm1.30$ & $17.88\pm0.78$ & $255.8$ \\
 $1023$ & $146.00\pm0.00^{\ddagger}$ & $18.00\pm0.00$ &
 $\boldsymbol{20.88\pm0.60}$ & $34.75\pm0.97$ &
 $32.00\pm0.71$ & $19.12\pm0.60$ & $406.2$ \\
 \bottomrule
 \end{tabular}}
\end{table}

\begin{figure}[!htbp]
  \centering
  \includegraphics[width=\textwidth]{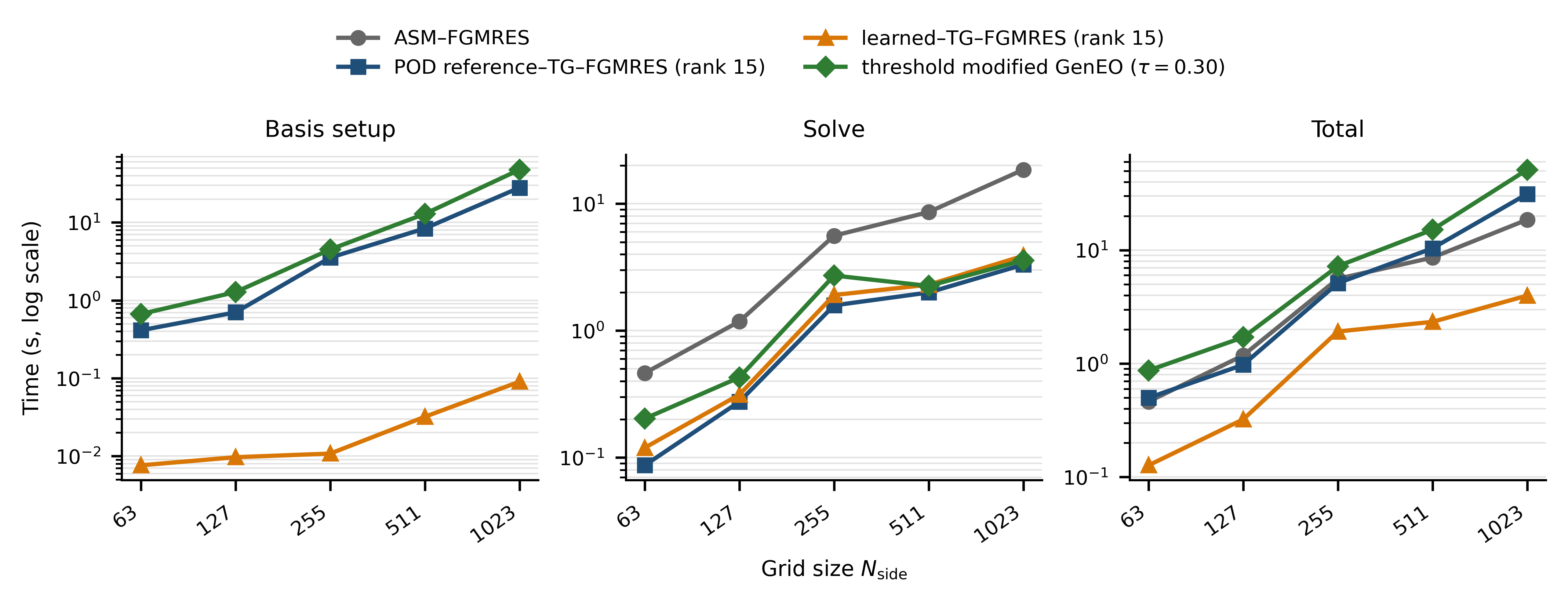}
  \caption{Coarse-space construction, FGMRES solution, and
  total wall-clock times for the static interface problem.
  All timings are obtained using the same right-preconditioned
  FGMRES configuration.}
  \label{fig:static-formal-timing}
\end{figure}

\Cref{fig:static-formal-timing} compares coarse-space construction costs,
iterative solution times, and total solver times. It shows that replacing
repeated snapshot generation and local spectral computations with network
evaluation reduces coarse-space construction costs.

The results in \Cref{tab:static-cross-grid-iterations} and \Cref{fig:static-formal-timing}
illustrate the four aspects introduced at the beginning of this section.
\paragraph{(i) Solver performance relative to the POD reference}
The learned coarse spaces achieve FGMRES iteration counts close to those
of the parameter-specific POD reference spaces across all tested grids.
These results indicate that the learned spaces provide effective coarse
corrections for the prescribed smoother and Krylov solver, even though their
geometric distances from the POD reference spaces may remain substantial.
\Cref{tab:preconditioner-perturbation} quantifies these geometric discrepancies.

\paragraph{(ii) Effect of parameter adaptation}
The proposed method consistently outperforms the
parameter-independent Global-POD and DeepONet trunk
basis approaches.
Notably, the learned coarse spaces achieve faster
convergence even when compared with the higher-dimensional
Global-POD spaces.
This highlights the advantage of adapting the coarse
space to the underlying coefficient field rather than
relying on a fixed space shared across parameters.

\paragraph{(iii) Compactness}
Compared with GenEO, the proposed method achieves
comparable or better convergence with substantially
smaller coarse spaces over the tested problems.
While the dimension of the GenEO coarse space increases
markedly with grid refinement, the learned space maintains
a fixed dimension.
The advantage becomes particularly pronounced on finer
grids, where the coarse-space dimension required by
GenEO is considerably larger.

\paragraph{(iv) Online construction efficiency}
A key practical advantage of the proposed method is its
efficient construction of parameter-dependent coarse
spaces through network evaluation.
Unlike snapshot-based or local spectral approaches,
the learned method does not require repeated snapshot
generation or spectral decomposition for each new
parameter instance.
The online stage consists of network evaluation,
basis preparation, and assembly of the small Galerkin coarse matrix, which is
computationally inexpensive compared with the iterative solution of the
fine-grid system.

\subsection{Moving-interface Hele--Shaw flow}
\label{sec:darcy-moving-interface}

This experiment considers a time-dependent sequence of high-contrast
operators generated by moving-interface Hele--Shaw flow. For each sampled
operator, we compare a newly generated learned parameter-adaptive space with
a rebuilt GenEO space and a GenEO space frozen at the first sampled operator.

The domain is $[0,1]^2$ with grid resolution
$N_{\mathrm{side}}=1023$. We impose $p=1$ and $p=0$ on the left and right
boundaries, respectively, together with homogeneous no-flow conditions on the
top and bottom boundaries, and set the body source to zero. The permeability
is $\kappa=1$. The background and droplet viscosities are $1$ and $10^{-4}$,
respectively, giving a mobility contrast of $10^4$ for
$\lambda=\kappa/\mu_f$.

The initial interface is centered at $(0.5,0.5)$ and has polar radius
$R(\vartheta)=0.18[1+0.2\cos(5\vartheta)]$. The level-set equation is advanced
for 1000 transport steps with a CFL number of $0.3$. Face Darcy fluxes are
averaged to obtain cell-centered velocities. The advection term
$\bv\cdot\nabla\phi$ is discretized by signed upwinding, and the resulting
semidiscrete equation is advanced with the SSPRK(2,2) scheme, equivalently
Heun's method \cite{gottlieb2001strong}. The learned model generates a
15-dimensional coarse space and is trained for 20 epochs on 5000 parameter
samples at $N_{\mathrm{side}}=64$. Each training sample uses $q=50$ random
initial errors smoothed by 50 ASM steps, with POD-whitening rank
$r_{\rm w}=15$.

To compare the three coarse-space strategies on exactly the same operator
sequence, we first generate the complete 1000-step trajectory with
learned TG--FGMRES. We then select 80 mobility operators at uniformly spaced
transport steps and solve this fixed sequence independently with each
method. The first linear system in the sampled sequence is initialized from
zero, and each subsequent solve uses the solution of the previous sampled
system as its initial guess. Trajectory-generation time is excluded from all
reported timings.
The rebuilt GenEO method reconstructs the threshold-selected
modified GenEO space with $\tau=0.30$ for every sampled operator. The
frozen GenEO method constructs the space only for the first sampled
operator and reuses the same span thereafter. The learned method
generates a new 15-dimensional space for the current operator by network
inference.

\begin{figure}[!htbp]
  \centering
  \includegraphics[width=\textwidth]{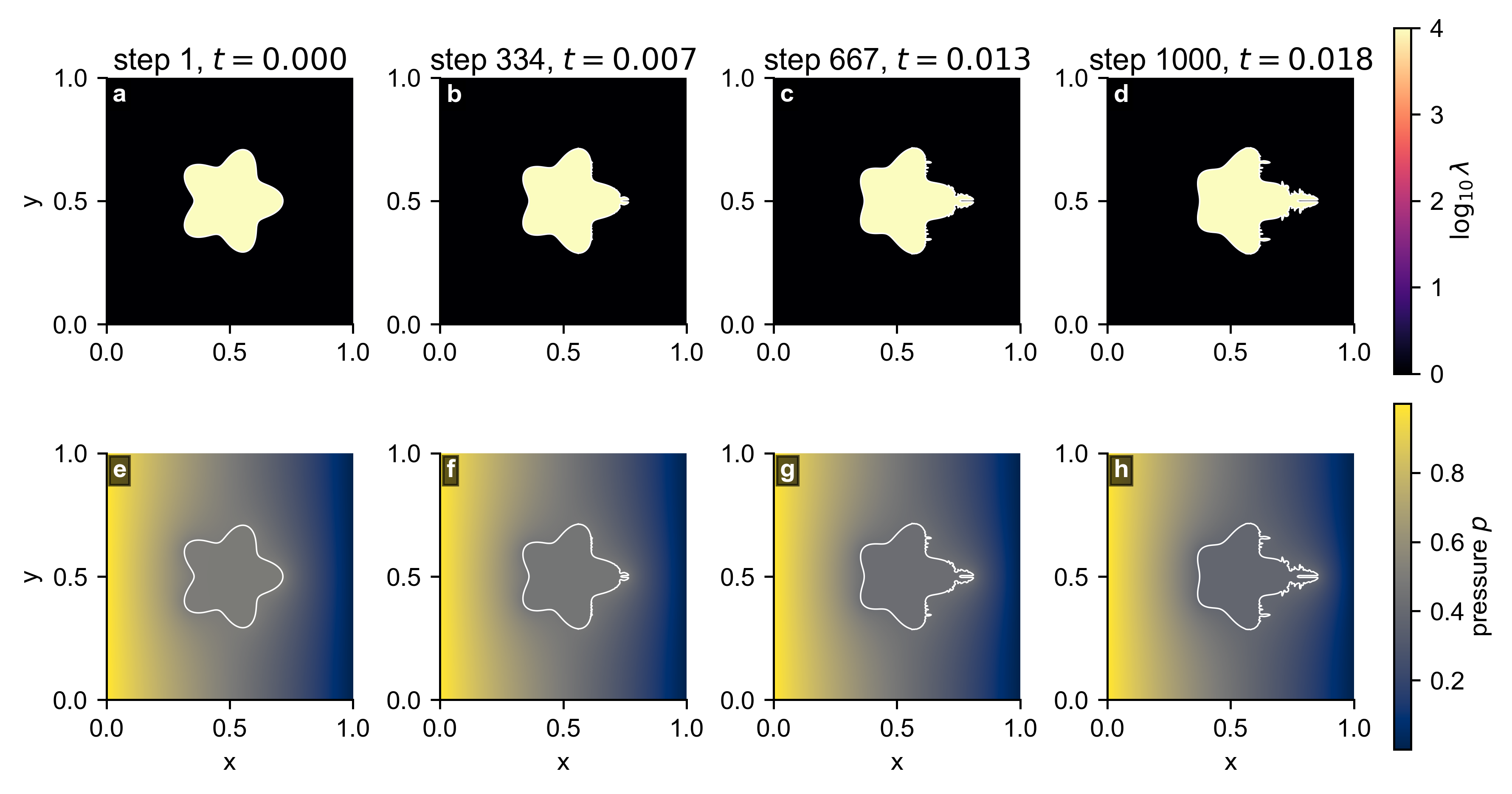}
  \caption{Mobility and pressure fields at four representative times in the
  moving-interface trajectory. White curves mark $\phi=0$}
  \label{fig:moving-representative-fields}
\end{figure}

\Cref{fig:moving-representative-fields} shows four representative states along
the trajectory. The droplet moves to the right and deforms as the mobility
field evolves. Over the full trajectory, the relative area change is about
$0.25\%$, and the centroid displacement is
$\Delta x\approx0.0513$.

\Cref{tab:hs-metrics} reports the results over the 80 sampled operators.
Iterations are given as the mean $\pm$ population standard deviation, while
the time columns retain cumulative measurements from a separate reference
timing run. The success column records the native FGMRES stopping-test outcome
over the 80 operators.

\begin{table}[!htbp]
  \centering
  \caption{FGMRES iterations, coarse dimensions, and wall-clock times for the
  moving-interface operator sequence.}
  \label{tab:hs-metrics}

  \begingroup
  \small
  \renewcommand{\arraystretch}{1.12}
  \setlength{\tabcolsep}{4.5pt}
  \resizebox{\textwidth}{!}{%
    \begin{tabular}{@{}lcccccc@{}}
      \toprule
      Method & FGMRES iterations & Recorded success & Mean coarse dim. &
      Construction (s) & Solve (s) & Total (s) \\
      \midrule
      Learned-15
      & $9.34\pm0.79$ & $100\%$ & $15.0$
      & $8.49$ & $84.87$ & $94.45$ \\
      Rebuilt GenEO
      & $8.32\pm0.93$ & $100\%$ & $418.8$
      & $8236.95$ & $61.26$ & $8299.35$ \\
      Frozen GenEO
      & $10.61\pm1.04$ & $100\%$ & $420.0$
      & $108.48$ & $73.59$ & $183.16$ \\
      \bottomrule
    \end{tabular}%
  }
  \endgroup
\end{table}

The learned method's mean FGMRES count lies between those of rebuilt and
frozen GenEO. All three methods report native stopping-test success on the
full 80-operator sequence. The principal differences in these records lie
in coarse-space dimension and construction cost. The learned space has
fixed dimension $15$, whereas the GenEO coarse spaces contain about $420$
global directions. Thus, similar native stopping-test outcomes are obtained
with a learned coarse space that is approximately 28 times smaller.

The separate reference timing measurements compare the construction costs
of repeated parameter adaptation. Frozen GenEO substantially reduces construction
cost relative to rebuilding the spectral space for every sampled operator, but
its total time is $183.16$ s, compared with $94.45$ s for the learned method.
Rebuilding GenEO for every sampled operator is much more expensive, with a
cumulative construction time of $8236.95$ s. In contrast, generating 80
parameter-adaptive learned spaces requires only $8.49$ s.

\subsection{Nonlinear Forchheimer flow in channelized media}
\label{sec:forchheimer-channelized}

This experiment studies coarse-space adaptation within Picard linearization
for high-contrast channelized media. Following
\cite{vasilyeva2026operator}, we use the domain $[0,1]^2$ with homogeneous
Dirichlet conditions on all four sides and grid resolution
$N_{\mathrm{side}}=1023$. The reference permeability $\kappa_0(\bx)$ is $1$
in the background and $10^6$ in the connected channels. We set
$\beta=10^{-4}$, $\mu_f=\rho=1$, and $f(\bx)\equiv10$. The
pressure-dependent permeability is
\[
  \kappa(\bx,p)
  =
  \kappa_0(\bx)\exp[\alpha(p-p_0)],
  \qquad p_0=0,
\]
and we consider both the Forchheimer-only case $\alpha=0$ and the
pressure-sensitive case $\alpha=10$. At each Picard step, the network
generates a coarse space for the linear system associated with the current
frozen mobility $\lambda^{[r]}$. The model outputs $n=30$ basis functions and
is trained for 30 epochs on 500 channel samples, with 450 samples at
$N_{\mathrm{side}}=63$ and 50 at $N_{\mathrm{side}}=127$. Each sample uses
$q=30$ random initial errors smoothed by 20 ASM steps. Because the spectral
tail varies among channel samples, the reliable whitening rank
$r_{\rm w}$ is selected between 4 and 30 using the relative singular-value
criterion in \Cref{sec:loss-spec}. The network still outputs a
30-dimensional space and is trained to contain the reliable POD directions.

The Picard iteration starts from $\mathbf p_h^{[0]}=\mathbf0$ and is capped
at 100 steps, with at most 300 FGMRES updates per linear system. At each step, $|\nabla_h \mathbf p_h^{[r]}|$ is evaluated by
centered differences in the interior and one-sided boundary differences, after which
$\lambda^{[r]}$ is frozen and $\mathbf A_h(\lambda^{[r]})$ is assembled. The previous
Picard iterate is used as the initial guess for the FGMRES solve. The
nonlinear iteration terminates when both
\[
  \frac{\|\mathbf p_h^{[r+1]}-\mathbf p_h^{[r]}\|_2}
  {\|\mathbf p_h^{[r+1]}\|_2}
  \le 10^{-6}
\]
and
\[
  \mathcal R_{\rm nl}(\mathbf p_h^{[r+1]})
  :=
  \frac{
    \|\mathbf f_h
    -\mathbf A_h(\lambda(\mathbf p_h^{[r+1]}))\mathbf p_h^{[r+1]}\|_2
  }{
    \|\mathbf f_h\|_2
  }
  \le 10^{-6}.
\]
We compare three coarse-space strategies. The learned method generates a
new 30-dimensional space from the current mobility at every Picard step.
Rebuilt GenEO reconstructs the threshold-selected modified GenEO
space with $\tau=0.30$ at every step. Frozen GenEO constructs the
GenEO span only at the first Picard step and reuses it thereafter. For all
three methods, the Galerkin coarse matrix is assembled from the current
operator at every linearization.

\begin{figure}[!htb]
  \centering
  \includegraphics[width=\textwidth]
  {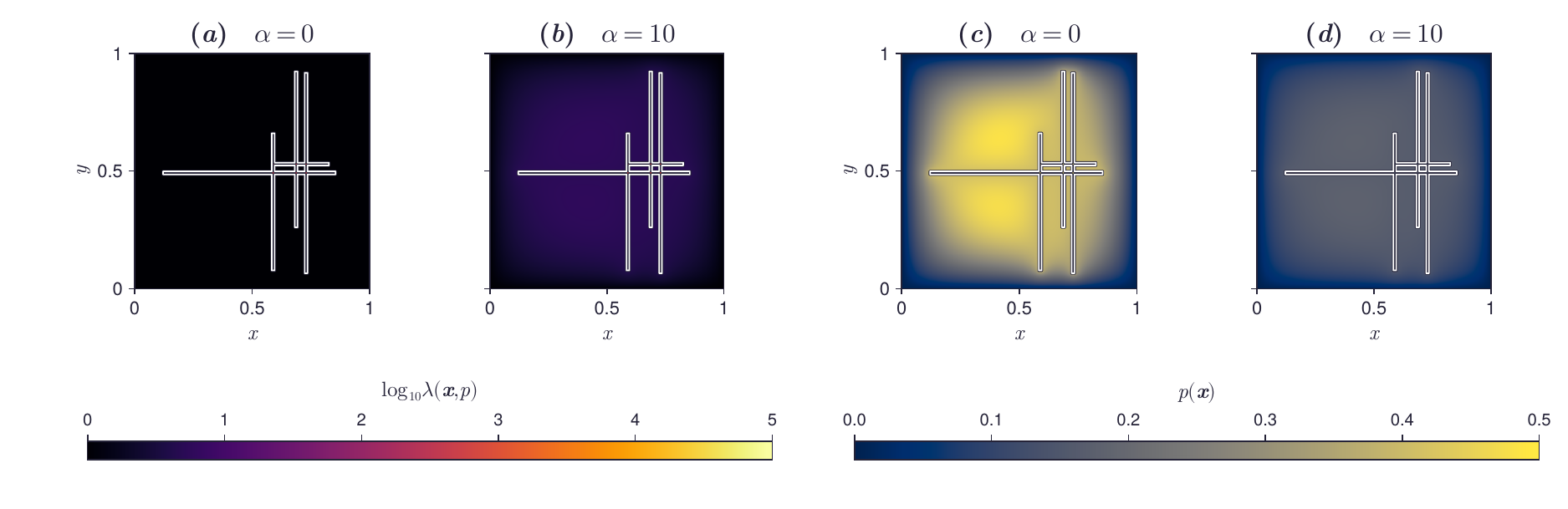}
  \caption{Terminal mobility and pressure fields for Learned-30 on channel
  field 0 at $\alpha=0$ and $\alpha=10$. White curves mark the channel boundaries.}
  \label{fig:forchheimer-representative-fields}
\end{figure}

\Cref{fig:forchheimer-representative-fields} shows terminal mobility and
pressure fields for Learned-30 on channel field 0. The two values of
$\alpha$ produce different effective mobilities while preserving the same
underlying channel geometry.

\Cref{tab:forchheimer-picard-ablation} compares convergence, coarse-space
dimensions, and wall-clock times for the three strategies. All methods
satisfy the prescribed linear and nonlinear stopping criteria in both
permeability settings. The learned method achieves iteration counts
comparable to rebuilt GenEO with a substantially smaller coarse space.
Network evaluation reduces the cost of adapting the coarse space to the
evolving mobility by avoiding repeated spectral construction. The separate
reference timings show that these construction savings yield lower total
solution times than both GenEO strategies.
\begin{table}[!htb]
  \centering
  \caption{Comparison of coarse-space strategies for Forchheimer flow using eight paired fields. $\pm$ denotes population standard
    deviation. }
  \label{tab:forchheimer-picard-ablation}

  \begingroup
  \small
  \renewcommand{\arraystretch}{1.12}
  \setlength{\tabcolsep}{3.5pt}
  \resizebox{\textwidth}{!}{%
    \begin{tabular}{@{}clccccccc@{}}
      \toprule
      $\alpha$ & Method & Picard steps & FGMRES/step &
      \shortstack{Final relative\\nonlinear residual} & Mean coarse dim. &
      Construction (s) & Solve (s) & Total (s) \\
      \midrule
      $0$ & Learned-30
      & $25.62\pm0.48$ & $21.60$ & $7.16\times10^{-7}$
      & $30.0$ & $2.63$ & $101.01$ & $105.61$ \\
      $0$ & Rebuilt GenEO
      & $25.38\pm0.48$ & $21.67$ & $7.96\times10^{-7}$
      & $440.8$ & $1536.28$ & $102.60$ & $1642.07$ \\
      $0$ & Frozen GenEO
      & $25.62\pm0.48$ & $22.23$ & $7.22\times10^{-7}$
      & $439.5$ & $58.28$ & $108.32$ & $168.96$ \\
      \addlinespace
      $10$ & Learned-30
      & $24.62\pm0.48$ & $19.27$ & $6.86\times10^{-7}$
      & $30.0$ & $2.37$ & $85.33$ & $89.48$ \\
      $10$ & Rebuilt GenEO
      & $24.38\pm0.48$ & $20.44$ & $7.60\times10^{-7}$
      & $439.0$ & $1265.34$ & $90.74$ & $1357.93$ \\
      $10$ & Frozen GenEO
      & $24.50\pm0.50$ & $21.23$ & $7.53\times10^{-7}$
      & $439.5$ & $58.33$ & $94.79$ & $154.94$ \\
      \bottomrule
    \end{tabular}%
  }
  \endgroup

\end{table}

\section{Subspace approximation and perturbation analysis}
\label{sec:convergence}

This section analyzes how approximation of the reference coarse space affects the
resulting TG preconditioner and the convergence of the associated
Krylov solver. The analysis proceeds from energy-norm subspace approximation
to operator perturbation estimates and finite-step GMRES bounds. Throughout
this section, the coefficient parameter and grid are fixed. The arguments
allow for nonsymmetric smoothing operators and are complemented by numerical
diagnostics across grid resolutions.

\subsection{Discrete setting and coarse-space approximation}
\label{sec:discrete-setting}

Let $\mathbf{A}=\mathbf{A}_h(\mu)\succ0$ and introduce the energy inner product
$\langle\mathbf{x},\mathbf{y}\rangle_{\mathbf{A}}=
\mathbf{x}^T\mathbf{A}\mathbf{y}$ and the induced operator norm
$\|\mathbf{K}\|_{\mathbf{A}}=
\|\mathbf{A}^{1/2}\mathbf{K}\mathbf{A}^{-1/2}\|_2$.
The error-propagation operator of the averaged ASM smoother is
\[
  \mathbf{S}=\mathbf{I}-\omega\mathbf{W}^{-1}
  \left(\sum_k\mathbf{R}_k^T\mathbf{B}_k\mathbf{R}_k\right)\mathbf{A},
\]
where $\mathbf{B}_k$ denotes a fixed linear local solution map.
Since each restriction selects a set of distinct global degrees of freedom,
$\mathbf{A}_k=\mathbf{R}_k\mathbf{A}\mathbf{R}_k^T$ is SPD.
The analysis covers both the exact local inverse
$\mathbf{B}_k=\mathbf{A}_k^{-1}$ and the diagonally shifted map
$\mathbf{B}_k=(\mathbf{A}_k+\varepsilon_{\rm loc}\mathbf{I})^{-1}$.
The left multiplicity weighting is retained, so $\mathbf{S}$ may be
nonsymmetric.

Write $\mathbf{R}=\mathbf{S}^{\nu_1}$ and
$\mathbf{L}=\mathbf{S}^{\nu_2}$ for the pre- and postsmoothing operators.
For a full-column-rank coarse basis $\mathbf{V}$, define its energy
projector and the TG error-propagation operator by
\begin{equation}
  \mathbf{P}_{\mathbf{V}}
  =\mathbf{V}(\mathbf{V}^T\mathbf{A}\mathbf{V})^{-1}\mathbf{V}^T\mathbf{A},
  \qquad
  \mathbf{E}_{\mathbf{V}}
  =\mathbf{L}(\mathbf{I}-\mathbf{P}_{\mathbf{V}})\mathbf{R}.
  \label{eq:discrete-tg-error-propagation}
\end{equation}
The online experiments use $(\nu_1,\nu_2)=(1,1)$, while the results below
apply to general fixed linear operators $\mathbf{L}$ and $\mathbf{R}$.

The reference space is constructed by energy-POD of the snapshot ensemble
$\mathbf{X}=\mathbf{S}^{\nu_{\rm snap}}\mathbf{E}^{(0)}$.
By \Cref{prop:weighted-pod-optimality}, this construction gives the optimal
energy-norm approximation of the snapshot ensemble and its associated
ellipsoid. In the static experiments, $\nu_{\rm snap}=50$, compared with a
single presmoothing step in the online TG method. The learning objective
therefore measures approximation of the selected reference subspace,
while the effect on online error propagation is determined by
\eqref{eq:discrete-tg-error-propagation}.

For completeness, we account for the regularization used in the practical coarse
solve. Let $\mathbf{V}$ denote the prepared basis used in that solve and set
$\mathbf{G}_{\mathbf{V}}=\mathbf{V}^T\mathbf{A}\mathbf{V}$.
The regularized coarse correction is
\[
  \mathbf{C}_{\mathbf{V}}
  =\mathbf{V}(\mathbf{G}_{\mathbf{V}}+
  \varepsilon_{\rm coarse}\mathbf{I})^{-1}\mathbf{V}^T\mathbf{A},
  \qquad \varepsilon_{\rm coarse}\ge0.
\]
Its difference from the energy projector is given exactly by
\begin{equation}
  \|\mathbf{P}_{\mathbf{V}}-\mathbf{C}_{\mathbf{V}}\|_{\mathbf{A}}
  =\frac{\varepsilon_{\rm coarse}}
  {\lambda_{\min}(\mathbf{G}_{\mathbf{V}})+\varepsilon_{\rm coarse}}
  =:\tau_{\mathbf{V}}.
  \label{eq:coarse-regularization-gap}
\end{equation}
To see this, the matrix
$\mathbf{U}=\mathbf{A}^{1/2}\mathbf{V}\mathbf{G}_{\mathbf{V}}^{-1/2}$
has orthonormal columns and
\[
  \mathbf{A}^{1/2}(\mathbf{P}_{\mathbf{V}}-\mathbf{C}_{\mathbf{V}})
  \mathbf{A}^{-1/2}
  =\mathbf{U}\varepsilon_{\rm coarse}
  (\mathbf{G}_{\mathbf{V}}+\varepsilon_{\rm coarse}\mathbf{I})^{-1}
  \mathbf{U}^T.
\]
Consequently, for
$\mathbf{E}_{\mathbf{V}}^{\rm reg}=
\mathbf{L}(\mathbf{I}-\mathbf{C}_{\mathbf{V}})\mathbf{R}$,
\[
  \|\mathbf{E}_{\mathbf{V}}^{\rm reg}
    -\mathbf{E}_{\mathbf{V}}\|_{\mathbf{A}}
  \le\|\mathbf{L}\|_{\mathbf{A}}
       \|\mathbf{R}\|_{\mathbf{A}}\tau_{\mathbf{V}}.
\]
The Gram matrix in \eqref{eq:coarse-regularization-gap} is formed from the
prepared basis, as scalar regularization depends on the basis
representation. We use $\varepsilon_{\rm coarse}=10^{-10}$ for the static and
moving-interface problems and $10^{-8}$ for the Forchheimer Picard problems.
The coarse regularization is separate from the local-solve shift
$\varepsilon_{\rm loc}$ and the linear residual tolerance
$\mathrm{tol}_{\rm lin}$.

\subsection{Perturbation of learned two-grid operators}
\label{sec:discrete-perturbation}

Let $\mathbf{P}_\theta$ and $\mathbf{P}_{\rm ref}$ be the energy projectors
onto the learned and reference spaces, both of dimension $n$.
We write $\mathcal{L}_{\rm white}$ for the projector loss in
\eqref{eq:white-projection-loss}, evaluated between these spaces on the
same target grid. The following result connects this loss with the
perturbation of the TG error-propagation operator.

\begin{theorem}[Discrete subspace perturbation]
  \label{thm:discrete-subspace-perturbation}
  For arbitrary fixed linear $\mathbf{L}$ and $\mathbf{R}$,
  \begin{align}
    \mathbf{E}_\theta-\mathbf{E}_{\rm ref}
      &=\mathbf{L}(\mathbf{P}_{\rm ref}-\mathbf{P}_\theta)\mathbf{R},
      \label{eq:exact-discrete-perturbation}\\
    \delta_{\mathbf{A}}
      :=\|\mathbf{P}_\theta-\mathbf{P}_{\rm ref}\|_{\mathbf{A}}
      &=\sin\vartheta_{\max}
      \le\min\{1,\sqrt{n\mathcal{L}_{\rm white}}\},
      \label{eq:discrete-projector-loss}\\
    \|\mathbf{E}_\theta-\mathbf{E}_{\rm ref}\|_{\mathbf{A}}
      &\le\|\mathbf{L}\|_{\mathbf{A}}
             \|\mathbf{R}\|_{\mathbf{A}}
             \min\{1,\sqrt{n\mathcal{L}_{\rm white}}\}.
      \label{eq:discrete-subspace-perturbation-bound}
  \end{align}
  In the online configuration $(\nu_1,\nu_2)=(1,1)$, the prefactor in
  \eqref{eq:discrete-subspace-perturbation-bound} is
  $\|\mathbf{S}\|_{\mathbf{A}}^2$.
\end{theorem}

\begin{proof}
  Subtracting the two instances of
  \eqref{eq:discrete-tg-error-propagation} gives
  \eqref{eq:exact-discrete-perturbation}.
  Conjugation by $\mathbf{A}^{1/2}$ transforms the energy projectors into
  Euclidean orthogonal projectors. The nonzero eigenvalues of their
  difference occur in pairs $\pm\sin\vartheta_j$.
  Hence its operator norm is $\sin\vartheta_{\max}\le1$, and
  \eqref{eq:principal-angle-loss} gives
  $\sin^2\vartheta_{\max}\le
    \sum_j\sin^2\vartheta_j=n\mathcal{L}_{\rm white}$.
  The final estimate follows by submultiplicativity.
\end{proof}

\Cref{thm:discrete-subspace-perturbation} establishes the dependence of
the TG operator on the learned coarse space. For regularized coarse
corrections, the perturbation bound acquires the additional term
$\|\mathbf{L}\|_{\mathbf{A}}\|\mathbf{R}\|_{\mathbf{A}}
(\tau_\theta+\tau_{\rm ref})$.
The bound based on the largest principal angle can be sharpened by
incorporating the actions of the smoother on the subspace discrepancy.

\begin{proposition}[Operator-weighted defect]
  \label[proposition]{prop:operator-weighted-defect}
  Let $\mathbf{D}=\mathbf{P}_\theta-\mathbf{P}_{\rm ref}$ and define
  $\Delta_{\rm op}:=\|\mathbf{L}\mathbf{D}\mathbf{R}\|_{\mathbf{A}}$
  and $\eta_{\mathbf{R}}:=\|\mathbf{D}\mathbf{R}\|_{\mathbf{A}}$.
  Then
  $\Delta_{\rm op}=\|\mathbf{E}_\theta-
  \mathbf{E}_{\rm ref}\|_{\mathbf{A}}$ and
  \begin{equation}
    \Delta_{\rm op}\le\ell_{\mathbf{D}}\eta_{\mathbf{R}}
    \le\|\mathbf{L}\|_{\mathbf{A}}
       \|\mathbf{R}\|_{\mathbf{A}}\delta_{\mathbf{A}},
    \qquad
    \ell_{\mathbf{D}}:=
    \sup_{0\ne\mathbf{y}\in\operatorname{range}(\mathbf{D})}
      \frac{\|\mathbf{L}\mathbf{y}\|_{\mathbf{A}}}
           {\|\mathbf{y}\|_{\mathbf{A}}},
    \label{eq:operator-weighted-defect}
  \end{equation}
  where $\ell_{\mathbf{D}}=0$ if $\mathbf{D}=0$.
  For energy-orthonormal bases $\mathbf{Q}_\theta$ and
  $\mathbf{Q}_{\rm ref}$, let
  $\mathbf{Z}=[\mathbf{Q}_\theta,\mathbf{Q}_{\rm ref}]$,
  $\mathbf{J}=\operatorname{diag}(\mathbf{I}_n,-\mathbf{I}_n)$, and
  $\mathbf{R}^*=\mathbf{A}^{-1}\mathbf{R}^T\mathbf{A}$.
  Define
  \[
    \mathbf{G}_{\mathbf{L}}
      =(\mathbf{L}\mathbf{Z})^T\mathbf{A}(\mathbf{L}\mathbf{Z}),
    \qquad
    \mathbf{G}_{\mathbf{R}}
      =(\mathbf{R}^*\mathbf{Z})^T\mathbf{A}(\mathbf{R}^*\mathbf{Z}).
  \]
  Then
  \begin{equation}
    \Delta_{\rm op}
      =\|\mathbf{G}_{\mathbf{L}}^{1/2}
         \mathbf{J}\mathbf{G}_{\mathbf{R}}^{1/2}\|_2,
    \label{eq:small-gram-perturbation}
  \end{equation}
  where the square roots are the positive semidefinite roots.
  Replacing $\mathbf{G}_{\mathbf{L}}$ by
  $\mathbf{Z}^T\mathbf{A}\mathbf{Z}$ gives the corresponding expression
  for $\eta_{\mathbf{R}}$.
\end{proposition}

\begin{proof}
  Because $\mathbf{D}\mathbf{R}\mathbf{x}\in
  \operatorname{range}(\mathbf{D})$,
  $\|\mathbf{L}\mathbf{D}\mathbf{R}\mathbf{x}\|_{\mathbf{A}}
  \le\ell_{\mathbf{D}}\|\mathbf{D}\mathbf{R}\mathbf{x}\|_{\mathbf{A}}$.
  The second inequality follows from
  $\ell_{\mathbf{D}}\le\|\mathbf{L}\|_{\mathbf{A}}$ and
  $\eta_{\mathbf{R}}\le\delta_{\mathbf{A}}\|\mathbf{R}\|_{\mathbf{A}}$.
  For \eqref{eq:small-gram-perturbation}, observe that
  $\mathbf{D}=\mathbf{Z}\mathbf{J}\mathbf{Z}^T\mathbf{A}$.
  Setting $\mathbf{F}=\mathbf{A}^{1/2}\mathbf{L}\mathbf{Z}$ and
  $\mathbf{G}=\mathbf{A}^{1/2}\mathbf{R}^*\mathbf{Z}$, the
  energy-coordinate matrix of $\mathbf{L}\mathbf{D}\mathbf{R}$ is
  $\mathbf{F}\mathbf{J}\mathbf{G}^T$.
  Polar factorization gives
  $\mathbf{F}=\mathbf{U}_{\mathbf{F}}\mathbf{G}_{\mathbf{L}}^{1/2}$
  and
  $\mathbf{G}=\mathbf{U}_{\mathbf{G}}\mathbf{G}_{\mathbf{R}}^{1/2}$.
  The partial isometries preserve the spectral norm of the product on
  the ranges of the two Gram roots, including when either Gram matrix
  is singular. This proves \eqref{eq:small-gram-perturbation}.
  The formula for $\eta_{\mathbf{R}}$ follows by taking
  $\mathbf{L}=\mathbf{I}$.
\end{proof}

The quantity $\Delta_{\rm op}$ measures the TG perturbation after
presmoothing and postsmoothing. It thus incorporates the interaction
between the coarse-space discrepancy and the smoother, which is absent
from $\delta_{\mathbf{A}}$ alone. Formula
\eqref{eq:small-gram-perturbation} expresses this quantity through Gram
matrices of order $2n$. Their construction involves the relevant operator
and adjoint actions. For $(\nu_1,\nu_2)=(1,1)$,
\[
  \mathbf{G}_{\mathbf{R}}
  =(\mathbf{S}^T\mathbf{A}\mathbf{Z})^T\mathbf{A}^{-1}
    (\mathbf{S}^T\mathbf{A}\mathbf{Z}),
\]
which requires solves with $\mathbf{A}$.
We use $\Delta_{\rm op}$ as an offline diagnostic of the TG operator
perturbation.

\begin{remark}[Unequal reliable ranks]
  \label{rem:discrete-unequal-ranks}
  Consider a learned space of dimension $n$ and a reference space of
  reliable rank $1\le r<n$, as in \Cref{rem:rank-deficient-targets}.
  The full projector distance is then
  $\|\mathbf{P}_\theta-\mathbf{P}_{\rm ref}\|_{\mathbf{A}}=1$.
  A one-sided estimate is more informative in this setting:
  \[
    \|(\mathbf{I}-\mathbf{P}_\theta)\mathbf{R}\|_{\mathbf{A}}
    \le
    \|(\mathbf{I}-\mathbf{P}_{\rm ref})\mathbf{R}\|_{\mathbf{A}}
    +\min\{1,\sqrt{r\mathcal{L}_{\rm white}^{(r)}}\}
      \|\mathbf{P}_{\rm ref}\mathbf{R}\|_{\mathbf{A}}.
  \]
  This follows by decomposing $\mathbf{R}$ with respect to
  $\mathbf{P}_{\rm ref}$ and using
  $\|(\mathbf{I}-\mathbf{P}_\theta)\mathbf{P}_{\rm ref}\|_{\mathbf{A}}^2
  \le\min\{1,r\mathcal{L}_{\rm white}^{(r)}\}$.
  The estimate characterizes the capture of reliable reference
  directions by the rank-adaptive Forchheimer coarse spaces.
\end{remark}

\subsection{Implications for right-preconditioned FGMRES}
\label{sec:discrete-krylov}

Let $\mathbf{B}_j$ be the fixed TG preconditioner associated with
$\mathbf{E}_j=\mathbf{I}-\mathbf{B}_j\mathbf{A}$, and put
$\mathbf{T}_j=\mathbf{A}\mathbf{B}_j$ for
$j\in\{\theta,{\rm ref}\}$.
In the dual energy norm
$\|\mathbf{r}\|_{\mathbf{A}^{-1}}=
(\mathbf{r}^T\mathbf{A}^{-1}\mathbf{r})^{1/2}$,
\begin{equation}
  \mathbf{T}_\theta-\mathbf{T}_{\rm ref}
  =-\mathbf{A}(\mathbf{E}_\theta-\mathbf{E}_{\rm ref})
    \mathbf{A}^{-1},
  \qquad
  \|\mathbf{T}_\theta-\mathbf{T}_{\rm ref}\|_{\mathbf{A}^{-1}}
  =\Delta_{\rm op}.
  \label{eq:residual-operator-isometry}
\end{equation}
The equality follows from the isometry between the energy norm on
errors and the dual energy norm on residuals. Equivalently,
$\|\mathbf{B}_\theta-\mathbf{B}_{\rm ref}\|_{
\mathbf{A}^{-1}\to\mathbf{A}}=\Delta_{\rm op}$.
The same relation holds for the actual regularized TG maps when
$\Delta_{\rm op}$ is replaced by their operator difference norm.

\begin{theorem}[Finite-step polynomial perturbation]
  \label{thm:krylov-polynomial-perturbation}
  Let
  $\varepsilon_T=
  \|\mathbf{T}_\theta-\mathbf{T}_{\rm ref}\|_{\mathbf{A}^{-1}}$
  and assume $\|\mathbf{T}_j\|_{\mathbf{A}^{-1}}\le M_j$.
  For any residual polynomial
  $p(z)=\sum_{k=0}^m a_kz^k$ with $a_0=1$, define
  \[
    K_p(M_\theta,M_{\rm ref})
    =\sum_{k=1}^m |a_k|\sum_{i=0}^{k-1}
      M_\theta^{k-1-i}M_{\rm ref}^{i}.
  \]
  Then
  \begin{equation}
    \|p(\mathbf{T}_\theta)-p(\mathbf{T}_{\rm ref})\|_{\mathbf{A}^{-1}}
    \le\varepsilon_TK_p(M_\theta,M_{\rm ref}).
    \label{eq:polynomial-perturbation}
  \end{equation}
  For exact-arithmetic GMRES minimizing the
  $\mathbf{A}^{-1}$ residual norm, the residual after $m$ steps from a
  common initial residual $\mathbf{r}_0\ne0$ satisfies
  \begin{equation}
    \frac{\|\mathbf{r}_{\theta,m}^{(\mathbf{A}^{-1})}\|_{\mathbf{A}^{-1}}}
         {\|\mathbf{r}_0\|_{\mathbf{A}^{-1}}}
    \le
    \frac{\|p(\mathbf{T}_{\rm ref})\mathbf{r}_0\|_{\mathbf{A}^{-1}}}
         {\|\mathbf{r}_0\|_{\mathbf{A}^{-1}}}
    +\varepsilon_TK_p(M_\theta,M_{\rm ref}).
    \label{eq:weighted-gmres-comparison}
  \end{equation}
\end{theorem}

\begin{proof}
  The noncommutative telescoping identity gives
  \[
    \mathbf{T}_\theta^k-\mathbf{T}_{\rm ref}^k
    =\sum_{i=0}^{k-1}
      \mathbf{T}_\theta^{k-1-i}
      (\mathbf{T}_\theta-\mathbf{T}_{\rm ref})
      \mathbf{T}_{\rm ref}^{i}.
  \]
  Taking norms and summing over the polynomial coefficients yields
  \eqref{eq:polynomial-perturbation}. The GMRES minimal-residual
  property, applied to the candidate polynomial $p$, and the triangle
  inequality then give \eqref{eq:weighted-gmres-comparison}.
\end{proof}

The estimate compares the learned operator with the action of a chosen
residual polynomial on the reference operator. For a fixed linear
preconditioner $\mathbf{B}$, exact-arithmetic FGMRES, with
$\mathbf{z}_j=\mathbf{B}\mathbf{v}_j$, generates the same correction space
as right-preconditioned GMRES for $\mathbf{T}=\mathbf{A}\mathbf{B}$
\cite{saad1993flexible,saad2003iterative}.
The implemented solver minimizes the Euclidean residual norm within
each restart cycle. The residual norms satisfy
\[
  \frac{\|\mathbf{r}\|_2}{\sqrt{\lambda_{\max}(\mathbf{A})}}
  \le\|\mathbf{r}\|_{\mathbf{A}^{-1}}
  \le\frac{\|\mathbf{r}\|_2}{\sqrt{\lambda_{\min}(\mathbf{A})}}.
\]
Applying the Euclidean minimal-residual property directly to the candidate
polynomial and using the norm equivalence yields the following bound:
\begin{equation}
  \frac{\|\mathbf{r}_{\theta,m}^{(2)}\|_2}{\|\mathbf{r}_0\|_2}
  \le
  \frac{\|p(\mathbf{T}_{\rm ref})\mathbf{r}_0\|_2}
       {\|\mathbf{r}_0\|_2}
  +c_{\mathbf{A}}(\mathbf{r}_0)
       \varepsilon_TK_p(M_\theta,M_{\rm ref}),
  \label{eq:euclidean-krylov-comparison}
\end{equation}
where
\[
  c_{\mathbf{A}}(\mathbf{r}_0)
  =\frac{\sqrt{\lambda_{\max}(\mathbf{A})}
          \|\mathbf{r}_0\|_{\mathbf{A}^{-1}}}{\|\mathbf{r}_0\|_2}
  \le\sqrt{\kappa_2(\mathbf{A})}.
\]
The norm equivalence is applied only to the polynomial perturbation term.
The reference-polynomial term is already measured in the Euclidean norm.
In particular,
$\|\mathbf{K}\|_2\le
\sqrt{\kappa_2(\mathbf{A})}\|\mathbf{K}\|_{\mathbf{A}^{-1}}$.
These finite-step bounds apply to a common cycle-initial residual and
quantify the sensitivity to the TG operator perturbation.
The factor $K_p$ captures polynomial amplification. If
$M=\max(M_\theta,M_{\rm ref})$, then
$K_p\le\sum_{k=1}^m k|a_k|M^{k-1}$.

In the static interface and Forchheimer implementations, $\mathbf{A}$,
the local factors, the predicted coarse basis, and the coarse factor remain
fixed during each FGMRES solve and are updated between parameter or Picard solves.
Convergence is assessed using the explicitly recomputed relative residual
$\|\mathbf{f}-\mathbf{A}\mathbf{x}\|_2/\|\mathbf{f}\|_2
\le\mathrm{tol}_{\rm lin}=10^{-6}$.
For high-contrast Forchheimer problems, the implementation also uses two
float64 pressure components, compensated accumulation, and residuals
formed from the stored system matrix. Additional finite-precision and
nonlinear residual checks are described in Appendix
\ref{sec:finite-precision-verification}.

\subsection{Numerical assessment}
\label{sec:discrete-cross-grid-discussion}

\Cref{tab:preconditioner-perturbation} compares geometric subspace errors,
TG operator perturbations, and FGMRES convergence for the learned and
POD-reference spaces. All diagnostics are evaluated on each target grid
using its energy metric, reference space, and smoother, with the static
solver settings of \Cref{sec:elliptic-interface-exp}. The paired fields
are drawn from the held-out static cohort, and the final column counts
actual FGMRES updates. Their geometries differ across grids, so the results
assess cross-grid performance rather than refinement of a single
coefficient field.

\begin{table}[!htbp]
  \centering
  \caption{Same-grid perturbation diagnostics, averaged over two paired
  samples. The product column uses the measured $\delta_{\mathbf{A}}$,
  and Ratio is the mean of
  $\|\mathbf{S}\|_{\mathbf{A}}^2\delta_{\mathbf{A}}/\Delta_{\rm op}$.
  The final column reports POD-reference and learned FGMRES updates,
  respectively, with independently verified relative residuals below
  $10^{-6}$.}
  \label{tab:preconditioner-perturbation}
  \small
  \setlength{\tabcolsep}{4pt}
  \begin{tabular}{@{}rcccccc@{}}
    \toprule
    $N_{\rm side}$ & $\mathcal{L}_{\rm white}$ & $\delta_{\mathbf{A}}$
    & $\Delta_{\rm op}$ &
    $\|\mathbf{S}\|_{\mathbf{A}}^2\delta_{\mathbf{A}}$
    & Ratio & Ref./Learned\\
    \midrule
    63   & 0.247 & 0.998 & 0.645 & 1.051  & 1.63 & 10.5/13.0\\
    127  & 0.519 & 0.997 & 0.769 & 1.472  & 1.91 & 14.0/16.5\\
    255  & 0.526 & 0.994 & 0.809 & 3.091  & 3.82 & 14.5/17.0\\
    511  & 0.625 & 0.995 & 0.986 & 6.399  & 6.50 & 16.5/19.5\\
    1023 & 0.687 & 0.993 & 1.180 & 11.690 & 9.91 & 18.0/21.0\\
    \bottomrule
  \end{tabular}
\end{table}

The learned spaces achieve FGMRES iteration counts close to those of the
parameter-specific POD references across the tested grids. This supports
the ability of the coordinate-based generator to provide effective coarse
spaces beyond the training resolution. The large worst-case projector
distances reveal discrepancies in some coarse-space directions, and the
loss-based estimate in \eqref{eq:discrete-projector-loss} is saturated
throughout the table. The largest principal angle alone therefore gives
an incomplete assessment of the learned spaces' value for Krylov
preconditioning.

The operator-weighted perturbation is consistently smaller than the
product-norm bound, and the gap widens on finer grids. Accounting for
pre- and postsmoothing thus gives a sharper assessment of how a
coarse-space discrepancy affects the TG operator. Coarse regularization
makes a negligible contribution to the perturbation bounds in these tests.
The finite-step polynomial perturbation bounds remain conservative, with amplification
terms too large to give quantitative predictions of the observed Krylov
performance. The perturbation analysis characterizes operator sensitivity,
while the paired FGMRES measurements provide direct evidence of solver
effectiveness.

We also examine the energy-norm behavior of the TG error operator.
Let $\mathbf{H}=\mathbf{A}^{1/2}\mathbf{S}\mathbf{A}^{-1/2}$ and
$\boldsymbol{\Pi}_{\mathbf{V}}
=\mathbf{A}^{1/2}\mathbf{P}_{\mathbf{V}}\mathbf{A}^{-1/2}$.
Because $\mathbf{H}\boldsymbol{\Pi}_{\mathbf{V}}\mathbf{H}$ has rank at
most $n$, the singular-value min--max principle yields
\[
  \|\mathbf{S}(\mathbf{I}-\mathbf{P}_{\mathbf{V}})
    \mathbf{S}\|_{\mathbf{A}}
  =\|\mathbf{H}^2-\mathbf{H}\boldsymbol{\Pi}_{\mathbf{V}}
    \mathbf{H}\|_2
  \ge\sigma_{n+1}(\mathbf{H}^2).
\]
The same rank bound applies to a regularized coarse correction.
Separate constant-field and high-contrast interface tests yield numerical
estimates of this singular-value lower bound above unity on the finer
grids. These numerical witnesses provide evidence against strict one-step
energy contraction for the chosen coarse dimension and smoother on those
configurations. They support the relevance of a perturbation framework
that accommodates nonsymmetric smoothing without a stationary contraction
assumption.

These diagnostics highlight the importance of assessing coarse spaces
through their interaction with the smoother and Krylov solver. The
analysis provides fixed-grid perturbation estimates, while the paired
FGMRES results support the cross-grid effectiveness of the learned spaces.

\section{Conclusions}
\label{sec:conclusion}

We developed a parameter-adaptive coarse-space learning method for
TG preconditioning of high-contrast parametric elliptic
problems. The method constructs reference coarse spaces from
smoothed-error snapshots through energy-inner-product POD and learns
the parameter-to-subspace map using a basis-invariant, whitened
subspace loss. The resulting coarse spaces can be evaluated directly
on unseen grids and incorporated into standard Galerkin coarse
corrections. Numerical experiments involving static interface
problems, moving-interface Darcy flow, and nonlinear Forchheimer
flow demonstrate effective iterative convergence across varying
parameters and grid resolutions, together with substantial savings
in repeated coarse-space construction.

The theoretical analysis relates coarse-space approximation errors
to perturbations of the generally nonsymmetric TG preconditioner.
An operator-weighted discrepancy characterizes the combined effects
of pre- and postsmoothing and admits a small-Gram representation.
Finite-step residual-polynomial estimates further establish
conditional perturbation bounds for fixed right-preconditioned
Krylov systems. These results provide a mathematical characterization
of learned coarse-space perturbations without requiring symmetry,
normality, or stationary energy contraction of the smoother.

Future work will focus on learning objectives more closely aligned
with Krylov convergence, adaptive selection of coarse-space
dimensions, and improved cross-grid generalization. Extensions to
indefinite and vector-valued problems, particularly Helmholtz and
Maxwell equations, will also be investigated, with coarse-space
constructions and stability analyses tailored to their respective
operator structures.

\section*{Acknowledgments}
This work was supported by the National Natural Science Foundation of China (No. 12625119, W2612010).

\section*{Data availability}
The data and source code supporting this study, including the Lean~4
formal verification files, will be made publicly available in the project
repository at \url{https://github.com/cuichen1996/ROM-TG}.

\newcounter{mainfigures}
\newcounter{maintables}
\setcounter{mainfigures}{\value{figure}}
\setcounter{maintables}{\value{table}}
\begin{appendices}
  \setcounter{figure}{\value{mainfigures}}
  \setcounter{table}{\value{maintables}}
  \renewcommand{\thefigure}{\arabic{figure}}
  \renewcommand{\thetable}{\arabic{table}}

  \section{Interface-aware cell-centered FVM and convergence verification}
  \label{sec:interface-aware-ccfvm}

  Let
  \[
    \Omega=(0,1)^2,\qquad h=\frac{1}{N_{\mathrm{side}}},
  \]
  with $N_{\mathrm{side}}$ uniform control volumes in each coordinate direction.
  The cell-center coordinates are
  \[
    x_i=\left(i-\tfrac12\right)h,\qquad
    y_j=\left(j-\tfrac12\right)h,\qquad
    1\le i,j\le N_{\mathrm{side}}.
  \]
  Integrating the governing equation over the control volume $K_{ij}$ with
  cell-center index $P=(i,j)$ and applying the divergence theorem gives
  \begin{equation}
    -\int_{\partial K_{ij}}\lambda\nabla p\cdot\mathbf n\,\mathrm ds
    =
    \int_{K_{ij}}f\,\mathrm d\boldsymbol x.
    \label{eq:ccfvm-control-volume-balance}
  \end{equation}
  Here $\mathbf n$ denotes the outward unit normal to $\partial K_{ij}$.
  On this uniform grid, each control-volume edge has length $h$, and the
  distance between adjacent cell centers is also $h$. The integrated normal
  flux across a shared face is approximated by the pressure difference between
  the two adjacent cell centers multiplied by the face transmissibility.
  The right-hand side of
  \eqref{eq:ccfvm-control-volume-balance} is approximated by $h^2f_{i,j}$,
  where $f_{i,j}=f(x_i,y_j)$. For adjacent control volumes indexed by $P$ and
  $Q$, denote their cell-center coordinates by $\boldsymbol x_P$ and
  $\boldsymbol x_Q$, respectively, and write
  $\lambda_P=\lambda(\boldsymbol x_P)$, $\lambda_Q=\lambda(\boldsymbol x_Q)$,
  $\phi_P=\phi(\boldsymbol x_P)$, and $\phi_Q=\phi(\boldsymbol x_Q)$.

  We detect a single interface crossing along the segment joining adjacent
  cell centers by a sign change of the level set values. If $\phi_P\phi_Q\ge0$, we use the standard
  harmonic transmissibility
  \begin{equation}
    T_{PQ}^{\mathrm{harm}}
    =
    \frac{2\lambda_P\lambda_Q}{\lambda_P+\lambda_Q}.
    \label{eq:harmonic-transmissibility}
  \end{equation}
  If $\phi_P\phi_Q<0$, the cell centers lie on opposite sides of the interface.
  Linear interpolation of the level set function along
  $[\boldsymbol x_P,\boldsymbol x_Q]$ gives the estimated interface location as
  a fraction of the distance between the cell centers:
  \begin{equation}
    \theta_{PQ}
    =
    \frac{|\phi_P|}{|\phi_P|+|\phi_Q|},
    \label{eq:theta-interface-location}
  \end{equation}
  where $\theta_{PQ}$ is measured from $\boldsymbol x_P$ toward
  $\boldsymbol x_Q$. Treating the diffusion resistances over the two segment
  lengths, $\theta_{PQ}h$ and $(1-\theta_{PQ})h$, as connected in series gives
  the interface-aware transmissibility
  \begin{equation}
    T_{PQ}^{\mathrm{ia}}
    =
    \left(
      \frac{\theta_{PQ}}{\lambda_P}
      +
      \frac{1-\theta_{PQ}}{\lambda_Q}
    \right)^{-1}.
    \label{eq:interface-aware-transmissibility}
  \end{equation}
  When $\theta_{PQ}=1/2$, the estimated interface location is midway between
  the two cell centers, and \eqref{eq:interface-aware-transmissibility} reduces
  exactly to the harmonic formula \eqref{eq:harmonic-transmissibility}.
  Thus, the interface-aware formula generalizes the standard harmonic
  transmissibility to account for an off-center interface.

  For a control volume indexed by $P=(i,j)$ that does not share a face with
  $\partial\Omega$, let $E,W,N,S$ denote the indices of its east, west, north,
  and south neighbors, respectively. Write the corresponding face
  transmissibilities as $T_E=T_{PE}$, $T_W=T_{PW}$, $T_N=T_{PN}$, and
  $T_S=T_{PS}$. Substituting the numerical fluxes into
  \eqref{eq:ccfvm-control-volume-balance} yields the five-point stencil
  \begin{equation}
    \begin{aligned}
      &(T_E+T_W+T_N+T_S)p_{i,j}
      -T_Ep_{i+1,j}
      -T_Wp_{i-1,j}\\
      &\hspace{4.2em}
      -T_Np_{i,j+1}
      -T_Sp_{i,j-1}
      =h^2f_{i,j}.
    \end{aligned}
    \label{eq:ccfvm-five-point-stencil}
  \end{equation}
  Thus, the nonzero entries in the row of $\mathbf{A}_h$ corresponding to $P$ are
  \[
    \begin{aligned}
      (\mathbf{A}_h)_{P,P}&=T_E+T_W+T_N+T_S,\\
      (\mathbf{A}_h)_{P,E}&=-T_E,\qquad (\mathbf{A}_h)_{P,W}=-T_W,\\
      (\mathbf{A}_h)_{P,N}&=-T_N,\qquad (\mathbf{A}_h)_{P,S}=-T_S.
    \end{aligned}
  \]

  For a control volume adjacent to a Dirichlet boundary, the distance from its
  center to each boundary face is $h/2$. Using $\lambda_P$ on this half-cell
  segment gives the boundary transmissibility
  \begin{equation}
    T_{P,\partial\Omega}=2\lambda_P.
    \label{eq:ccfvm-dirichlet-boundary-transmissibility}
  \end{equation}
  For each Dirichlet boundary face, $T_{P,\partial\Omega}$ is added to the
  diagonal entry $(\mathbf{A}_h)_{P,P}$, and $T_{P,\partial\Omega}g_D$ is added to the
  corresponding right-hand-side entry, with $g_D$ evaluated at the face
  midpoint. Here $g_D=0$, so only the diagonal contribution remains.
  For a corner control volume, both boundary faces contribute.

  Let $\mathcal E_h$ denote the set of interior faces, each identified with
  an unordered pair $\{P,Q\}$ of adjacent cell indices, and let $\Gamma_D$
  denote the Dirichlet boundary. Define $\mathcal B_{D,h}$ as the set of pairs
  $(P,F)$ such that $F$ is a face of the control volume indexed by $P$ and
  $F\subset\Gamma_D$. For each such face, let
  $T_{P,F}=T_{P,\partial\Omega}=2\lambda_P$. For any cell-centered vector
  $\mathbf v$, the discrete energy under homogeneous Dirichlet conditions is
  \begin{equation}
    \mathbf v^T\mathbf{A}_h\mathbf v
    =
    \sum_{\{P,Q\}\in\mathcal E_h}
    T_{PQ}(v_P-v_Q)^2
    +
    \sum_{(P,F)\in\mathcal B_{D,h}}
    T_{P,F}v_P^2.
    \label{eq:ccfvm-discrete-energy}
  \end{equation}
  Each interior face is counted once, and each Dirichlet boundary face
  contributes separately, including both boundary faces of a corner control
  volume. Since $T_{PQ}=T_{QP}$, the matrix $\mathbf{A}_h$ is symmetric. All
  transmissibilities are positive, so zero energy forces equal values on
  adjacent cells and zero values on cells adjacent to $\Gamma_D$. The grid is
  connected and has Dirichlet boundary faces. Hence zero energy implies
  $\mathbf v=\mathbf 0$. Therefore, $\mathbf{A}_h$ is SPD.

  We compare the convergence behavior of the two transmissibility formulas
  using a manufactured solution with a circular interface. On
  $\Omega=(0,1)^2$, let
  \[
    g(x,y)=\sin(\pi x)\sin(\pi y),\qquad
    \phi(x,y)=\left(x-\tfrac12\right)^2+\left(y-\tfrac12\right)^2-R^2,
  \]
  where $R=0.23$ is the interface radius. Define the diffusion coefficient and
  manufactured solution by
  \[
    \lambda(x,y)=
    \begin{cases}
      \kappa,&\phi(x,y)<0,\\
      1,&\phi(x,y)\ge0,
    \end{cases}
    \qquad
    p(x,y)=g(x,y)\frac{\phi(x,y)}{\lambda(x,y)},
    \qquad \kappa=10^4.
  \]
  Since $\lambda$ is constant on each side of the interface, the source term
  is $f=-\Delta(g\phi)$ in each subregion. Using $\Delta g=-2\pi^2g$ gives
  \begin{equation}
    f=2\pi^2\phi g-2\nabla g\cdot\nabla\phi-g\,\Delta\phi.
    \label{eq:manufactured-interface-source}
  \end{equation}
  Because $g=0$ on $\partial\Omega$, the manufactured solution satisfies the
  homogeneous Dirichlet condition. On the interface, $\phi=0$, so the pressure
  traces from both sides vanish. Moreover, in each subregion,
  $\lambda\nabla p=\phi\nabla g+g\nabla\phi$, whose traces agree across the
  interface. Thus, both the pressure jump and the normal-flux jump vanish,
  and the homogeneous interface conditions are satisfied.

  We use uniform cell-centered grids with
  $N_{\mathrm{side}}\in\{32,64,128,256,512\}$. The discrete systems are solved
  using an independent sparse direct solver, with maximum relative algebraic residual $7.9334\times10^{-12}$. Let $d\in\{\mathrm{harm},\mathrm{ia}\}$ identify the
  discretization, and let $\mathbf p_h^{(d)}$ and $\mathbf p_h^{\mathrm{ex}}$
  denote the numerical solution and the exact solution sampled at the cell
  centers, respectively. Define the error vector by
  $\mathbf e_h^{(d)}=\mathbf p_h^{(d)}-\mathbf p_h^{\mathrm{ex}}$.
  We measure this error in the discrete $L_h^2$, maximum, and energy norms
  \[
    E_{2,h}^{(d)}=h\|\mathbf e_h^{(d)}\|_2,\qquad
    E_{\infty,h}^{(d)}=\|\mathbf e_h^{(d)}\|_\infty,\qquad
    E_{\mathbf{A},h}^{(d)}
    =\bigl((\mathbf e_h^{(d)})^T\mathbf{A}_h^{(d)}\mathbf e_h^{(d)}\bigr)^{1/2},
  \]
  where $\mathbf{A}_h^{(d)}$ is the stiffness matrix assembled using discretization
  $d$.

  \begin{table}[!htbp]
 \centering
 \caption{$L_h^2$ and $L^\infty$ errors and adjacent-grid orders
 for the two transmissibilities on an unfitted circular interface, $\kappa=10^4$.}
 \label{tab:transmissibility-circle-orders}
 \footnotesize\setlength{\tabcolsep}{3pt}
 \resizebox{\textwidth}{!}{\begin{tabular}{rcccccccc}\toprule
 & \multicolumn{4}{c}{$L_h^2$ error} & \multicolumn{4}{c}{$L^\infty$ error}\\
 \cmidrule(lr){2-5}\cmidrule(lr){6-9}
 $N_{\rm side}$ & Harmonic & Order & Interface-aware & Order & Harmonic & Order & Interface-aware & Order\\
 \midrule
$32$ & $6.0794\times10^{-4}$ & -- & $3.8349\times10^{-4}$ & -- & $2.0057\times10^{-3}$ & -- & $7.2906\times10^{-4}$ & -- \\
$64$ & $3.8227\times10^{-4}$ & $0.6693$ & $9.5918\times10^{-5}$ & $1.9993$ & $1.7034\times10^{-3}$ & $0.2357$ & $1.8703\times10^{-4}$ & $1.9628$ \\
$128$ & $9.0415\times10^{-5}$ & $2.0800$ & $2.3981\times10^{-5}$ & $1.9999$ & $1.0649\times10^{-3}$ & $0.6777$ & $4.7353\times10^{-5}$ & $1.9817$ \\
$256$ & $7.9383\times10^{-5}$ & $0.1877$ & $5.9944\times10^{-6}$ & $2.0002$ & $4.3301\times10^{-4}$ & $1.2982$ & $1.1912\times10^{-5}$ & $1.9910$ \\
$512$ & $2.3719\times10^{-5}$ & $1.7428$ & $1.4987\times10^{-6}$ & $1.9999$ & $2.1716\times10^{-4}$ & $0.9956$ & $2.9870\times10^{-6}$ & $1.9956$ \\ \bottomrule\end{tabular}}
\end{table}

  \begin{figure}[!htbp]
    \centering
    \includegraphics[width=\textwidth]{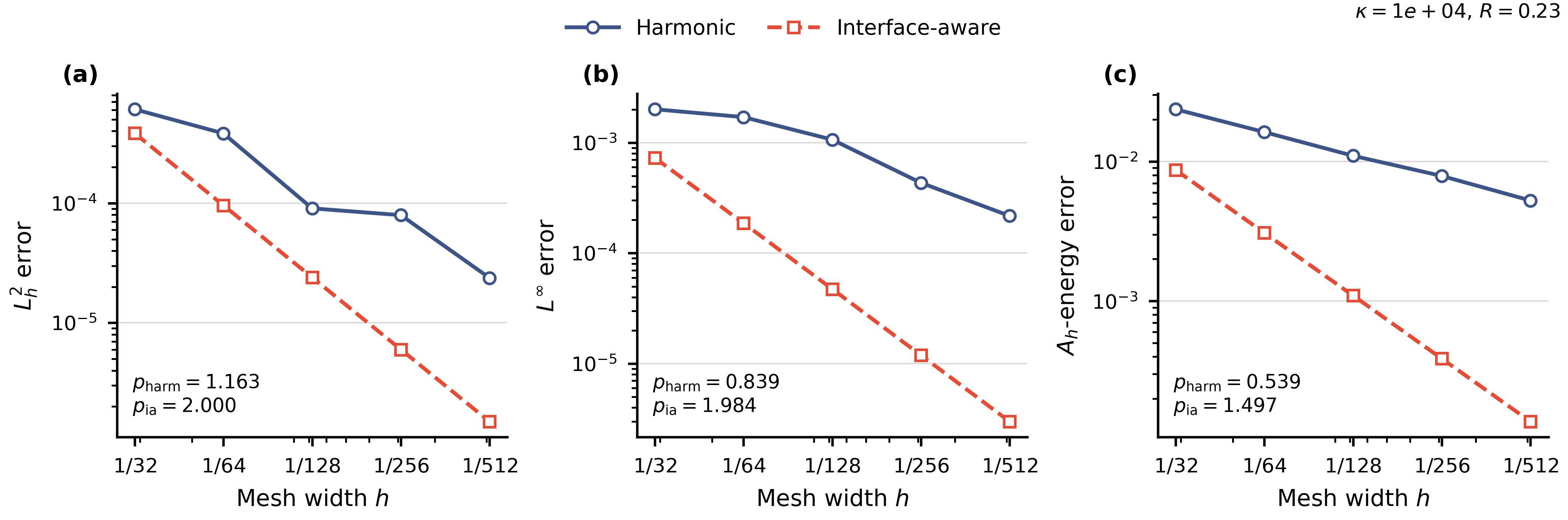}
    \caption{Mesh-convergence curves on an unfitted circular interface for
      $\kappa=10^4$. Panels~(a), (b), and (c) show the $L_h^2$, $L^\infty$, and
      discretization-specific $\mathbf{A}_h^{(d)}$ errors, respectively. The value $p$ in
    each panel is the global fitted order over the five grids}
    \label{fig:transmissibility-circle-convergence}
  \end{figure}

  \Cref{tab:transmissibility-circle-orders} and
  \Cref{fig:transmissibility-circle-convergence} show that the interface-aware
  scheme achieves nearly second-order convergence in both $L_h^2$ and
  $L^\infty$ for this test. In contrast, the harmonic scheme exhibits lower
  overall convergence rates and more variable adjacent-grid orders.
  These results support the use of the interface-aware transmissibility
  in the unfitted curved-interface experiments.

  \section{Implementation details for the additive Schwarz and smoothing effect}
  \label{sec:asm-implementation}

  This section details the ASM used in this work.
  We define the restriction and prolongation operators, describe the local
  problems and their boundary treatment, and explain how multiplicity weights
  combine the local corrections into a global update. We also outline the
  tensorized implementation and illustrate the smoothing effect.

  \subsection{Restriction and prolongation}

  Let
  $\mathcal{I}_k=\{i_{k,1},\ldots,i_{k,N_k}\}\subset\{1,\ldots,N_h\}$
  denote the set of global indices of the physical degrees of freedom in
  subdomain $D_k$, where
  $N_k=|\mathcal I_k|$. Define the Boolean restriction operator
  $\mathbf{R}_k\in\{0,1\}^{N_k\times N_h}$ by
  \begin{equation}
    (\mathbf{R}_k)_{\ell i}=
    \begin{cases}
      1, & i=i_{k,\ell},\\
      0, & \text{otherwise},
    \end{cases}
    \qquad 1\le \ell\le N_k,\quad 1\le i\le N_h.
    \label{eq:asm-restriction-entry}
  \end{equation}
  For any global vector $\mathbf r\in\mathbb R^{N_h}$,
  \begin{equation}
    (\mathbf{R}_k\mathbf{r})_{\ell}=r_{i_{k,\ell}},
    \qquad \ell=1,\ldots,N_k,
    \label{eq:asm-restriction-action}
  \end{equation}
  so $\mathbf{R}_k$ extracts the components associated with $D_k$. Its transpose
  $\mathbf{R}_k^T\in\{0,1\}^{N_h\times N_k}$ is the corresponding prolongation
  operator. For any local vector $\mathbf v_k\in\mathbb R^{N_k}$,
  \begin{equation}
    (\mathbf{R}_k^T\mathbf{v}_k)_i=
    \begin{cases}
      (\mathbf{v}_k)_{\ell}, & i=i_{k,\ell},\quad 1\le\ell\le N_k,\\
      0, & i\notin\mathcal{I}_k,
    \end{cases}
    \qquad 1\le i\le N_h,
    \label{eq:asm-prolongation-action}
  \end{equation}
  which extends the local vector by zero to the global degree-of-freedom space.
  Because the indices within a subdomain are distinct,
  $\mathbf{R}_k\mathbf{R}_k^T=\mathbf{I}_{N_k}$.

  \subsection{Local problems and boundary treatment}

  For the current pressure approximation $\mathbf p_h^{(\ell)}$, define the
  residual
  \[
    \mathbf{r}=\mathbf{f}_h-\mathbf{A}_h(\mu)\mathbf{p}_h^{(\ell)}
  \]
  and its restriction $\mathbf r_k=\mathbf{R}_k\mathbf r$ to $D_k$. The corresponding
  local stiffness matrix is
  \begin{equation}
    \mathbf{A}_k(\mu)
    =
    \mathbf{R}_k\,\mathbf{A}_h(\mu)\,\mathbf{R}_k^T
    \in\mathbb{R}^{N_k\times N_k}.
    \label{eq:local-asm-matrix}
  \end{equation}
  Thus, $\mathbf{A}_k(\mu)$ is the principal submatrix of $\mathbf{A}_h(\mu)$ associated with the
  subdomain degrees of freedom. This construction imposes homogeneous
  Dirichlet truncation at the artificial boundary
  $\partial D_k\setminus\partial\Omega$ in the algebraic sense: the local
  correction is extended by zero to all degrees of freedom outside
  $\mathcal I_k$. The portion of the subdomain boundary that coincides with
  the physical boundary inherits the original boundary discretization of
  $\mathbf{A}_h(\mu)$.

  Let $\mathbf{B}_k(\mu)$ denote the local solution operator. In the online experiments,
  $\mathbf{B}_k(\mu)=(\mathbf{A}_k(\mu)+\varepsilon_{\rm loc}\mathbf{I})^{-1}$, with
  $\varepsilon_{\rm loc}=10^{-6}$, is applied by exact block-tridiagonal
  Cholesky factorization. This is the implemented local diagonal shift,
  distinct from coarse Gram regularization and the relative linear residual tolerance.
  The local correction is
  \begin{equation}
    \mathbf{z}_k=\mathbf{B}_k(\mu)\,\mathbf{r}_k.
    \label{eq:local-asm-solve}
  \end{equation}
  The local solve is based on $\mathbf{A}_k(\mu)$ and therefore uses
  the same FVM face transmissibilities as the fine-grid discretization.

  \subsection{Weighted assembly and implementation}

  Because a global degree of freedom may receive contributions from several
  local corrections, we average these contributions using the multiplicity
  of the overlapping subdomain cover. Define the diagonal multiplicity matrix
  $\mathbf{W}\in\mathbb R^{N_h\times N_h}$ by
  \begin{equation}
    \mathbf{W}
    =
    \sum_{k=1}^K \mathbf{R}_k^T \mathbf{R}_k.
    \label{eq:partition-of-unity-weight}
  \end{equation}
  The diagonal entry $\mathbf{W}_{ii}$ is the number of subdomains containing global
  degree of freedom $i$. Every degree of freedom is covered at least once, so
  $\min_i \mathbf{W}_{ii}\ge1$ and $\mathbf{W}$ is invertible. Consequently,
  \[
    \sum_{k=1}^K \mathbf{W}^{-1}\mathbf{R}_k^T \mathbf{R}_k=\mathbf{I}_{N_h}.
  \]
  Thus, $\{\mathbf{W}^{-1}\mathbf{R}_k^T\mathbf{R}_k\}_{k=1}^K$ forms a discrete partition of unity, and
  $\mathbf{W}^{-1}\mathbf{R}_k^T$ is the corresponding weighted prolongation operator.

  The overlapping subdomain cover is determined by the window side length $w$
  and fixed stride $d$, both measured in grid cells, with $w>d$ so that adjacent
  windows overlap. If
  $K_{\mathrm{1d}}$ is the number of subdomains in each coordinate direction,
  then $K=K_{\mathrm{1d}}^2$, and each subdomain has at most $w^2$ physical
  degrees of freedom. In the tensorized implementation, windows that extend
  beyond the physical domain are zero-padded to a common $w\times w$ shape.
  A mask identifies entries corresponding to physical degrees of freedom.
  Padded tensor entries do not add restricted global unknowns to $\mathbf{A}_k$.

  The local corrections are extended by zero, summed, and weighted to give
  the global correction
  \begin{equation}
    \mathbf{z}
    =
    \mathbf{W}^{-1}\sum_{k=1}^K \mathbf{R}_k^T\mathbf{z}_k.
    \label{eq:averaged-asm-correction}
  \end{equation}
  For the residual
  $\mathbf r=\mathbf f_h-\mathbf{A}_h(\mu)\mathbf p_h^{(\ell)}$, this gives
  $\mathbf z=\mathbf{M}_{\mathrm{ASM}}(\mu)^{-1}\mathbf r$. The damped ASM update is
  \begin{equation}
    \mathbf{p}_h^{(\ell+1)}=\mathbf{p}_h^{(\ell)}
    +\omega \mathbf{M}_{\mathrm{ASM}}(\mu)^{-1}
    \bigl(\mathbf{f}_h-\mathbf{A}_h(\mu)\mathbf{p}_h^{(\ell)}\bigr),
    \label{eq:asm-iter}
  \end{equation}
  where $\omega$ is the damping parameter and $\mathbf{M}_{\mathrm{ASM}}(\mu)^{-1}$
  is the operator in \eqref{eq:asm-preconditioner}, with its dependence on
  $\mu$ shown explicitly.

  \subsection{ASM smoothing effect}
  \label{sec:asm-smoothing-visualization}

  To illustrate the smoothing action described in \Cref{sec:two-grid}, we
  consider the homogeneous error equation $\mathbf{A}_h\mathbf e=\mathbf0$ for a fixed
  parameter sample, whose unique solution is $\mathbf e=\mathbf0$. The experiment uses
  $N_{\mathrm{side}}=63$, a random high-conductivity inclusion with contrast
  $10^4$, and the interface-aware transmissibility in
  \eqref{eq:interface-aware-transmissibility}. The grid is decomposed into
  a $10\times10$ array of overlapping subdomains with target overlap
  ratio $0.3$ (patch side 10 and stride 7 on this grid). Each local residual-correction problem is solved
  directly using a Cholesky factorization. Starting from a standard Gaussian
  random initial error $\mathbf e^{(0)}$, we apply ASM iteration with
  $\omega=0.8$.

\begin{figure}[!htb]
    \centering
    \includegraphics[width=0.98\textwidth]{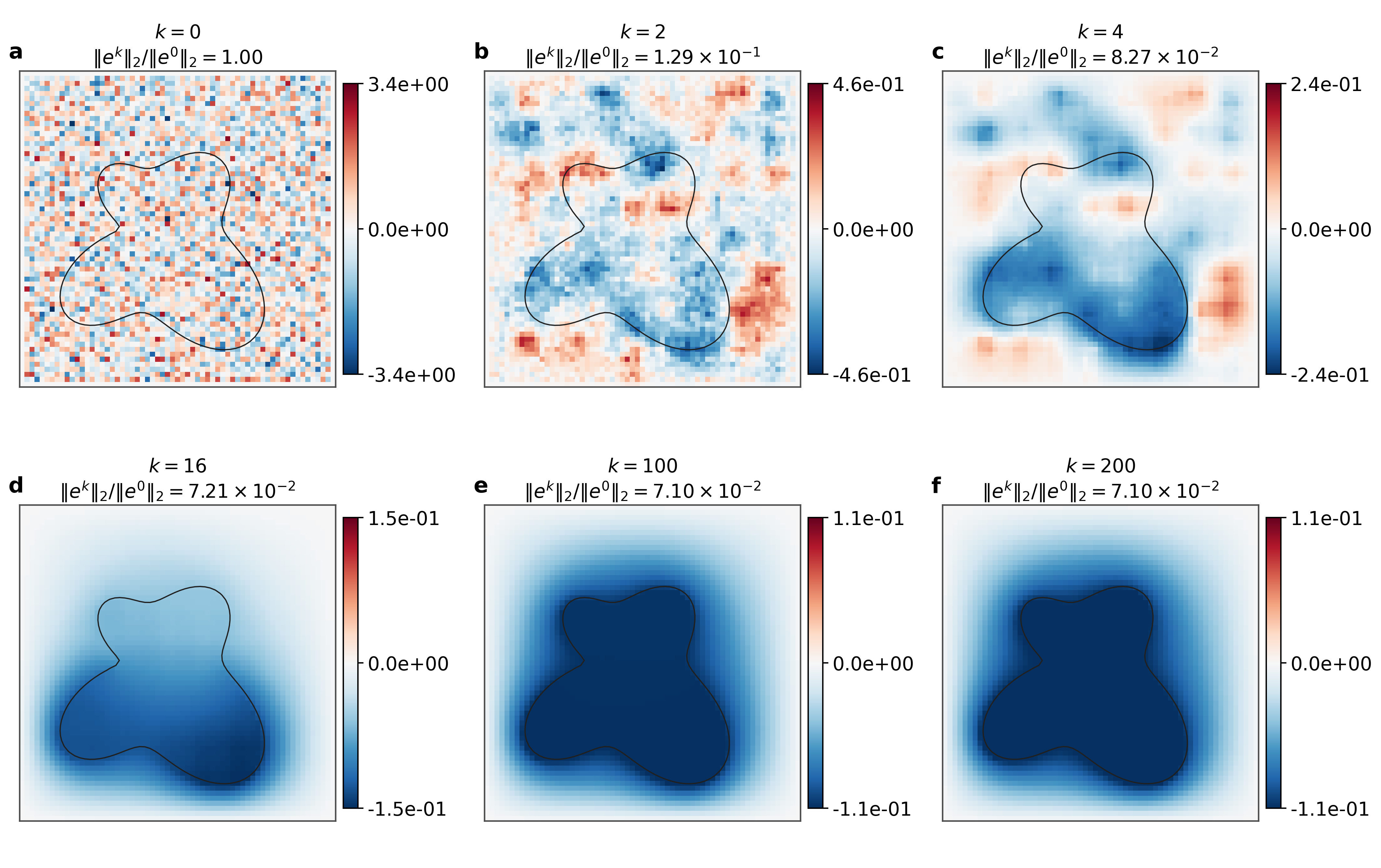}
    \caption{Evolution of a random initial error under ASM after
    $k=0,2,4,16,100,200$ updates}
    \label{fig:asm-random-error-smoothing}
  \end{figure}

  \Cref{fig:asm-random-error-smoothing} shows the error fields after
  $k=0,2,4,16,100,200$ updates. As the iteration proceeds, the initially
  oscillatory and fragmented pattern gives way to larger, more spatially
  coherent structures that persist after 200 updates. This behavior
  illustrates the smoothing effect of ASM: local oscillations are damped,
  while larger-scale error components decay more slowly, motivating the use
  of a coarse-grid correction.

\section{Ablation of snapshot smoothing and sample count}
\label{sec:snapshot-strategy-ablation}

We investigate the effects of snapshot smoothing and sample count
on the performance of POD-based coarse spaces. All experiments use
exact block-tridiagonal Cholesky ASM with TG(1,1) on a grid with
$N_{\rm side}=63$, coefficient contrast $10^4$, damping parameter
$\omega=0.8$, and coarse-space dimension $n=15$.
Two parameter sweeps are conducted using the same eight held-out
coefficient fields. First, we vary the number of snapshot smoothing
steps $\nu_{\rm snap}\in\{5,10,20,30,50,75\}$ with $q=20$.
Second, we vary the number of Gaussian random initial errors
$q\in\{15,20,30,40,50,75\}$ with $\nu_{\rm snap}=50$.
The random initial error sets are nested, with each smaller set
drawn from the first $q$ samples of a common set of 75 realizations.
The two sweeps share the configuration $(\nu_{\rm snap},q)=(50,20)$,
resulting in 176 solver runs across both methods.
All solves use a relative residual stopping criterion and a maximum
of 200 updates, with the final residuals independently verified.

\begin{table}[!htb]
  \centering
  \caption{Effects of snapshot smoothing and sample count on
  POD--TG and POD--TG--FGMRES convergence. Iteration counts are
  reported as mean $\pm$ population standard deviation over
  converged solves. Runs exceeding the 200-update limit are
  classified as failures. All entries correspond to 8/8
  converged solves unless otherwise indicated:
  $\dagger$ denotes 0/8 and $\ddagger$ denotes 4/8.}
  \label{tab:snapshot-strategy-ablation}
  \footnotesize
  \setlength{\tabcolsep}{4pt}
  \begin{tabular}{@{}rcc@{\qquad}rcc@{}}
    \toprule
    \multicolumn{3}{c}{Fixed $q=20$}
    & \multicolumn{3}{c}{Fixed $\nu_{\rm snap}=50$}\\
    \cmidrule(lr){1-3}\cmidrule(lr){4-6}
    $\nu_{\rm snap}$ & POD--TG & POD--TG--FGMRES
    & $q$ & POD--TG & POD--TG--FGMRES\\
    \midrule
    5  & \textemdash$^{\dagger}$ & $21.62\pm0.86$
       & 15 & $26.25\pm2.82$ & $10.62\pm0.48$ \\
    10 & $87.75\pm11.01^{\ddagger}$ & $17.00\pm1.00$
       & 20 & $25.88\pm2.57$ & $10.50\pm0.50$ \\
    20 & $36.38\pm3.84$ & $11.88\pm0.60$
       & 30 & $24.62\pm2.78$ & $10.25\pm0.43$ \\
    30 & $29.38\pm3.57$ & $11.00\pm0.00$
       & 40 & $24.88\pm2.42$ & $10.38\pm0.48$ \\
    50 & $25.88\pm2.57$ & $10.50\pm0.50$
       & 50 & $25.12\pm2.71$ & $10.38\pm0.48$ \\
    75 & $23.88\pm2.03$ & $10.12\pm0.33$
       & 75 & $24.25\pm2.68$ & $10.38\pm0.48$ \\
    \bottomrule
  \end{tabular}
\end{table}

\Cref{tab:snapshot-strategy-ablation} shows that snapshot smoothing
has a pronounced effect on the convergence of both solvers.
With only 5 or 10 smoothing steps, standalone POD--TG exhibits
poor convergence, whereas POD--TG--FGMRES converges on all
eight fields but requires more iterations. Increasing
$\nu_{\rm snap}$ from 10 to 50 substantially improves
convergence, reducing the mean POD--TG--FGMRES iteration count
from 17.00 to 10.50. Further increasing $\nu_{\rm snap}$ to 75
provides only modest additional improvement.
This behavior is consistent with the role of smoothing in
attenuating oscillatory error components, allowing the
POD basis to capture the slowly decaying error modes
that benefit most from coarse-grid correction.

In contrast, increasing the number of random initial errors
has a relatively small effect once sufficient smoothing
is applied. At $\nu_{\rm snap}=50$, increasing $q$ from 15
to 75 reduces the mean standalone POD--TG iteration count
from 26.25 to 24.25, while the POD--TG--FGMRES count remains
close to 10 throughout the sweep. The small, nonmonotone
variations indicate that the dominant error directions
are already well represented by a moderate number of
smoothed snapshots. These results motivate the choice
$(\nu_{\rm snap},q)=(50,20)$ for constructing the reference
coarse spaces.

  \section{Continuous-coordinate implementation}
  \label{sec:network-implementation}

  This section describes the implementation of the continuous-coordinate
  network for generating  coarse spaces. We detail the
  branch encoder, coordinate trunk, and coefficient head, then explain the
  boundary treatment and basis representations used during training and in
  the online Galerkin correction. We also clarify the full-rank requirement
  and the prepared basis representation used by the online solver.

  \subsection{Branch encoder}
  The network first resamples $[\log_{10}\lambda,\phi]$ by bilinear
  interpolation to obtain a $64\times64$ two-channel image. A $5\times5$
  convolution increases the channel count from 2 to 32 and is followed by
  three stages with 32, 64, and 128 channels, respectively.
  Each stage contains two residual blocks, each with two $3\times3$
  convolutions. The first block in each of the last two stages downsamples
  with stride two. The convolutional branch uses GroupNorm with eight groups, GELU,
  and squeeze-and-excitation with a reduction ratio of four, together with a
  $1\times1$ residual projection whenever the channel count or resolution
  changes. Global average pooling and a final $128\to128\to128$ MLP produce
  the parameter representation $\mathbf{z}_\mu$ of dimension $d_z=128$.

  \subsection{Coordinate trunk and coefficient head}
  For each cell center $\boldsymbol x_i=(x_i,y_i)$, the trunk input is
  \begin{equation}
    \begin{aligned}
      \boldsymbol{\chi}_i
      &=\Bigl[x_i,\,y_i,\,\phi(\boldsymbol x_i),\,
        \boldsymbol{\gamma}_0(x_i,y_i),\,\ldots,\,
        \boldsymbol{\gamma}_5(x_i,y_i)\Bigr]\in\mathbb R^{27},\\[3pt]
      \boldsymbol{\gamma}_r(x,y)
      &=\bigl[\sin(\omega_r x),\,\cos(\omega_r x),\,
        \sin(\omega_r y),\,\cos(\omega_r y)\bigr],\\[2pt]
      \omega_r&=2^r\pi,\qquad r=0,\ldots,5.
    \end{aligned}
    \label{eq:trunk-point-features}
  \end{equation}
  The resulting 27-dimensional input includes coordinate Fourier features
  to enhance the representation of high-frequency content
  \cite{tancik2020fourier}. The trunk has width 128 and contains three residual
  MLP blocks with internal width 256, GELU activations, and LayerNorm.
  Each block is modulated by $\mathbf{z}_\mu$ through FiLM \cite{perez2018film}
  \begin{equation}
    \mathbf{u}_{\mathrm{FiLM}}
    =\mathbf{u}\odot\bigl(1+0.1\,\mathbf a(\mathbf{z}_\mu)\bigr)+0.1\,\mathbf b(\mathbf{z}_\mu),
    \label{eq:film-implementation}
  \end{equation}
  where $\odot$ denotes elementwise multiplication. The nonlinear increment
  of each block is further modulated by the gate
  $\operatorname{sigmoid}(\mathbf{W}_g\mathbf{z}_\mu)$. The trunk output dimension is
  $d_{\mathrm{dict}}=64$. The coefficient head concatenates a learned
  32-dimensional column embedding with $\mathbf{z}_\mu$ and uses a $160\to128$ input
  layer, two $128\to256\to128$ residual blocks, and a columnwise
  $128\to64$ output layer to generate $\mathbf C_\theta(\mu)$. For the static
  experiment with $n=15$, the network has $1{,}432{,}208$ trainable parameters.

  \subsection{Basis representation and the online coarse problem}
  For the unit-square problems with homogeneous Dirichlet conditions on the
  entire boundary, the implementation uses the boundary envelope
  \begin{equation}
    \begin{aligned}
      b_{\partial}(x,y)&=16x(1-x)y(1-y),\\
      \mathbf D_{\partial,h}
      &=\operatorname{diag}\bigl(b_{\partial}(x_i,y_i)\bigr)_{i=1}^{N_h}.
    \end{aligned}
    \label{eq:boundary-envelope-formula}
  \end{equation}
  This is the formula evaluated by the code and applied once to the
  assembled tentative basis in \eqref{eq:boundary-envelope} during the
  network forward pass. It vanishes on all four sides of $[0,1]^2$ and is
  sampled at the target grid's cell centers. For the mixed-boundary
  moving-interface problem, the envelope is omitted and
  $\mathbf D_{\partial,h}=\mathbf I$.
  In both cases, $\mathbf{V}_{\theta,h}^{\rm raw}$ defines the coarse space
  used in the loss and the Galerkin correction.

  During training, the Cholesky change of basis in
  \eqref{eq:training-energy-orthogonalization} yields an $\mathbf{A}_h$-orthonormal
  basis for stable evaluation of the projector without changing the predicted
  space. The ideal online Galerkin correction depends only on the raw full-rank
  span and is expressed using $\mathbf Q^\star$ in
  \eqref{eq:online-raw-direct-galerkin}. The online coarse solve factors
  the Gram matrix of the prepared full-rank representative plus
  $\varepsilon_{\rm coarse}\mathbf{I}$. When used, regularized Cholesky preparation
  need not yield exact energy orthonormality.
  The static and moving-interface preparations use $10^{-10}$.
  The Forchheimer Picard preparations use $10^{-8}$.
  The regularized correction is not literally an
  orthogonal projector. Its deviation is bounded by
  \eqref{eq:coarse-regularization-gap}. A sample that fails the full-rank check is
  counted as a coarse-space construction failure rather than repaired by
  column truncation, which would alter the predicted space.
  Before regularization, the two basis representations give the same
  projector. A scalar diagonal shift in an arbitrary raw basis is not
  basis invariant. The bound must use the Gram of the basis actually
  prepared. The online setup cost includes this preparation work.

\section{Finite-precision solves and independent residual verification}
\label{sec:finite-precision-verification}

\subsection{Fixed-map FGMRES and acceptance}

The static problems are solved using right-preconditioned
FGMRES(50), with preconditioned vectors
$\mathbf{z}_j=\mathbf{B}\mathbf{v}_j$ retained explicitly
and two-pass Euclidean orthogonalization.
The discrete operator, local Cholesky factors, and
Galerkin coarse factors remain fixed throughout each
linear solve. Residual estimates obtained from the
Arnoldi least-squares problem are used to identify
candidate stopping points, while final acceptance
is determined by the explicitly evaluated residual
\[
  \rho_{\rm lin}
  =
  \frac{\|\mathbf{f}-\mathbf{A}_{\rm st}\mathbf{p}\|_2}
       {\|\mathbf{f}\|_2}
  \le \mathrm{tol}_{\rm lin}=10^{-6},
\]
where $\mathbf{A}_{\rm st}$ is the stored float64
discrete matrix. For the static experiments, accepted
solutions are additionally verified using an
independently assembled CPU sparse-matrix action.

\subsection{Retained pressure components for Forchheimer flow}

To improve residual accuracy in the high-contrast
Forchheimer simulations, pressure is represented
as the sum of two float64 components,
\[
  \mathbf{p}
  =
  \mathbf{p}_{\rm high}
  +
  \mathbf{p}_{\rm low}.
\]
Compensated accumulation preserves low-order
contributions during solution updates and FGMRES
restarts. Residuals are evaluated directly from
both pressure components using the original
stored five-point coefficients, including their
diagonal roundoff contributions. This representation
avoids the loss of low-order information associated
with forming a single rounded float64 pressure vector.
Arnoldi vectors and TG preconditioner applications
remain in float64, while network inference uses
the float32 arithmetic of the trained checkpoint.
The linear solves employ a relative residual tolerance
of $5\times10^{-7}$, followed by independent
verification against a tolerance of $10^{-6}$.

Both pressure components are retained throughout
the Picard iteration, including warm starts and
the evaluation of pressure gradients, permeability,
and mobility. At each Picard step, the discrete
operator and TG preconditioner are constructed
from the current mobility and held fixed during
the corresponding linear solve. Nonlinear convergence
is assessed using both the relative split-pressure
change and the updated-mobility residual defined
in \Cref{sec:forchheimer-channelized}, with a
tolerance of $10^{-6}$ for each criterion.
A pressure update is accepted only after the
associated linear residual check is satisfied.
The maximum budgets are 300 linear updates and
100 Picard iterations. Full-state checkpoints
preserve both pressure components, configuration
hashes, and complete iteration histories.

\subsection{Independent residual verification}

The independent Forchheimer validation covers all eight paired channel
fields, both $\alpha=0$ and $\alpha=10$, and all three coarse-space methods
in \Cref{tab:forchheimer-picard-ablation}. The 48 nonlinear trajectories
comprise 1201 accepted Picard steps. At each accepted step, the stored
data include both pressure components, the frozen and updated nonlinear
stencils, mobility fields, and coarse Galerkin matrices. An independent CPU evaluator
uses double-double arithmetic for products and
accumulations, supplemented by selected-row
crosschecks with 100-digit Decimal arithmetic.
These evaluations verify the linear and nonlinear
residuals, relative pressure changes, and
consistency of the mobility updates against
the stored discrete operators.

The independent checks provide additional
numerical validation of the reported solver
trajectories and stopping criteria.
Together with the retained pressure representation,
they improve the reliability of residual evaluation
in high-contrast regimes while preserving the
fixed-operator setting underlying the perturbation
analysis in \Cref{sec:convergence}.

\end{appendices}

\bibliography{ref}

@article{meng2026overview,
  title={An overview of machine learning-enhanced iterative methods for systems of linear and nonlinear equations},
  author={Meng, Yuhuang and Zhao, Jing and Heinlein, Alexander},
  journal={arXiv preprint arXiv:2610.07211},
  year={2026},
  url = {https://arxiv.org/abs/2610.07211}
}

@inproceedings{adam,
  author    = {Diederik P. Kingma and
               Jimmy Ba},
  editor    = {Yoshua Bengio and
               Yann LeCun},
  title     = {Adam: {A} Method for Stochastic Optimization},
  booktitle = {3rd International Conference on Learning Representations, {ICLR} 2015,
               San Diego, CA, USA, May 7-9, 2015, Conference Track Proceedings},
  year      = {2015},
  url       = {http://arxiv.org/abs/1412.6980},
  bibsource = {dblp computer science bibliography, https://dblp.org}
}

@incollection{bennequintechniques,
  title     = {Techniques to study spectral properties of the {Robin--Schwarz} operator at the continuous level},
  author    = {Bennequin, Daniel and Cuvelier, Fran{\c{c}}ois and Gander, Martin J. and Halpern, Laurence},
  booktitle = {Domain Decomposition Methods in Science and Engineering XXVIII},
  series    = {Lecture Notes in Computational Science and Engineering},
  volume    = {155},
  publisher = {Springer},
  address   = {Cham},
  year      = {forthcoming},
  url       = {https://www.ddm.org/DD28/proceedings/DD28_375.pdf}
}

@article{brandt2011bootstrap,
  title     = {Bootstrap {AMG}},
  author    = {Brandt, Achi and Brannick, James and Kahl, Karsten and Livshits, Irene},
  journal   = {SIAM Journal on Scientific Computing},
  volume    = {33},
  number    = {2},
  pages     = {612--632},
  year      = {2011},
  publisher = {SIAM},
  doi = {10.1137/090752973}
}

@article{brezina2005adaptive,
  title     = {Adaptive Smoothed Aggregation ($\alpha${SA}) Multigrid},
  author    = {Brezina, Marian and Falgout, Robert and MacLachlan, Scott and Manteuffel, Tom and McCormick, Steve and Ruge, John},
  journal   = {SIAM Review},
  volume    = {47},
  number    = {2},
  pages     = {317--346},
  year      = {2005},
  publisher = {SIAM},
  doi = {10.1137/050626272}
}

@article{cai1999restricted,
  author  = {Cai, Xiao-Chuan and Sarkis, Marcus},
  title   = {A Restricted Additive {Schwarz} Preconditioner for General Sparse Linear Systems},
  journal = {SIAM Journal on Scientific Computing},
  volume  = {21},
  number  = {2},
  pages   = {792--797},
  year    = {1999},
  doi     = {10.1137/S106482759732678X}
}

@article{caldana2024deep,
  title     = {A deep learning algorithm to accelerate algebraic multigrid methods in finite element solvers of {3D} elliptic {PDEs}},
  author    = {Caldana, Matteo and Antonietti, Paola F. and Ded{\`e}, Luca},
  journal   = {Computers \& Mathematics with Applications},
  volume    = {167},
  pages     = {217--231},
  year      = {2024},
  publisher = {Elsevier},
  doi = {10.1016/j.camwa.2024.05.013}
}

@article{chen2022meta,
  author    = {Chen, Yuyan and Dong, Bin and Xu, Jinchao},
  journal   = {Journal of Computational Physics},
  title     = {{Meta-MgNet}: Meta multigrid networks for solving parameterized partial differential equations},
  year      = {2022},
  pages     = {110996},
  volume    = {455},
  publisher = {Academic Press},
  ranking   = {rank4},
  doi = {10.1016/j.jcp.2022.110996}
}

@inproceedings{chen2025graph,
  title     = {Graph Neural Preconditioners for Iterative Solutions of Sparse Linear Systems},
  author    = {Jie Chen},
  booktitle = {The Thirteenth International Conference on Learning Representations},
  year      = {2025},
  url = {https://openreview.net/forum?id=Tkkrm3pA35}
}

@article{chow2006aggregation,
  title     = {An aggregation multilevel method using smooth error vectors},
  author    = {Chow, Edmond},
  journal   = {SIAM Journal on Scientific Computing},
  volume    = {27},
  number    = {5},
  pages     = {1727--1741},
  year      = {2006},
  publisher = {SIAM},
  doi = {10.1137/040608192}
}

@article{CSIAM-AM-4-13,
  author  = {Cui, Chen and Jiang, Kai and Shu, Shi},
  title   = {Solving Time-Dependent Parametric {PDEs} by Multiclass Classification-Based Reduced Order Model},
  journal = {CSIAM Transactions on Applied Mathematics},
  year    = {2023},
  volume  = {4},
  number  = {1},
  pages   = {13--40},
  issn    = {2708-0579},
  doi = {10.4208/csiam-am.SO-2021-0042}
}

@article{cui2022fourier,
  author  = {Cui, Chen and Jiang, Kai and Liu, Yun and Shu, Shi},
  journal = {Mathematics},
  title   = {{Fourier} Neural Solver for Large Sparse Linear Algebraic Systems},
  year    = {2022},
  issn    = {2227-7390},
  number  = {21},
  volume  = {10},
  pages   = {4014},
  ranking = {rank5},
  doi = {10.3390/math10214014}
}

@article{cui2024hybrid,
  title     = {A Hybrid Iterative Neural Solver Based on Spectral Analysis for Parametric {PDEs}},
  author    = {Cui, Chen and Jiang, Kai and Liu, Yun and Shu, Shi},
  journal   = {Journal of Computational Physics},
  volume    = {538},
  pages     = {114165},
  year      = {2025},
  doi       = {10.1016/j.jcp.2025.114165},
  publisher = {Elsevier}
}

@article{cui2025neural,
  title     = {A neural multigrid solver for {Helmholtz} equations with high wavenumber and heterogeneous media},
  author    = {Cui, Chen and Jiang, Kai and Shu, Shi},
  journal   = {SIAM Journal on Scientific Computing},
  volume    = {47},
  number    = {3},
  pages     = {C655--C679},
  year      = {2025},
  publisher = {Society for Industrial and Applied Mathematics University City, Philadelphia},
  doi = {10.1137/24M1654397}
}

@incollection{cuvelier2023fundamental,
  address   = {Cham},
  title     = {Fundamental coarse space components for {Schwarz} methods with crosspoints},
  author    = {Cuvelier, Fran{\c{c}}ois and Gander, Martin J and Halpern, Laurence},
  booktitle = {Domain Decomposition Methods in Science and Engineering XXVI},
  pages     = {41--52},
  year      = {2022},
  publisher = {Springer},
  series = {Lecture Notes in Computational Science and Engineering},
  volume = {145},
  doi = {10.1007/978-3-030-95025-5_4}
}

@article{dimola2026numerical,
  title     = {Numerical solution of mixed-dimensional {PDEs} using a neural preconditioner},
  author    = {Dimola, Nunzio and Franco, Nicola Rares and Zunino, Paolo},
  journal   = {Computers \& Mathematics with Applications},
  volume    = {206},
  pages     = {58--79},
  year      = {2026},
  publisher = {Elsevier},
  doi = {10.1016/j.camwa.2025.12.013}
}

@incollection{dohrmann2008family,
  address   = {Berlin, Heidelberg},
  title     = {A family of energy minimizing coarse spaces for overlapping {Schwarz} preconditioners},
  author    = {Dohrmann, Clark R and Klawonn, Axel and Widlund, Olof B},
  booktitle = {Domain decomposition methods in science and engineering XVII},
  pages     = {247--254},
  year      = {2008},
  publisher = {Springer},
  doi = {10.1007/978-3-540-75199-1_28}
}

@article{efendiev2013generalized,
  title     = {Generalized multiscale finite element methods ({GMsFEM})},
  author    = {Efendiev, Yalchin and Galvis, Juan and Hou, Thomas Y},
  journal   = {Journal of Computational Physics},
  volume    = {251},
  pages     = {116--135},
  year      = {2013},
  publisher = {Elsevier},
  doi = {10.1016/j.jcp.2013.04.045}
}

@inproceedings{fanaskov2026deep,
  title     = {Deep Learning for Subspace Regression},
  author    = {Fanaskov, Vladimir and Trifonov, Vladislav and Rudikov, Alexander and Muravleva, Ekaterina and Oseledets, Ivan},
  booktitle = {International Conference on Learning Representations},
  year      = {2026},
  url       = {https://openreview.net/forum?id=HF60Lu1Maj}
}

@article{fresca2022pod,
  title     = {{POD-DL-ROM}: enhancing deep learning-based reduced order models for nonlinear parametrized {PDEs} by proper orthogonal decomposition},
  author    = {Fresca, Stefania and Manzoni, Andrea},
  journal   = {Computer Methods in Applied Mechanics and Engineering},
  volume    = {388},
  pages     = {114181},
  year      = {2022},
  publisher = {Elsevier},
  doi = {10.1016/j.cma.2021.114181}
}

@article{fu2024efficient,
  title     = {An efficient multiscale preconditioner for large-scale highly heterogeneous flow},
  author    = {Fu, Shubin and Chung, Eric and Zhao, Lina},
  journal   = {SIAM Journal on Scientific Computing},
  volume    = {46},
  number    = {2},
  pages     = {S352--S377},
  year      = {2024},
  publisher = {SIAM},
  doi = {10.1137/22M1502859}
}

@article{gander2025hierarchical,
  title   = {Hierarchical Coarse Basis by Randomised {SVD}: the {Helmholtz} Problem},
  author  = {Gander, Martin J and Jiang, Yao-Lin and Zhang, Hui},
  journal = {arXiv preprint arXiv:2512.12538},
  year    = {2025},
  url = {https://arxiv.org/abs/2512.12538}
}

@article{gelu,
  title   = {Gaussian error linear units (gelus)},
  author  = {Hendrycks, Dan and Gimpel, Kevin},
  journal = {arXiv preprint arXiv:1606.08415},
  year    = {2016}
}

@article{girault2008darcyforchheimer,
  author  = {Girault, Vivette and Wheeler, Mary F.},
  title   = {Numerical Discretization of a {Darcy--Forchheimer} Model},
  journal = {Numerische Mathematik},
  volume  = {110},
  number  = {2},
  pages   = {161--198},
  year    = {2008},
  doi     = {10.1007/s00211-008-0157-7}
}

@article{gottlieb2001strong,
  title     = {Strong stability-preserving high-order time discretization methods},
  author    = {Gottlieb, Sigal and Shu, Chi-Wang and Tadmor, Eitan},
  journal   = {SIAM Review},
  volume    = {43},
  number    = {1},
  pages     = {89--112},
  year      = {2001},
  publisher = {SIAM},
  doi = {10.1137/S003614450036757X}
}

@inproceedings{greenfeld2019learning,
  title        = {Learning to Optimize Multigrid {PDE} Solvers},
  author       = {Greenfeld, Daniel and Galun, Meirav and Basri, Ronen and Yavneh, Irad and Kimmel, Ron},
  booktitle    = {International Conference on Machine Learning},
  pages        = {2415--2423},
  year         = {2019},
  organization = {PMLR},
  volume = {97},
  series = {Proceedings of Machine Learning Research},
  url = {https://proceedings.mlr.press/v97/greenfeld19a.html}
}

@book{hackbusch1985multi,
  address   = {Berlin, Heidelberg},
  title     = {Multi-grid methods and applications},
  author    = {Hackbusch, Wolfgang},
  year      = {1985},
  publisher = {Springer},
  doi = {10.1007/978-3-662-02427-0}
}

@book{hesthaven2016certified,
  author    = {Hesthaven, Jan S. and Rozza, Gianluigi and Stamm, Benjamin},
  title     = {Certified Reduced Basis Methods for Parametrized Partial Differential Equations},
  series    = {SpringerBriefs in Mathematics},
  publisher = {Springer},
  address   = {Cham},
  year      = {2016},
  doi       = {10.1007/978-3-319-22470-1}
}

@article{hu2025hybrid,
  title   = {A hybrid iterative method based on {MIONet} for {PDEs}: Theory and numerical examples},
  author  = {Hu, Jun and Jin, Pengzhan},
  journal = {Mathematics of Computation},
  volume  = {95},
  number  = {359},
  pages   = {1327--1359},
  year    = {2026},
  doi     = {10.1090/mcom/4086}
}

@article{huang2022learning,
  title   = {Learning Optimal Multigrid Smoothers via Neural Networks},
  author  = {Huang, Ru and Li, Ruipeng and Xi, Yuanzhe},
  journal = {SIAM Journal on Scientific Computing},
  volume  = {45},
  number  = {3},
  pages   = {S199--S225},
  year    = {2023},
  doi     = {10.1137/21M1430030}
}

@article{huang2024reducing,
  title     = {Reducing Operator Complexity of {Galerkin} Coarse-grid Operators with Machine Learning},
  author    = {Huang, Ru and Chang, Kai and He, Huan and Li, Ruipeng and Xi, Yuanzhe},
  journal   = {SIAM Journal on Scientific Computing},
  volume    = {46},
  number    = {5},
  pages     = {S296--S316},
  year      = {2024},
  publisher = {Society for Industrial and Applied Mathematics University City, Philadelphia},
  doi = {10.1137/23M1583533}
}

@article{katrutsa2020black,
  author    = {Katrutsa, Alexandr and Daulbaev, Talgat and Oseledets, Ivan},
  journal   = {Journal of Computational and Applied Mathematics},
  title     = {Black-box learning of multigrid parameters},
  year      = {2020},
  pages     = {112524},
  volume    = {368},
  publisher = {Elsevier},
  ranking   = {rank4},
  doi = {10.1016/j.cam.2019.112524}
}

@article{kopanivcakova2025deeponet,
  title     = {{DeepONet} Based Preconditioning Strategies for Solving Parametric Linear Systems of Equations},
  author    = {Kopani{\v{c}}{\'a}kov{\'a}, Alena and Karniadakis, George Em},
  journal   = {SIAM Journal on Scientific Computing},
  volume    = {47},
  number    = {1},
  pages     = {C151--C181},
  year      = {2025},
  publisher = {SIAM},
  doi = {10.1137/24M162861X}
}

@article{liu2026approximating,
  title     = {Approximating generalized multiscale coarse spaces with application in high-contrast Darcy flow problems},
  author    = {Liu, Yucheng and Fu, Shubin and Zhou, Yingjie and Ye, Changqing and Chung, Eric T},
  journal   = {Journal of Computational Physics},
  pages     = {114950},
  year      = {2026},
  publisher = {Elsevier},
  volume = {560},
  doi = {10.1016/j.jcp.2026.114950}
}

@inproceedings{luz2020learning,
  author       = {Luz, Ilay and Galun, Meirav and Maron, Haggai and Basri, Ronen and Yavneh, Irad},
  booktitle    = {International Conference on Machine Learning},
  title        = {Learning algebraic multigrid using graph neural networks},
  year         = {2020},
  organization = {PMLR},
  pages        = {6489--6499},
  ranking      = {rank5},
  volume = {119},
  series = {Proceedings of Machine Learning Research},
  url = {https://proceedings.mlr.press/v119/luz20a.html}
}

@article{mao2026accurate,
  title     = {Accurate and scalable deep {Maxwell} solvers},
  author    = {Mao, Chenkai and Fan, Jonathan A},
  journal   = {Proceedings of the National Academy of Sciences},
  volume    = {123},
  number    = {18},
  pages     = {e2530330123},
  year      = {2026},
  publisher = {National Academy of Sciences},
  doi = {10.1073/pnas.2530330123}
}

@article{moore2025graph,
  title     = {Graph neural networks and applied linear algebra},
  author    = {Moore, Nicholas S and Cyr, Eric C and Ohm, Peter and Siefert, Christopher M and Tuminaro, Raymond S},
  journal   = {SIAM Review},
  volume    = {67},
  number    = {1},
  pages     = {141--175},
  year      = {2025},
  publisher = {SIAM},
  doi = {10.1137/23M1609786}
}

@article{osher1988fronts,
  title     = {Fronts Propagating with Curvature-Dependent Speed: Algorithms Based on {Hamilton--Jacobi} Formulations},
  author    = {Osher, Stanley and Sethian, James A},
  journal   = {Journal of Computational Physics},
  volume    = {79},
  number    = {1},
  pages     = {12--49},
  year      = {1988},
  publisher = {Elsevier},
  doi = {10.1016/0021-9991(88)90002-2}
}

@inproceedings{perez2018film,
  author    = {Perez, Ethan and Strub, Florian and de Vries, Harm and Dumoulin, Vincent and Courville, Aaron},
  title     = {{FiLM}: Visual Reasoning with a General Conditioning Layer},
  booktitle = {Proceedings of the AAAI Conference on Artificial Intelligence},
  volume    = {32},
  pages     = {3942--3951},
  year      = {2018},
  doi       = {10.1609/aaai.v32i1.11671}
}

@book{quarteroni2016reduced,
  author    = {Quarteroni, Alfio and Manzoni, Andrea and Negri, Federico},
  title     = {Reduced Basis Methods for Partial Differential Equations: An Introduction},
  publisher = {Springer},
  address   = {Cham},
  year      = {2016},
  doi       = {10.1007/978-3-319-15431-2}
}

@article{saad1993flexible,
  title     = {A flexible inner-outer preconditioned {GMRES} algorithm},
  author    = {Saad, Youcef},
  journal   = {SIAM Journal on Scientific Computing},
  volume    = {14},
  number    = {2},
  pages     = {461--469},
  year      = {1993},
  publisher = {SIAM},
  doi = {10.1137/0914028}
}

@book{saad2003iterative,
  title     = {Iterative methods for sparse linear systems},
  author    = {Saad, Yousef},
  year      = {2003},
  publisher = {SIAM},
  edition = {2},
  address = {Philadelphia},
  doi = {10.1137/1.9780898718003}
}

@book{smith1996domain,
  address   = {Cambridge},
  title     = {Domain Decomposition: Parallel Multilevel Methods for Elliptic Partial Differential Equations},
  author    = {Smith, Barry F. and Bj{\o}rstad, Petter E. and Gropp, William D.},
  year      = {1996},
  publisher = {Cambridge University Press}
}

@article{spillane2014abstract,
  title   = {Abstract Robust Coarse Spaces for Systems of {PDEs} via Generalized Eigenproblems in the Overlaps},
  author  = {Spillane, Nicole and Dolean, Victorita and Hauret, Patrice and Nataf, Fr{\'e}d{\'e}ric and Pechstein, Clemens and Scheichl, Robert},
  journal = {Numerische Mathematik},
  volume  = {126},
  number  = {4},
  pages   = {741--770},
  year    = {2014},
  doi     = {10.1007/s00211-013-0576-y}
}

@phdthesis{spillane2014robust,
  title  = {Robust Domain Decomposition Methods for Symmetric Positive Definite Problems},
  author = {Spillane, Nicole},
  school = {Universit{\'e} Pierre et Marie Curie},
  year   = {2014}
}

@article{strehlow2026fast,
  title     = {Fast-convergent two-level restricted additive {Schwarz} methods based on optimal local approximation spaces},
  author    = {Strehlow, Arne and Ma, Chupeng and Scheichl, Robert},
  journal   = {SIAM Journal on Scientific Computing},
  volume    = {48},
  number    = {1},
  pages     = {B31--B58},
  year      = {2026},
  publisher = {SIAM},
  doi = {10.1137/24M1695786}
}

@article{sun2026learning,
  title     = {Learning singularity-encoded {Green}'s functions with application to iterative methods},
  author    = {Sun, Qi and Li, Shengyan and Zheng, Bowen and Ju, Lili and Xu, Xuejun},
  journal   = {Journal of Computational Physics},
  pages     = {114894},
  year      = {2026},
  publisher = {Elsevier},
  volume = {559},
  doi = {10.1016/j.jcp.2026.114894}
}

@article{sussman1994levelset,
  author  = {Sussman, Mark and Smereka, Peter and Osher, Stanley},
  title   = {A Level Set Approach for Computing Solutions to Incompressible Two-Phase Flow},
  journal = {Journal of Computational Physics},
  volume  = {114},
  number  = {1},
  pages   = {146--159},
  year    = {1994},
  doi     = {10.1006/jcph.1994.1155}
}

@article{taghibakhshi2021optimization,
  author  = {Taghibakhshi, Ali and MacLachlan, Scott and Olson, Luke and West, Matthew},
  journal = {Advances in Neural Information Processing Systems},
  title   = {Optimization-based algebraic multigrid coarsening using reinforcement learning},
  year    = {2021},
  pages   = {12129--12140},
  volume  = {34},
  ranking = {rank5}
}

@article{taghibakhshi2022learning,
  title   = {Learning interface conditions in domain decomposition solvers},
  author  = {Taghibakhshi, Ali and Nytko, Nicolas and Zaman, Tareq Uz and MacLachlan, Scott and Olson, Luke and West, Matthew},
  journal = {Advances in Neural Information Processing Systems},
  volume  = {35},
  pages   = {7222--7235},
  year    = {2022}
}

@inproceedings{taghibakhshi2023mg,
  title        = {{MG-GNN}: Multigrid graph neural networks for learning multilevel domain decomposition methods},
  author       = {Taghibakhshi, Ali and Nytko, Nicolas and Zaman, Tareq Uz and MacLachlan, Scott and Olson, Luke and West, Matthew},
  booktitle    = {International conference on machine learning},
  pages        = {33381--33395},
  year         = {2023},
  organization = {PMLR}
}

@article{tancik2020fourier,
  title   = {Fourier features let networks learn high frequency functions in low dimensional domains},
  author  = {Tancik, Matthew and Srinivasan, Pratul and Mildenhall, Ben and Fridovich-Keil, Sara and Raghavan, Nithin and Singhal, Utkarsh and Ramamoorthi, Ravi and Barron, Jonathan and Ng, Ren},
  journal = {Advances in neural information processing systems},
  volume  = {33},
  pages   = {7537--7547},
  year    = {2020}
}

@book{toselli2005domain,
  address   = {Berlin, Heidelberg},
  title     = {Domain Decomposition Methods---Algorithms and Theory},
  author    = {Toselli, Andrea and Widlund, Olof},
  series    = {Springer Series in Computational Mathematics},
  volume    = {34},
  year      = {2005},
  publisher = {Springer}
}

@inproceedings{trifonov2025efficient,
  title     = {Efficient preconditioning for iterative methods with graph neural networks},
  author    = {Trifonov, Vladislav and Rudikov, Alexander and Iliev, Oleg and Laevsky, Yuri M and Oseledets, Ivan and Muravleva, Ekaterina},
  booktitle = {AI4X 2025 International Conference},
  year      = {2025},
  url = {https://openreview.net/forum?id=PjDJonBo9T}
}

@book{trottenberg2000multigrid,
  address   = {London},
  author    = {Trottenberg, Ulrich and Oosterlee, Cornelius W and Schuller, Anton},
  publisher = {Elsevier},
  title     = {Multigrid},
  year      = {2000}
}

@article{vasilyeva2026operator,
  title   = {Operator Learning on the Data-Driven Multiscale Space for Nonlinear Flow in Random Heterogeneous Porous Media},
  author  = {Vasilyeva, Maria and Pangilinan, Raphael},
  journal = {arXiv preprint arXiv:2606.25820},
  year    = {2026},
  url = {https://arxiv.org/abs/2606.25820}
}

@article{whitaker1996forchheimer,
  author  = {Whitaker, Stephen},
  title   = {The {Forchheimer} Equation: A Theoretical Development},
  journal = {Transport in Porous Media},
  volume  = {25},
  number  = {1},
  pages   = {27--61},
  year    = {1996},
  doi     = {10.1007/BF00141261}
}

@article{xu1992iterative,
  title     = {Iterative methods by space decomposition and subspace correction},
  author    = {Xu, Jinchao},
  journal   = {SIAM Review},
  volume    = {34},
  number    = {4},
  pages     = {581--613},
  year      = {1992},
  publisher = {SIAM},
  doi = {10.1137/1034116}
}

@article{zhang2024blending,
  title     = {Blending Neural Operators and Relaxation Methods in {PDE} Numerical Solvers},
  author    = {Zhang, Enrui and Kahana, Adar and Kopani{\v{c}}{\'a}kov{\'a}, Alena and Turkel, Eli and Ranade, Rishikesh and Pathak, Jay and Karniadakis, George Em},
  journal   = {Nature Machine Intelligence},
  volume    = {6},
  number    = {11},
  pages     = {1303--1313},
  year      = {2024},
  publisher = {Nature Publishing Group UK London},
  doi = {10.1038/s42256-024-00910-x}
}

@article{eckart1936approximation,
  title   = {The approximation of one matrix by another of lower rank},
  author  = {Eckart, Carl and Young, Gale},
  journal = {Psychometrika},
  volume  = {1},
  number  = {3},
  pages   = {211--218},
  year    = {1936},
  doi     = {10.1007/BF02288367}
}

@book{pinkus1985nwidths,
  title     = {{$n$}-Widths in Approximation Theory},
  author    = {Pinkus, Allan},
  series    = {Ergebnisse der Mathematik und ihrer Grenzgebiete (3)},
  volume    = {7},
  publisher = {Springer},
  address   = {Berlin},
  year      = {1985},
  doi = {10.1007/978-3-642-69894-1}
}

\end{document}